\documentclass[pdflatex,sn-mathphys-num]{sn-jnl}% Math and Physical Sciences Numbered Reference Style
\usepackage{graphicx}%
\usepackage{multirow}%
\usepackage{amsmath,amssymb,amsfonts}%
\usepackage{amsthm}%
\usepackage{mathrsfs}%
\usepackage[title]{appendix}%
\usepackage{xcolor}%
\usepackage{textcomp}%
\usepackage{manyfoot}%
\usepackage{booktabs}%
\usepackage{algorithm}%
\usepackage{algorithmicx}%
\usepackage{algpseudocode}%
\usepackage{listings}%

\numberwithin{equation}{section}

\theoremstyle{thmstyleone}
\newtheorem{theorem}{Theorem}[section]
\newtheorem{lemma}[theorem]{Lemma}

\newtheorem{corollary}[theorem]{Corollary}

\theoremstyle{thmstyletwo}
\newtheorem{definition}{Definition}[section]

\theoremstyle{thmstylethree}
\newtheorem{remark}{Remark}[section]
\newtheorem{example}{Example}[section]

\newcommand{\R}{\mathbb{R}}
\newcommand{\C}{\mathbb{C}}
\newcommand{\I}{\mathcal{I}}

\newcommand{\bx}{\boldsymbol{x}}
\newcommand{\by}{\boldsymbol{y}}
\newcommand{\br}{\boldsymbol{r}}
\newcommand{\btau}{\boldsymbol{\tau}}
\newcommand{\bn}{\boldsymbol{n}}
\newcommand{\bphi}{\boldsymbol{\phi}}
\newcommand{\calK}{\mathcal{K}}
\newcommand{\Sinc}{\operatorname{sinc}}

\begin{document}

\title[Corner Spectra and SGN for 2D Elastostatics]{A Singularity Guided Nystr\"om Method for Elastostatics on Two Dimensional Domains with Corners}

%%=============================================================%%
%% GivenName	-> \fnm{Joergen W.}
%% Particle	-> \spfx{van der} -> surname prefix
%% FamilyName	-> \sur{Ploeg}
%% Suffix	-> \sfx{IV}
%% \author*[1,2]{\fnm{Joergen W.} \spfx{van der} \sur{Ploeg} 
%%  \sfx{IV}}\email{iauthor@gmail.com}
%%=============================================================%%

\author[1]{\fnm{Baoling} \sur{Xie}}\email{baolingxie@zju.edu.cn}

\author*[1]{\fnm{Jun} \sur{Lai}}\email{laijun6@zju.edu.cn}

% \author[1]{\fnm{Baoling} \sur{Xie}}\email{baolingxie@zju.edu.cn}

% \author[1]{\fnm{Jun} \sur{Lai}}\email{laijun6@zju.edu.cn}

\affil[1]{\orgdiv{School of Mathematical Sciences},
          \orgname{Zhejiang University},
          \orgaddress{\city{Hangzhou}, \state{Zhejiang}, \postcode{310027}, \country{China}}}
          
%%==================================%%
%% Sample for unstructured abstract %%
%%==================================%%

\abstract{We develop a comprehensive analytical and numerical framework for boundary integral equations (BIEs) of the 2D Lam\'e system on cornered domains. By applying local Mellin analysis on a wedge, we obtain a factorizable characteristic equation for the singular exponents of the boundary densities, and clarify their dependence on boundary conditions. The Fredholm well-posedness of the BIEs on cornered domains is proved in weighted Sobolev spaces. We further construct an explicit density-to-Taylor mapping for the BIE and show its invertibility for all but a countable set of angles. Based on these analytical results, we propose a singularity guided Nystr\"om (SGN) scheme for the numerical solution of BIEs on cornered domains. The SGN uses the computed corner exponents and a Legendre-tail indicator to drive panel refinement. An error analysis that combines this refinement strategy with an exponentially accurate far-field quadrature rule is provided. Numerical experiments across various cornered geometries demonstrate that SGN obtains higher order accuracy than uniform Nystr\"om method and reveal a crowding-limited regime for domains with re-entrant angles.}

\keywords{Boundary integral equations; Lam\'e system; corner singularities; Mellin symbol; Nystr\"om methods; Fredholm theory}

%%\pacs[JEL Classification]{D8, H51}
%%\pacs[MSC Classification]{65R20, 74B05, 45E05}
\pacs[MSC Classification]{65R20, 65N38, 74B05, 45E05, 65N80}

\maketitle
% --------------------------
% Section 1: Introduction
% --------------------------
\section{Introduction}\label{sec1}

Geometric discontinuities such as corners, edges, and junctions are fundamental features of models in continuum mechanics. They occur across a wide range of scales, from crack tips and material interfaces to tectonic faults, boundaries between geological strata, and subsurface scatterers in geophysics~\cite{Christensen1979composites,He1989interface,Leguillon2002delamination,epstein2016smoothed}. In elastostatics and related elliptic problems, these nonsmooth geometries generate stress fields with corner singularities that may be oscillatory and strongly localized. The resulting behavior influences scattering patterns, delamination, and failure. Consequently, developing a high order numerical scheme is essential for analyzing these singular structures.

Traditional numerical techniques, particularly the finite element method, are widely used in structural analysis, but they struggle to resolve the branch-type singular behavior near corners. Asymptotics at the material junctions often involve fractional and occasionally complex corner exponents, which demand excessive local refinement and can produce unstable artifacts~\cite{Babuska1994pFEM}. Beyond volume-based discretizations, there exists the boundary-oriented class of approaches that represent or approximate the solution from boundary data. Popular schemes include the method of fundamental solutions~\cite{Alves2024} and the Lightning solver~\cite{gopal2019solving, Gopal2019, Brubeck2022}, which place sources or poles that cluster exponentially near corners and capture the expected local behavior. Although these methods are meshless and can achieve high-order accuracy in certain geometries, their underlying theoretical foundations are still being established. 

In this paper, we focus on boundary integral equation methods (BIEs) formulated via layer potential theory in elastostatics~\cite{Linkov2002, hsiao2008boundary, kanwal2013linear}. BIEs not only reduce dimensionality but also naturally handle the far-field conditions~\cite{atkinson1997numerical}. However, while BIEs are highly accurate on smooth geometries~\cite{Lai2019,Dong2021,Dong2022}, their application to nonsmooth domains with cracks, interfaces, and sharp re-entrant corners is complicated by two issues. First, the kernel function for the second-kind formulation is a tensor and involves Cauchy singular integrals, which complicate the numerical analysis. Second, one must account for branch singularities in the unknown boundary density, which can be oscillatory and strongly coupled near corners~\cite{costabel2010corner, blaasten2014corners}.

To address these issues, various numerical strategies have been proposed. Approaches include kernel splitting techniques combined with generalized Gaussian quadrature~\cite{ma1996generalized, helsing2009integral, yao2024robust}, which have proven to be effective for handling logarithmic and Cauchy singularities. Another prominent approach for cornered domains involves local mesh refinement, often combined with compression techniques to control the size of the linear system, including the $L^2$ weighting method~\cite{bremer2012nystrom} and recursive compressed inverse preconditioning~\cite{helsing2008corner}. Alternatively, geometric regularization via corner rounding to modify the domain so that  the most severe singular features are suppressed has also been proposed in~\cite{epstein2016smoothed}. Recently, Rokhlin and Serkh constructed an effective quadrature method by utilizing explicit representations of the solution's asymptotic behavior near corners~\cite{serkh2016solution, serkh2016solution2}. They observed that solutions to these integral equations can be represented by rapidly convergent singular power series for Laplace problems, fractional order Bessel functions for Helmholtz cases, and expressions of the form $\sum_j\left(c_j t^{\mu_j} \sin \left(\beta_j \ln (t)\right)+d_j t^{\mu_j} \cos \left(\beta_j \ln (t)\right)\right)$ for Stokes flow, as summarized in Table~\ref{tab:2d_expansions}. Numerically, explicit singular basis methods require the prior determination of the characteristic exponents. To ensure stability and accuracy, the resulting bases are orthogonalized by SVD, yielding a well-conditioned set of functions. These expansions are then discretized using generalized Gaussian quadrature rules tailored to the singular behavior. 

\begin{table}[!ht]
\centering
\setlength{\tabcolsep}{5pt}
\caption{Corner singularity expansions for 2D elliptic PDEs. Here \texorpdfstring{$t$}{t} is the distance to the corner, \texorpdfstring{$\pi\alpha$}{pi·alpha} is the opening angle, and \texorpdfstring{$\mu_j,\beta_j$}{mu\_j, beta\_j} are characteristic exponents.}
\label{tab:2d_expansions}
\begin{tabular}{@{}llp{6.5cm}@{}}
\toprule
\textbf{Equation} & \textbf{Green's Function} & \textbf{Density Expansion} \\
\midrule
Laplace &
$-\dfrac{1}{2\pi}\ln r$ &
\raggedright $\displaystyle\sum_j \Bigl(c_j t^{\frac{j}{\alpha}} + d_j t^{\frac{j}{2-\alpha}}\Bigr)$\par \\[1pt]
\addlinespace[3pt]
Helmholtz &
$\dfrac{i}{4}H_0^{(1)}(kr)$ &
\raggedright $\displaystyle\sum_j\Bigl(c_j J_{\frac{j}{\alpha}}(kt)+d_j J_{\frac{j}{2-\alpha}}(kt)\Bigr)$\par \\[1pt]
\addlinespace[3pt]
Stokes &
$\dfrac{1}{4\pi}\!\left(-\ln r\,\mathbf{I} + \dfrac{\mathbf{r}\otimes\mathbf{r}}{r^2}\right)$ &
\raggedright $\displaystyle\sum_j\!\Bigl(c_j t^{\mu_j}\sin(\beta_j\ln t) + d_j t^{\mu_j}\cos(\beta_j\ln t)\Bigr)$\par \\[1pt]
\bottomrule
\end{tabular}
\end{table}

Although the singularity analysis proposed in~\cite{serkh2016solution} has been successfully extended to vectorial systems such as the 2D Stokes equations~\cite{rachh2020solution}, its application to elastic equations remains underdeveloped. In particular, the Cauchy-type singularities of the elastic Green's kernel prevent the parity transform from decoupling the system. A fully coupled $4\times4$ boundary formulation makes both the analysis of the corner exponents and the design of efficient algorithms substantially more involved. This paper fills this gap by establishing a complete analytical and numerical framework for corner singularities in two dimensional elastostatics. Our contributions are threefold:
\begin{enumerate}
    \item \textbf{Corner exponents from a factorizable Mellin symbol:} 
    On a wedge, we derive a characteristic equation for the boundary density's singular exponents that factorizes into two branches, identifying real and complex exponents and explaining which modes are compatible with interior vs.\ exterior conditions. This connects the BIE symbol calculus directly to classical PDE corner analysis.
    \item \textbf{Fredholm criterion and density-to-Taylor mapping:}
    In the weighted Sobolev space $H^s_\nu(\Gamma)^2$, we obtain a computable Fredholm criterion $s-\nu\notin\{\Re z_{n,j}\}\cup\mathbb Z$ for the interior double-layer operator, with $z_{n,j}$ denoting the zeros of the Mellin symbol. We also give an explicit density-to-Taylor map $\mathbf B(\theta)$ and prove its invertibility for all but a discrete set of angles, thereby linking corner singular coefficients to smooth boundary data.
    \item \textbf{A high order singularity guided Nystr\"om (SGN) scheme:}
    Rather than spanning the density with singular powers, whose hybrid variant is severely ill-conditioned, we use the computed exponents to drive panel adaptivity via a multi-exponent Legendre-tail indicator. Moreover, we prove an a priori error bound that couples this indicator with exponentially small far-field quadrature error. In numerical experiments, we observe that the SGN scheme achieves higher accuracy than a uniform Nystr\"om baseline. Moreover, on re-entrant corners, the attainable accuracy is crowding-limited by near singular integration.
\end{enumerate}
    
The remainder of the paper is organized as follows. Section~\ref{sec2} formulates the elastostatic BIEs and develops the well-posedness analysis using the Fredholm index theory in weighted Sobolev spaces. Section~\ref{sec3} presents the mathematical apparatus, including Mellin-type integrals and a characterization of the root distribution for $\Sinc$ functions. Section~\ref{sec4} performs the local corner analysis, derives the characteristic equation, and explicitly constructs the density-to-Taylor mapping. Section~\ref{sec5} presents the singularity guided Nystr\"om scheme, error analysis, and several numerical experiments. Section~\ref{sec6} summarizes the paper and discusses future directions. Technical details are deferred to Appendices~\ref{app:AF_calculation} and~\ref{app:invertibility}.

% -------------------------------
% Section 2: Problem formulation
% -------------------------------

\section{Problem formulation}\label{sec2}

\begin{figure}[!ht]
    \centering
    \includegraphics[width=.4\linewidth]{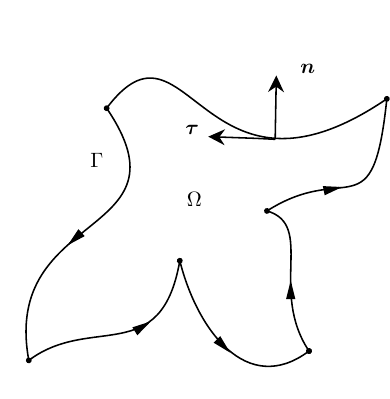} 
    \caption{A cornered domain $\Omega$ with boundary $\Gamma$.}
    \label{fig:sampledomain}
\end{figure}

Consider a bounded, simply connected planar elastic domain $\Omega\subset \R^2$ with a piecewise smooth boundary $\Gamma$. The boundary is assumed to contain a finite number of corners, each locally diffeomorphic to a wedge with an opening angle (cusps excluded). Let $\bn=(n_1, n_2)^\top$ and $\btau=(\tau_1, \tau_2)^\top$ denote the unit outward normal and the positively oriented unit tangent on $\Gamma$ respectively, where $\tau_1 = -n_2$ and $\tau_2 = n_1$, as shown in Figure~\ref{fig:sampledomain}. 
% \begin{figure}[!ht]
% \centering
% \scalebox{0.8}{\input{figs/SampleDomain.tikz}}
% \caption{A sample domain with six corners.}
% \label{fig:sampledomain}
% \end{figure}

Assume the domain $\Omega$ is filled with a homogeneous, isotropic elastic material. The displacement field $\boldsymbol{u}=(u_1,u_2)^\top$ under static loading satisfies the homogeneous Lam\'e system:
\begin{equation}\label{eq:lame_equation}
    \mathcal{L} \boldsymbol{u} \triangleq \mu \Delta \boldsymbol{u} + (\lambda + \mu) \nabla (\nabla \cdot \boldsymbol{u}) = \boldsymbol{0}, \quad \bx \in \Omega,
\end{equation}
where $\lambda$ and $\mu$ are the Lam\'e parameters. The conditions $\mu > 0$ and $\lambda + \mu > 0$ ensure the strong ellipticity of $\mathcal{L}$~\cite{costabel2010corner}. Throughout, we assume that all material parameters are constant and focus on static elasticity, although the methodology extends naturally to quasi-static or time-harmonic elastic wave scattering~\cite{yao2024robust}.

Two standard types of boundary conditions are considered. Under Dirichlet boundary conditions, the displacement field satisfies
\begin{equation}\label{eq:HC_bc}
    \boldsymbol{u}(\bx) = \boldsymbol{f}(\bx), \quad \bx \in \Gamma
\end{equation}
with $\boldsymbol{f} \in H^{1/2}(\Gamma)^2$. Under Neumann boundary conditions, the traction field satisfies
\begin{equation}\label{eq:SF_bc}
    \mathcal{T} \boldsymbol{u}(\bx) = \boldsymbol{g}(\bx), \quad \bx \in \Gamma
\end{equation}
with $\boldsymbol{g} \in H^{-1/2}(\Gamma)^2$. The traction operator $\mathcal{T}: H^{1}(\Omega)^2 \to H^{-1/2}(\Gamma)^2$ is defined by $\mathcal{T} \boldsymbol{u} \triangleq \boldsymbol{\sigma}(\boldsymbol{u}) \bn$, where the stress tensor $\boldsymbol{\sigma}(\boldsymbol{u})$ satisfies
\begin{equation}\label{eq:stress_tensor}
    \sigma_{ij}(\boldsymbol{u}) = \lambda(\nabla \cdot \boldsymbol{u})\delta_{ij} + \mu\left( \frac{\partial u_i}{\partial x_j} + \frac{\partial u_j}{\partial x_i} \right), \quad i,j=1,2.
\end{equation}
Notably, the analytic framework admits natural extensions to various boundary conditions, which are omitted here for clarity. The central difficulty in such problems lies in the singular behavior of the stress field near corners, which is essential for predicting failure and ensuring structural integrity.

\subsection{Boundary integral formulation}
%% BIE formulations: elasto-static equations
\label{subsec21}
To efficiently resolve the singularity near corners, we employ a boundary integral formulation based on layer potential theory. Let $\bphi:\Gamma \to \R^2$ be a vector-valued density function defined on the boundary. The fundamental solution (Kelvin matrix~\cite{McLean2000}) of the 2D Lam\'e system in plane strain is given by
\begin{equation}\label{eq:fundamental_solution}
    \mathbf{G}(\bx, \by) = \frac{1}{4\pi\mu} \left[ -\frac{\lambda + 3\mu}{\lambda + 2\mu} \ln r \cdot \mathbf{I} + \frac{\lambda + \mu}{\lambda + 2\mu} \frac{\br \otimes \br}{r^2} \right],
\end{equation}
where $\br = \bx - \by$ and $r = \|\br\|$. This kernel defines the single-layer potential (SLP):
\begin{equation}\label{eq:SLP}
    S[\bphi](\bx) = \int_\Gamma \mathbf{G}(\bx, \by) \bphi(\by) \, \mathrm{d}s_{\by}.
\end{equation}
The corresponding traction fundamental solution $\mathbf{D}(\bx,\by)$, derived via the traction operator $\mathcal{T}_{\by}$ acting on $\mathbf{G}$, is 
\begin{equation}\label{eq:fundamental_solution_traction}
\begin{aligned}
\mathbf{D}(\bx,\by)&=\frac{\mu}{\lambda+2\mu}\left[\frac{\bn(\by)\cdot\br}{2\pi r^2}\mathbf{I}+\frac{\btau(\by)\cdot\br}{2\pi r^2}\mathbf{L}\right]
+\frac{\lambda+\mu}{\lambda+2\mu}\frac{(\bn(\by)\cdot\br)\br\otimes\br}{\pi r^4},
\end{aligned}
\end{equation}
where $\mathbf{L}\triangleq\begin{bmatrix}
		           0 & 1 \\
		          -1 & 0
	              \end{bmatrix}$.
The double-layer potential (DLP) is subsequently defined by
\begin{equation}\label{eq:DLP}
    D[\bphi](\bx)=\int_\Gamma\mathbf{D}(\bx,\by)\bphi(\by)\mathrm{d}s_{\by}.
\end{equation}

The single and double-layer potentials satisfy Lam\'e equations in $\Omega$ and $\R^2 \setminus \overline{\Omega}$, and their limiting values on the boundary exhibit the following well-known jump behaviors.

\begin{lemma}[Jump relations \cite{atkinson1997numerical, kanwal2013linear}]\label{lem2.1}
    When $\bx$ approaches the boundary $\Gamma$, the double layer potential $D[\bphi](\bx)$ and the traction of single-layer potential $\mathcal{T}_{\bx}S[\bphi]$ exhibit the following jumps:
    \begin{equation}\label{eq:jumpBIE}
        \begin{aligned}
            \lim _{\epsilon \rightarrow 0^{+}} D[\bphi](\bx \pm \epsilon \bn) & = \pm \frac{\theta(\bx)}{2\pi} \bphi(\bx)+\calK[\bphi](\bx), \\
            \lim _{\epsilon \rightarrow 0^{+}} \mathcal{T}_{\bx} S[\bphi](\bx \pm \epsilon \bn) & =\mp \frac{\theta(\bx)}{2\pi} \bphi(\bx)+\calK^{*}[\bphi](\bx),
        \end{aligned}
    \end{equation}
    where $\theta(\bx)$ denotes the interior angle $\theta_{int}(\bx)$ at $\bx \in \Gamma$ for interior problems, and the exterior angle $(2\pi-\theta_{int}(\bx))$ for exterior problems. Here, $\calK$ is the corresponding boundary integral operator of $D$ and $\calK^{*}$ is the adjoint of $\calK$. Both of them are understood in the Cauchy principal value sense.
\end{lemma}

Based on these potentials, we formulate the second-kind boundary integral equations. Throughout, we focus on the interior Dirichlet problem (IDP) for the Lam\'e system. Using the double-layer potential representation $\boldsymbol{u}(\bx)=D[\bphi](\bx)$ and the jump relation from Lemma~\ref{lem2.1}, we obtain
\begin{equation}~\label{eq:bie_original}
    \left(-\frac{\theta(\bx)}{2\pi}\I + \calK\right)\bphi(\bx) = \boldsymbol{f}(\bx), \quad \bx \in \Gamma.
\end{equation}
\begin{remark}
    Although our analysis below focuses on interior Dirichlet problems, the same second-kind formulation extends to exterior elastostatic problems, provided one imposes suitable decay conditions at infinity to ensure uniqueness~\cite{hsiao2008boundary}. Similar extensions hold for the Neumann case via the single-layer potential and $\calK^\ast$.
\end{remark}

\subsection{Well-posedness analysis}
\label{subsec22}
The well-posedness of the boundary integral equation~\eqref{eq:bie_original} is determined by the spectral properties of the associated double-layer potential operator $\calK$, also known as the Neumann-Poincar\'e (NP) operator, within an appropriate function space. For domains with corners, the analysis is non-trivial and requires specialized techniques beyond classical theories. Taking the Laplace equation as an example, the operator $\calK$ is compact in $L^2(\Gamma)$ on a $C^{1,\alpha}$ boundary, where the Fredholm alternative readily ensures the existence and uniqueness of solutions. However, when the boundary $\Gamma$ includes corners, the NP operator exhibits a continuous spectrum, and classical Fredholm theory requires local regularization techniques~\cite{atkinson1997numerical,kanwal2013linear}. The situation becomes even more subtle in the Lam\'e system, since even on smooth boundaries, the NP operator $\calK$ is non-compact and admits the decomposition:
\begin{equation}\label{eq:decompNP_lame}
    \calK[\bphi](\bx) = \int_{\Gamma} \left(\mathcal{T}_{\by} \mathbf{G}(\bx,\by)\right)^{\top} 
    \bphi(\by) \mathrm{d} s_{\by} 
    :=\calK_{0}[\bphi](\bx)+\calK_{1}[\bphi](\bx),
    \quad \bx \in \Gamma,
\end{equation}
where the kernels of $\calK_{0}$ and $\calK_{1}$ are given by
\begin{align}
    K_{0}(\bx,\by) &= \frac{\mu}{\lambda+2 \mu} \frac{\btau(\by) \cdot \br}{2 \pi r^{2}} \mathbf{L}, \\
    K_{1}(\bx,\by) &= \frac{\bn(\by) \cdot \br}{2 \pi r^{2}} \mathbf{I} -\frac{\lambda+\mu}{\lambda+2 \mu}\left(-\frac{(\bn(\by) \cdot \br) \br \otimes \br}{\pi r^{4}}+\frac{\bn(\by) \cdot \br}{2 \pi r^{2}} \mathbf{I}\right). 
\end{align}
Here, $K_1$ is only weakly singular, whereas $K_0$ contains the Cauchy-type singularity responsible for non-compactness. More precisely, using Stokes' theorem and the Hilbert transform~\cite{Ando2018Spectral}, it is known that $\calK_0^2-k_0^2 \I$ is compact in $H^{1/2}(\Gamma)^2$ with $k_0=\frac{\mu}{2(\lambda+2\mu)}$. Although specific regularization techniques enable well-posedness for smooth domains, they become insufficient when the boundary contains corners.

To address this challenge, we adopt a framework developed for elliptic boundary value problems in nonsmooth domains. This theory, initiated by Kondratiev~\cite{Kondratiev1967Boundary} and significantly extended by Costabel, Dauge, and others~\cite{Dauge1988Elliptic, CostabelDauge2002Weighted}, provides a general operator-theoretic setting for boundary integral formulations on domains with corners. The key component of this framework involves weighted Sobolev spaces, specifically designed to capture the singular behavior of solutions near corners.
\begin{definition}[Weighted Sobolev space]\label{def:weighted_sobolev}
Let $V = \{v_1, \ldots, v_M\}$ be the set of corner vertices on the boundary $\Gamma$. Let $r_k(\bx) = |\bx - v_k|$ denote the Euclidean distance from a point $\bx \in \Gamma$ to the corner $v_k$. For any real smoothness order $s\ge0$ and weight exponent $\nu$, the weighted Sobolev space $H_\nu^s(\Gamma)$ is defined as the space of distributions $\phi$ on $\Gamma$ for which the following norm is finite:
    $$\|\phi\|_{H_\nu^s(\Gamma)}^2 = \|\phi\|_{H^s(\Gamma \setminus \cup U_k)}^2 + \sum_{k=1}^{M} \int_{U_k \cap \Gamma} \sum_{|\gamma| \le \lfloor s \rfloor} r_k(\bx)^{2(\nu - s + |\gamma|)} |D^\gamma \phi(\bx)|^2 \mathrm{d}s_{\bx}.$$
    Here, $U_k$ is a small neighborhood of the corner $v_k$, $D^\gamma$ denotes the $|\gamma|$-th order tangential derivative along $\Gamma$, and the definition is extended to noninteger $s$ via interpolation. The space for 2D vector-valued functions is denoted by $H_\nu^s(\Gamma)^2$. This norm incorporates the weight function $r_k(\bx)$ to characterize the local behavior of functions near the corner $v_k$. The parameter $\nu$ precisely controls the admissible singularity, where larger positive values of $\nu$ permit stronger singularities.
\end{definition}

Equipped with this functional framework, we can now state the main result of this section, which connects the well-posedness of the operator to the singular exponents derived in Section~\ref{sec4}.
%% 定理给出了算子成为 Fredholm 算子的充要条件，即权重参数 $\nu$ 必须避开由奇性指数决定的临界值。
\begin{theorem}\label{thm:fredholm}
    Let $\{z_{n,j}\}$ be the set of singular exponents defined as the roots of the transcendental equation~\eqref{eq:trans1}, as established in Theorem~\ref{thm4.1}. Fix $s \in (1/2, 3/2)$ and $\nu \in \R$. If $s - \nu \notin \{\operatorname{Re} z_{n,j}\} \cup\mathbb Z$, then the boundary integral operator
    \begin{equation}\label{eq:NP_operator}
        \mathcal{A} := \left(-\frac{\theta(\bx)}{2\pi}\I + \calK\right) : H_\nu^s(\Gamma)^2 \to H_\nu^s(\Gamma)^2
    \end{equation}
    is a Fredholm operator of index zero. 
\begin{proof}
    % 1. 局部化算子：采用单位分解，将原始算子分解为两类。
    % 这里的$\chi_\ell$是作用在y上的单位分解，$\psi_\ell$是作用在x上的截断算子。并且远离对角块部分的核函数没有奇异性，是个无限光滑算子。
    \emph{Step 1 (Localization).}
    Let $\{\chi_\ell\}$ be a smooth partition of unity on $\Gamma$ subordinate to an open cover by smooth arcs and corner patches, and set $\psi_\ell:=\chi_\ell$. For each $\ell$, define the set of neighboring indices
    $N(\ell):=\bigl\{\,m:\ \operatorname{supp}\psi_\ell\cap \operatorname{supp}\chi_m\neq\emptyset\,\bigr\}$.
    Define the localized near-diagonal piece and the far remainder by
    $$
    \mathcal A_\ell \;:=\; \psi_\ell\,\mathcal A \Bigl(\sum_{m\in N(\ell)}\chi_m\Bigr),
    \qquad
    \mathcal R \;:=\! \sum_{\substack{\ell,m:\ \operatorname{dist}(\operatorname{supp}\psi_\ell,\,\operatorname{supp}\chi_m)>0}}
    \psi_\ell\,\mathcal A\,\chi_m \,.
    $$
    Using $I=\sum_\ell\psi_\ell=\sum_m\chi_m$, we obtain the exact decomposition
    $$
    \mathcal A \;=\; \sum_\ell \mathcal{A}_\ell \;+\; \mathcal{R}.
    $$
    For every pair $(\ell,m)$ with $\operatorname{dist}(\operatorname{supp}\psi_\ell,\operatorname{supp}\chi_m)>0$, the kernel of $\psi_\ell\,\mathcal A\,\chi_m$ is smooth, hence each term in $\mathcal R$ lies in $\Psi^{-\infty}$\cite{atkinson1997numerical}. Therefore, by Rellich--Kondrachov~\cite{McLean2000} near vertices and standard embeddings on smooth arcs, $\mathcal R$ is compact on $H_\nu^s(\Gamma)^2$ with $s \in (1/2, 3/2)$ and $\nu \in \R$. Thus it suffices to analyze the localized operators $\mathcal A_\ell$, each supported either on a smooth arc or on a single corner patch.

    % 2. 光滑部分：0阶椭圆算子。
    % 证明了算子-1/2I+K的主符号可逆，即0阶椭圆算子，这就保证了存在一个零阶的微分算子P，使得P(-1/2I+K) = I-compact。也就是说在某个局部区域中，这个算子可以被一个似逆“反演”，前后相乘只相差一个紧算子。
    \medskip\noindent
    \emph{Step 2 (Smooth arcs).}
    On a smooth arc, it holds $\theta(\bx)=\pi$ and hence $\mathcal A_\ell = -\tfrac12 \I + \calK$. With $\calK=\mathcal K_0+\mathcal K_1$ in Eq.~\eqref{eq:decompNP_lame}, one has $\mathcal K_1\in\Psi^{-1}$. In arclength coordinates, the kernel of $\mathcal K_0$ is
    $$
    K_0(s,t)\;=\;\mathrm{p.v.}\,\frac{k_0}{\pi(s-t)}\,\mathbf L \;+\; \text{smoother terms},
    $$
    where $\mathbf L=\begin{bmatrix}0&1\\ -1&0\end{bmatrix}$ and $k_0=\dfrac{\mu}{2(\lambda+2\mu)}\in(0,1/2)$. Therefore, by the pseudodifferential calculus~\cite{Ando2018Spectral},
    $$
    \mathcal K_0 \;=\; k_0\,\mathcal H\,\mathbf L \;+\; \Psi^{-1},
    $$
    where $\mathcal H$ is the 1D Hilbert transform with symbol $-i\,\mathrm{sgn}\,\xi$ and $\mathbf L$ has eigenvalues $\pm i$. Thus the principal symbol of $-\tfrac12 I+\calK$ has eigenvalues $-\frac12\pm k_0\neq 0$ for $\xi\neq0$. Hence, $\mathcal A_\ell\in\Psi^{0}$ is elliptic and admits a microlocal parametrix modulo compact operators on $H^s_\nu$.
  
    % 3. 角点部分：Mellin象征。
    % 主要推导参见超越方程部分～
    \medskip\noindent
    \emph{Step 3 (Corner patches).}
    For the corner patches, the corner degenerate calculus~\cite{Dauge1988Elliptic} reduces $\mathcal A_\ell$ to an operator family via the Mellin transform $\mathfrak M$ in the radial variable:
    $$
    \mathfrak M\mathcal A_\ell \;\equiv\; \mathbf A(z,\theta),
    $$
    as detailed in Section~\ref{subsec41}. Then $\mathcal A_\ell:H^s_\nu\to H^s_\nu$ is Fredholm if and only if $\mathbf A(z,\theta)$ is invertible along the vertical line $\Re z \;=\; s-\nu$. Here $\mathbf A(z,\theta)$ is meromorphic in $z$, with simple poles at $z\in\mathbb Z$ coming from $\cot(\pi z)$ and $\csc(\pi z)$, and zeros at the singular exponents $z_{n,j}$ characterized by $\det\mathbf A(z,\theta)=0$. Hence, under the exclusions
    $$
    s-\nu \;\notin\; \{\operatorname{Re} z_{n,j}\} \cup \mathbb Z
    $$
    the matrix $\mathbf A(z,\theta)$ is invertible for all $z$ with $\Re z=s-\nu$, which produces a corner parametrix modulo compact operators on $H^s_\nu$.

    % 4. 算子的 Fredholmness
    % 核心在于利用第2/3步的结论把局部近似逆粘起来
    \medskip\noindent
    \emph{Step 4 (Gluing parametrices).}
    For each $\ell$, set $R_\ell:=\psi_\ell\mathcal A-\mathcal A_\ell
    =\psi_\ell\mathcal A\bigl(I-\sum_{m\in N(\ell)}\chi_m\bigr)$. Each summand in $R_\ell$ has separated supports, hence $R_\ell\in\Psi^{-\infty}$ and is compact on $H^s_\nu(\Gamma)^2$. By Steps~2–3 there exists a patchwise parametrix $P_\ell$ with $P_\ell\mathcal A_\ell=\psi_\ell-K_\ell$, where $K_\ell$ is compact. Define the global parametrix by
    $\mathcal P:=\sum_\ell P_\ell\,\psi_\ell$, then
    $$
    \mathcal P\mathcal A
    =\sum_\ell P_\ell\psi_\ell\mathcal A
    =\sum_\ell P_\ell(\mathcal A_\ell+R_\ell)
    =\sum_\ell(\psi_\ell-K_\ell)+\sum_\ell P_\ell R_\ell
    =I-\mathcal C_1,
    $$
    with $\mathcal C_1$ compact on $H^s_\nu(\Gamma)^2$. The right parametrix is analogous, hence $\mathcal A$ is Fredholm by Atkinson's theorem~\cite{McLean2000}.

    % 5. 根据同伦证明Fredholm index为 0
    \medskip\noindent
    \emph{Step 5 (Index $0$).}
    For $t\in[0,1]$, continuously open each corner angle by $\theta(t):=(1-t)\theta+t\pi$, and denote the resulting boundary by $\Gamma_t$. Choose a smooth family of boundary diffeomorphisms $\Phi_t:\Gamma\to\Gamma_t$ and use pullback to identify the weighted spaces on $\Gamma_t$ with the fixed space $H^s_\nu(\Gamma)^2$. Define $\widetilde{\mathcal A}_t := \Phi_t^*\,\mathcal A_t\,(\Phi_t^{-1})^* : H^s_\nu(\Gamma)^2 \to H^s_\nu(\Gamma)^2$, then $t\mapsto \widetilde{\mathcal A}_t$ is norm-continuous~\cite{Kato1995}. For the interior Dirichlet problem of the Lam\'e system, one can fix $\nu\in\R$ such that the weight line $\Re z=s-\nu$ stays at a uniform positive distance from $\{\Re z_{n,j}(\theta(t))\}\cup\mathbb Z$ for all $t\in[0,1]$. Hence, by Step~3, $\widetilde{\mathcal A}_t$ is Fredholm for all $t$ and $\operatorname{ind}(\mathcal A)=\operatorname{ind}(\widetilde{\mathcal A}_0)=\operatorname{ind}(\widetilde{\mathcal A}_1)$. By Step~2, the principal symbol of $\widetilde{\mathcal A}_1$ has eigenvalues $-\tfrac12\pm k_0\neq 0$, so $\det\sigma(\widetilde{\mathcal A}_1)$ is a positive constant and thus has zero winding number. The index formula~\cite{Boettcher2006} yields $\operatorname{ind}(\widetilde{\mathcal A}_1)=0$ and therefore $\operatorname{ind}(\mathcal A)=0$.
\end{proof}
\end{theorem}
\begin{remark}
    According to Kondratiev's theory~\cite{Kondratiev1967Boundary}, the lower bound $s>1/2$ yields the continuous embedding $H^s(\Gamma)\hookrightarrow C^0(\Gamma)$, which is a convenient setting for the localization of corner singularities and the use of Mellin-type transforms. The upper bound $s<3/2$ restricts the weighted norms to involve at most first-order tangential derivatives near the vertices. 
    %Moreover, for smooth boundaries, we can obtain $H_\nu^s(\Gamma)=H^s(\Gamma)\hookrightarrow H^{1/2}(\Gamma)$ for all $s>1/2$. When corners exist, through an appropriate choice of the weight parameter $\nu$, this weighted space is compatible with the classical energy space $H^{1/2}(\Gamma)$.
\end{remark}

We now arrive at the well-posedness of the boundary integral equation for the interior Dirichlet problem.
\begin{theorem}\label{thm:wellposedness}
    Let the operator $\mathcal{A}$ be defined as in Theorem~\ref{thm:fredholm}. If the condition $s - \nu \notin \{\operatorname{Re} z_{n,j}\}\cup\mathbb Z$ holds, then
    $$
    \mathcal{A}: H_\nu^s(\Gamma)^2 \to H_\nu^s(\Gamma)^2
    $$
    is an isomorphism. Consequently, the corresponding boundary integral equation for the interior Dirichlet problem is well-posed.
\begin{proof}
    Let $\bphi \in \ker(\mathcal{A})$, then the corresponding double-layer potential $\boldsymbol{u}(\bx) = D[\bphi]$ satisfies the boundary condition
    $$
    \boldsymbol{u}|_{\Gamma^{-}}(\bx) = \left(-\frac{\theta(\bx)}{2\pi}\I + \calK\right)[\bphi](\bx) = \mathbf{0}.
    $$
    By the uniqueness theorem for the interior Dirichlet problem~\cite{McLean2000}, $\boldsymbol{u}\equiv\mathbf{0}$ in $\Omega$. For the Lam\'e double layer, the traction is continuous across $\Gamma$, hence $\mathcal T\boldsymbol u|_{\Gamma^+}=\mathcal T\boldsymbol u|_{\Gamma^-}=\mathbf 0$. Observing that the double-layer potential decays at infinity, the exterior uniqueness theorem implies $\boldsymbol u\equiv\mathbf 0$ in the exterior domain. Using the displacement jump relation $\boldsymbol u|_{\Gamma^+}-\boldsymbol u|_{\Gamma^-}=\boldsymbol\phi$ yields $\boldsymbol\phi=\mathbf 0$, which implies $\ker(\mathcal A)=\{\mathbf 0\}$. Since $\mathcal A$ is Fredholm of index $0$, injectivity implies bijectivity with a bounded inverse, which completes the proof.
\end{proof}
\end{theorem}
\begin{remark}
    For the interior Neumann problem, where $\mathcal{T}\boldsymbol{u} = \boldsymbol{g}$ on $\Gamma$, a natural representation is the single-layer potential $\boldsymbol{u} = S[\boldsymbol{\psi}]$. This leads to the adjoint equation $(\frac{\theta(\bx)}{2\pi}\I + \calK^*) [\boldsymbol{\psi}] = \boldsymbol{g}$. The well-posedness of this formulation follows a slightly different path, as the kernel of the operator $\mathcal{A}^*$ is non-trivial and corresponds to the space of densities that generate rigid body motions. In 2D, this space is 3-dimensional, spanned by two translations and one rotation~\cite{hsiao2008boundary}.
\end{remark}

The sequence of results in this subsection provides a theoretical framework for the analysis of the Lam\'e system on 2D domains with corners. It confirms that by choosing an appropriate weighted Sobolev space which avoids the corner exponents and the integers, the boundary integral formulation of the interior Dirichlet problem is well-posed, guaranteeing a unique solution.

% ----------------------------------
% Section 3: Mathematical Apparatus
% ----------------------------------
\section{Mathematical Preliminaries}\label{sec3}
%% Section 3 所有数学准备，如下所示：
% 1. 应用复变Cauchy integral theorem，计算需要的三类积分形式；
% 2. 求解超越方程Sinc(w)+Sinc(\alpha)=0在复数域中的所有根 (with figures)；
\subsection{Three types of singular integrals}
\label{subsec31}
In this subsection, we compute three Mellin-type singular integrals, referred to as the Cauchy, Laplace, and Stokes cases, respectively. These integrals play a crucial role in the calculations presented in Section~\ref{sec4}. 

\begin{figure}[!ht]
\centering
  \begin{minipage}[b]{.47\linewidth}
    \centering\includegraphics[width=\linewidth]{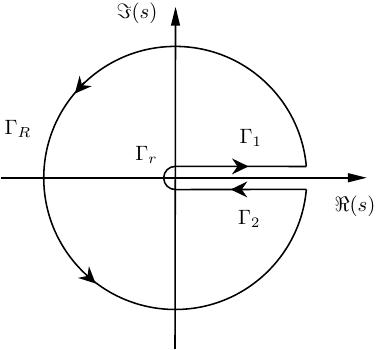}\\
    \footnotesize (a) Contour path \texorpdfstring{$\Gamma_A$}{GammaA}
  \end{minipage}\hfill
  \begin{minipage}[b]{.47\linewidth}
    \centering\includegraphics[width=\linewidth]{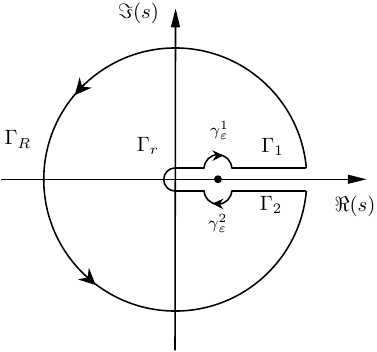}\\
    \footnotesize (b) Contour path \texorpdfstring{$\Gamma_B$}{GammaB}
  \end{minipage}
\caption{Two contour paths used for the Mellin-type integrals. Contour paths $\Gamma_A$ and $\Gamma_B$ are designed for the Laplace and Cauchy cases, respectively.}
\label{fig:contour}
\end{figure}

To begin with, we construct two contour paths, as shown in Figure \ref{fig:contour}. Unless indicated otherwise, we fix the branch cut of the logarithm along the positive real axis, i.e., $\ln s = \ln|s| + i\,\arg(s)$ with $\arg(s)\in[0,2\pi)$. Consider a keyhole contour $\Gamma_{A}=\Gamma_1 \cup \Gamma_2 \cup \Gamma_R \cup \Gamma_r$, as shown in Figure \ref{fig:contour}(a), where $\theta_0=\arcsin(r/R)$, $r<R$, and
\begin{align*}
    \Gamma_1 & =\left\{s+i r, 0 \leq s \leq \sqrt{R^2-r^2}\right\}, &\Gamma_2 & = \left\{s-i r, 0 \leq s \leq \sqrt{R^2-r^2}\right\}, \\
    \Gamma_R & =\left\{R e^{i \theta}, \theta_0<\theta<2 \pi-\theta_0\right\}, &\Gamma_{r} & =\left\{-r e^{i \theta},-\frac{\pi}{2}<\theta<\frac{\pi}{2}\right\}.
\end{align*}
Given $t>0$, we modify the contours $\Gamma_1$ and $\Gamma_2$ by making small semi-circular detours $\gamma_{\varepsilon}^1$ and $\gamma_{\varepsilon}^2$ around the point $s=t$:
\begin{align*}
    \gamma_{\varepsilon}^1 & =\left\{t + ir + \varepsilon e^{i \theta}, 0<\theta<\pi \right\}, \\
    \gamma_{\varepsilon}^2 & =\left\{t - ir + \varepsilon e^{i \theta}, -\pi<\theta<0 \right\}.
\end{align*}
Then we construct $\Gamma_{B}=\Gamma_1^{\varepsilon} \cup \Gamma_2^{\varepsilon} \cup \Gamma_R \cup \Gamma_r \cup \gamma_{\varepsilon}^1 \cup \gamma_{\varepsilon}^2$, as shown in Figure \ref{fig:contour}(b).

\begin{lemma}[Cauchy case]\label{lem3.1}
    Suppose $z \in \C$ with $\Re(z)>-1$ and $z \neq 0, 1, 2, \ldots$. Then
    \begin{align}\label{eqn31}
        I_0(z,t)&\triangleq\frac{1}{\pi}\,\mathrm{p.v.}\!\int_0^1\frac{s^z}{s-t}\mathrm{d}s=-\cot(\pi z) \cdot t^z + \frac{1}{\pi} \sum_{k=0}^{\infty}\frac{t^k}{z-k}
    \end{align}
    for $0<t<1$.
\begin{proof}
    Suppose $-1<\Re(z)<0$. Using Cauchy's integral theorem along the contour path $\Gamma_B$, we obtain
    $$\int_{\Gamma_B} \frac{s^z}{s-t} \mathrm{d}s = \left(\int_{\Gamma_1^{\varepsilon}}+\int_{\Gamma_2^{\varepsilon}}+\int_{\Gamma_R}+\int_{\Gamma_r}+\int_{\gamma_{\varepsilon}^1}+\int_{\gamma_{\varepsilon}^2}\right)\frac{s^z}{s-t} \mathrm{d}s=0.$$
    For $\Gamma_R$ and $\Gamma_r$, we have
    \begin{align*}
        & \left|\int_{\Gamma_R} \frac{s^z}{s-t} \mathrm{d}s\right| \leq \frac{2 \pi R^{\Re(z)+1}}{R-t} \rightarrow 0, \quad R \rightarrow \infty, \\
        & \left|\int_{\Gamma_r} \frac{s^z}{s-t} \mathrm{d}s\right| \leq \frac{ \pi r^{\Re(z)+1}}{t-r} \rightarrow 0, \quad r \rightarrow 0 .
    \end{align*}
    For the two small semi-circles, we first consider the limit of the contour integral as $r \rightarrow 0$. In this limit, the dumbbell contour collapses into a keyhole contour indented at $s=t$ on both the upper and lower sides of the branch cut. Subsequently, taking the limit as $\varepsilon \rightarrow 0$, we have
    \begin{align*}
        \int_{\gamma_{\varepsilon}^1}\frac{s^z}{s-t} \mathrm{d}s \rightarrow -i\pi t^z,\ \int_{\gamma_{\varepsilon}^2}\frac{s^z}{s-t} \mathrm{d}s \rightarrow -i\pi (te^{2\pi i})^z.
    \end{align*}
    Therefore, as $R \rightarrow \infty$, we can derive 
    \begin{equation*}
        (1-e^{2\pi iz})\int_0^\infty \frac{s^z}{s-t} \mathrm{d}s - i\pi t^z - i\pi (te^{2\pi i})^z=0,
    \end{equation*}
    which implies
    \begin{equation*}
        \int_0^\infty \frac{s^z}{s-t} \mathrm{d}s = -\pi\cot(\pi z) \cdot t^z,\ t>0.
    \end{equation*}
    
    By the geometric series expansion, we obtain
    $$
    \int_1^{\infty} \frac{s^z}{s-t} \mathrm{d}s=\int_1^{\infty} \sum_{k=0}^{\infty} s^{z-1}(t/s)^k \mathrm{d}s=-\sum_{k=0}^{\infty} \frac{t^k}{z-k},\ 0<t<1.
    $$
    Combining the above two equations, we conclude that
    $$
    \pi I_0(z,t)=\int_0^\infty\frac{s^z}{s-t}\mathrm{d}s-\int_1^\infty\frac{s^z}{s-t}\mathrm{d}s =-\pi\cot(\pi z) \cdot t^z + \sum_{k=0}^{\infty}\frac{t^k}{z-k},\ 0<t<1.
    $$
    By analytic continuation, the identity holds for all $z$ such that $\Re(z)>-1$ and $z\neq 0,1,2, \ldots$.
\end{proof}
\end{lemma}

\begin{lemma}[Laplace case]\label{lem3.2}
   Suppose $z\in\C$ with $\Re(z)>-1$, $z\neq 0,1,2, \ldots$, and $\theta\in(0,2\pi)$. Then 
    \begin{align}\label{eqn32}I_1(z,\theta,t)
    &\triangleq\frac{1}{\pi}\int_0^1\frac{s^z}{s^2+t^2-2st\cos\theta}\mathrm{d}s\nonumber\\
    &=\frac{\sin\left((\pi-\theta)z\right)}{\sin\left(\pi z\right)\sin\theta}t^{z-1}+\frac{1}{\pi}\sum_{k=1}^\infty\frac{\sin\left(k\theta\right)}{(z-k)\sin\theta}t^{k-1}.
    \end{align}
    for $0<t<1$.
\begin{proof}
    Suppose $-1<\Re(z)<0$. Using the residue theorem along the contour path $\Gamma_A$, we obtain
    \begin{align*}
        \int_{\Gamma_A} \frac{s^z}{s^2+t^2-2 s t \cos \theta}\mathrm{d}s = &\left(\int_{\Gamma_1}+\int_{\Gamma_2}+\int_{\Gamma_R}+\int_{\Gamma_r}\right) \frac{s^z}{\left(s-t e^{i\theta}\right)\left(s-t e^{i (2\pi-\theta)}\right)} \mathrm{d}s \\
        = &\ 2 \pi i\left(\left.\frac{s^z}{s-t e^{i (2 \pi-\theta)}}\right|_{s=t e^{i \theta}}+\left.\frac{s^z}{s-t e^{i \theta}}\right|_{s=t e^{i (2\pi-\theta)}}\right) \\
        = &\ \frac{\pi t^{z-1}}{\sin\theta}\left(e^{i \theta z}-e^{i (2\pi-\theta) z}\right).
    \end{align*}
    Therefore, we can derive
    \begin{equation*}
        (1-e^{2\pi iz})\int_0^{\infty} \frac{s^z}{s^2+t^2-2 s t \cos\theta}\mathrm{d}s = \frac{\pi t^{z-1}}{\sin\theta}\cdot\left(e^{i \theta z}-e^{i (2\pi-\theta) z}\right),
    \end{equation*}
    which implies
    \begin{equation*}
        \int_0^{\infty} \frac{s^z}{s^2+t^2-2 s t \cos\theta}\mathrm{d}s = \frac{\pi\sin\left((\pi-\theta)z\right)}{\sin\left(\pi z\right)\sin\theta}t^{z-1},\ 0<t<1.
    \end{equation*}
    Moreover, using the expansion
    \begin{align*}
        \int_1^{\infty} \frac{s^z}{s^2+t^2-2 s t \cos\theta}\mathrm{d}s &= \int_1^{\infty} \frac{s^z}{st \sin\theta}\frac{st \sin\theta}{s^2+t^2-2 s t \cos\theta} \mathrm{d}s\\
        &= \int_1^{\infty} \frac{s^{z-1}}{t \sin\theta}\sum_{k=1}^{\infty}(t/s)^k \sin(k\theta) \mathrm{d}s\\
        &= - \sum_{k=1}^\infty\frac{\sin\left(k\theta\right)}{(z-k)\sin\theta}t^{k-1}, \ 0<t<1,
    \end{align*}
    we can conclude that
    $$
    \pi I_1(z,\theta,t)=\frac{\pi\sin\left((\pi-\theta)z\right)}{\sin\left(\pi z\right)\sin\theta}t^{z-1}+\sum_{k=1}^\infty\frac{\sin\left(k\theta\right)}{(z-k)\sin\theta}t^{k-1},\ 0<t<1.
    $$
    By analytic continuation, the identity holds for all $z$ such that $\Re(z)>-1$ and $z\neq 0,1,2,\ldots$.
\end{proof}
\end{lemma}

\begin{lemma}[Stokes case]\cite{rachh2020solution}\label{lem3.3}
    Suppose $z\in\C$ with $\Re(z)>-1$, $z\neq 0,1,2,\ldots$, and $\theta\in(0,2\pi)$. Then
    \begin{align}\label{eqn33}
    I_2(z,\theta,t)&\triangleq\frac{1}{\pi}\int_0^1\frac{s^z}{(s^2+t^2-2st\cos\theta)^2}\mathrm{d}s\nonumber\\
    &=a(z,\theta)t^{z-3}-\sum_{k=1}^{\infty}F(k,z,\theta)t^k,
    \end{align}
    for $0<t<1$, where 
    \begin{align*}
    &a(z, \theta) = \frac{z \sin\left(2\theta\right) \cos\left(\left(\pi-\theta\right)z\right) + 2 \sin\left(\left(\pi-\theta\right)z\right)\left(1-z\sin^2\left(\theta\right)\right)}{4 \sin\left(\pi z\right) \sin^3\left(\theta\right)},\\
    &F(k,z,\theta)=\frac{(k+1) \sin\left((k+3)\theta\right) - (k+3) \sin\left((k+1)\theta\right)}{4\pi(-z+k+3)\sin^3\left(\theta\right)}.
    \end{align*}
\end{lemma}

\begin{remark}
    For Lemmas \ref{lem3.1}, and \ref{lem3.2}, \ref{lem3.3}, when $z$ is a non-negative integer or $\theta=\pi$, we can interpret the right-hand sides of Eq. \eqref{eqn31}-\eqref{eqn33} as their limiting values.
\end{remark}

\subsection{Roots of sinc function in the complex domain}
\label{subsec32}
In order to analyze the exponents of the density functions near corners, we need to consider the transcendental equation 
\begin{equation}\label{eq:trans2}
    G(\omega,\alpha)=\Sinc(\omega)+\mathrm{K} \Sinc(\alpha)=0,
\end{equation}
where $\omega\in\C$, $\alpha\in(0,2\pi)$, and $\mathrm{K}=\pm 1$ or $\pm \frac{\lambda+\mu}{\lambda+3\mu}$. 

\begin{figure}[!ht]
  \centering
  \begin{minipage}[b]{.6\linewidth}
    \centering
    \includegraphics[width=\linewidth]{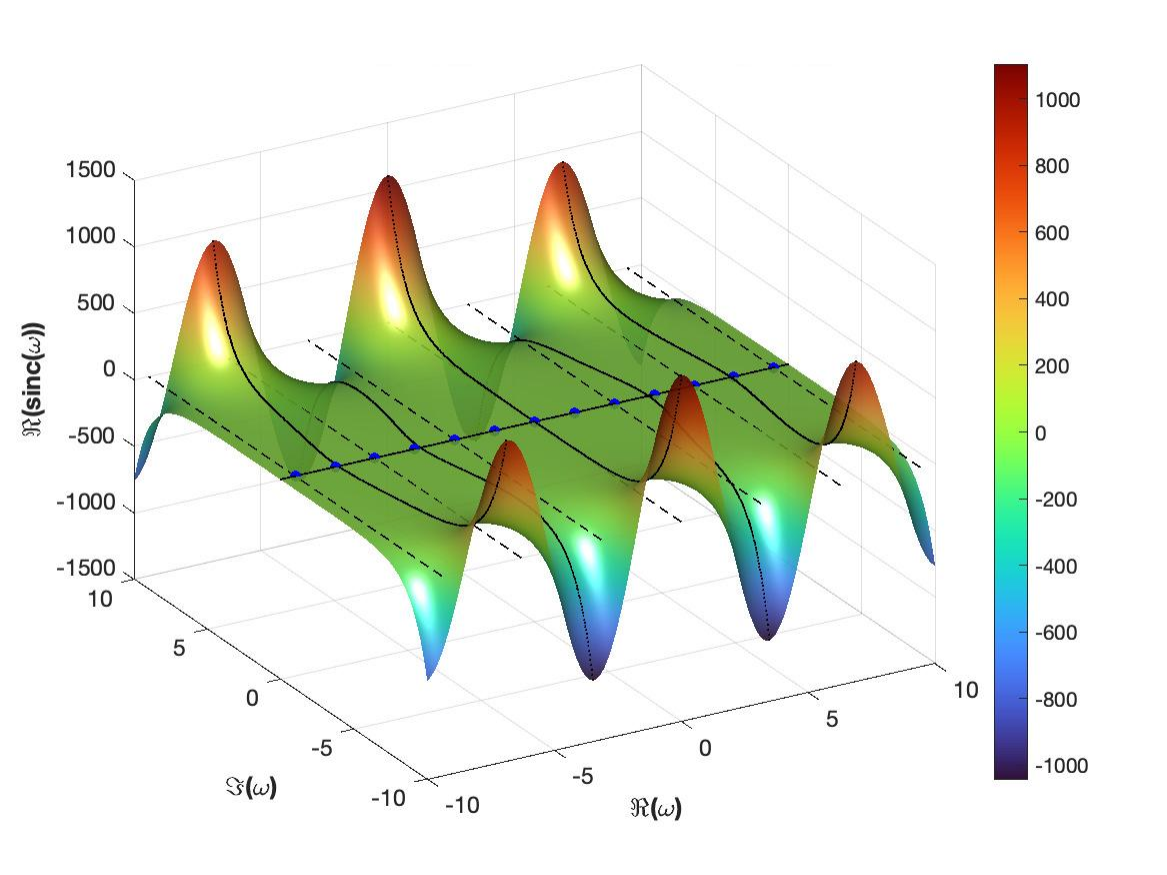}\\
    \footnotesize (a) Global view
  \end{minipage}

  \vspace{0.6\baselineskip} % 
  \begin{minipage}[b]{.47\linewidth}
    \centering
    \includegraphics[width=\linewidth]{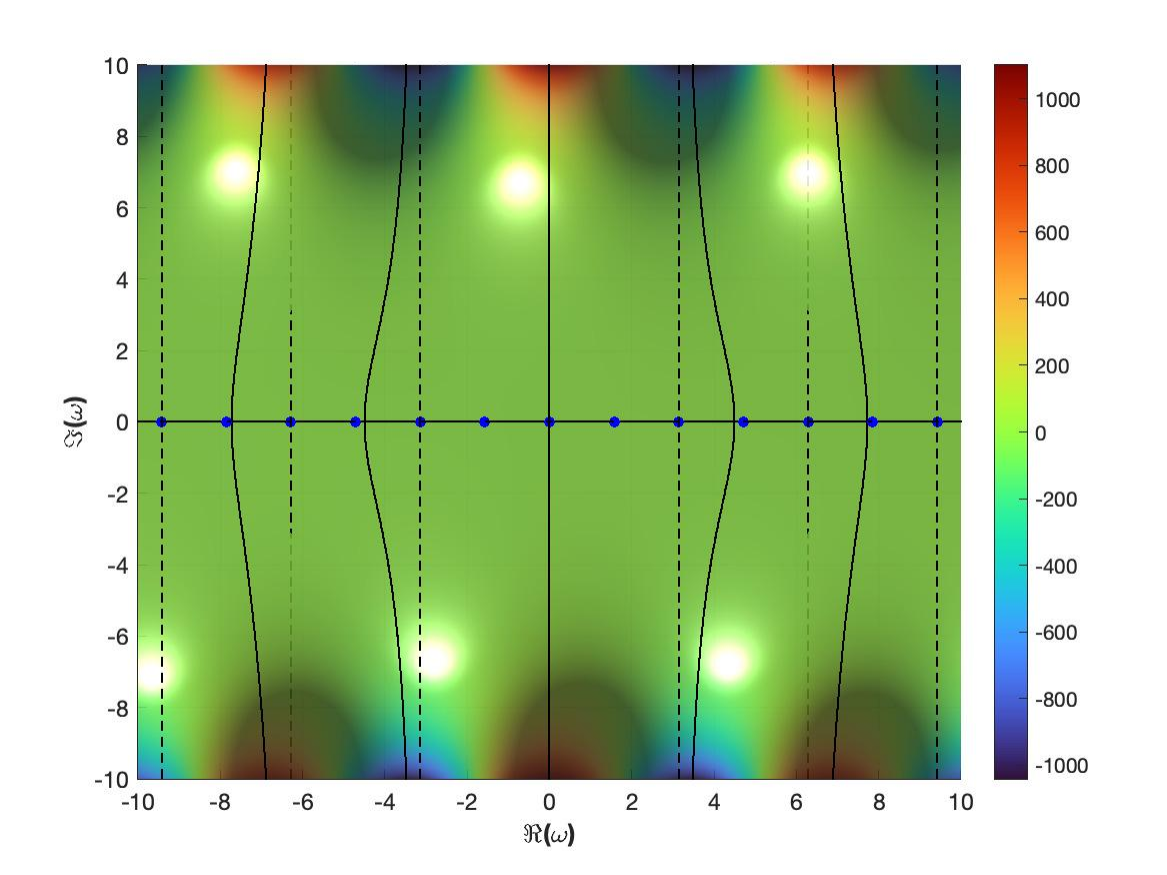}\\
    \footnotesize (b) Top view
  \end{minipage}\hfill
  \begin{minipage}[b]{.47\linewidth}
    \centering
    \includegraphics[width=\linewidth]{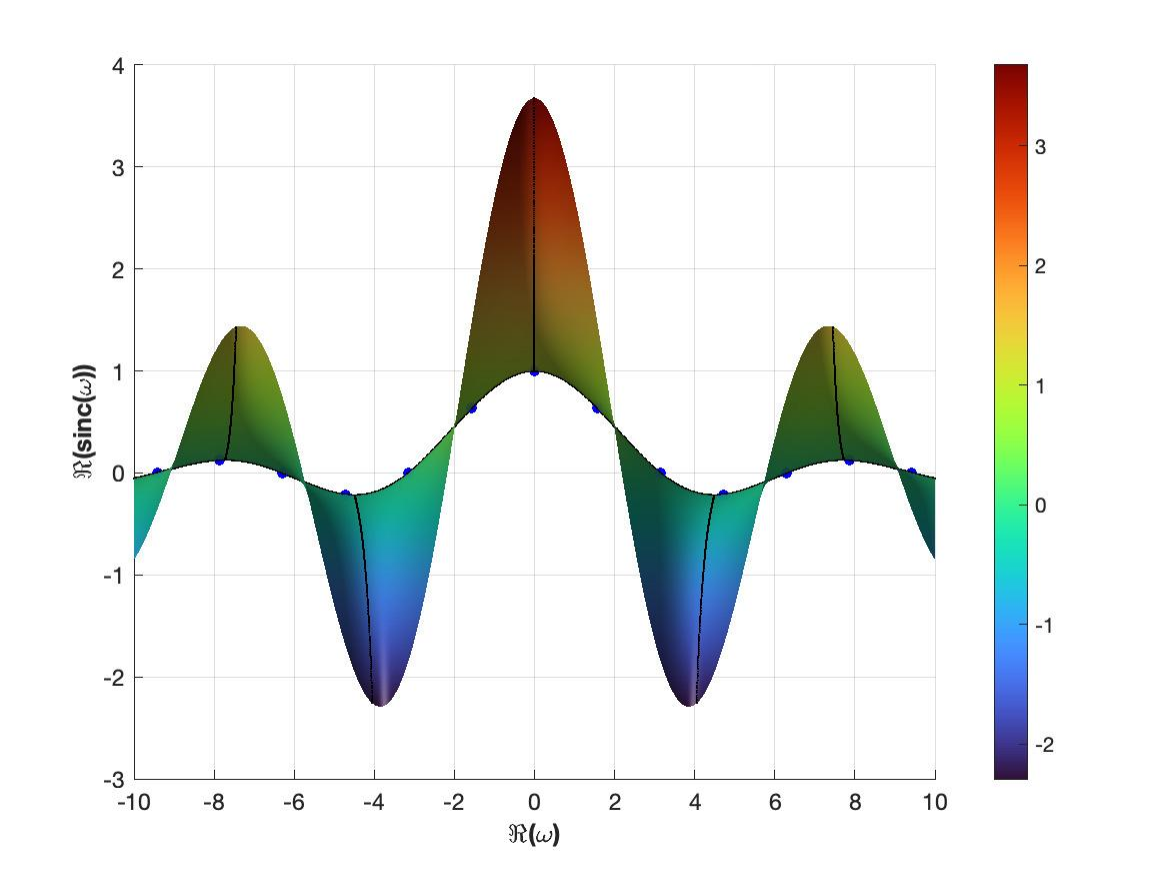}\\
    \footnotesize (c) Front view
  \end{minipage}
  \caption{Roots of the \texorpdfstring{$\Sinc$}{sinc} function in the complex plane.
  (a) 3D surface of \texorpdfstring{$\Sinc(x+iy)$}{sinc(x+iy)}.
  (b) Top view with the level set \texorpdfstring{$\Im(\Sinc(\omega))=0$}{Im(sinc(omega))=0} highlighted.
  (c) Front view when $|\Re(\Sinc(\omega))|\le \pi$.}
  \label{fig:Sinc_Roots}
\end{figure}

Since $\Sinc(\alpha)$ is real-valued, any solution $\omega$ of~\eqref{eq:trans2} must satisfy $\Im(\Sinc(\omega))=0$. Writing $\omega=x+iy$, we obtain
\begin{align*}
    &\Re(\Sinc(\omega))=\frac{x\sin{x}\cosh{y}+y\cos{x}\sinh{y}}{x^2+y^2},\\
    &\Im(\Sinc(\omega))=\frac{x\cos{x}\sinh{y}-y\sin{x}\cosh{y}}{x^2+y^2}.
\end{align*}
Thus $\Sinc(\omega)$ is purely real if and only if
\begin{equation}\label{eq:im_sinc_zero}
    x \cos x \sinh y - y \sin x \cosh y = 0.
\end{equation}
For $y \neq 0$ and $ x \neq 0$, it can be rearranged into a more compact form
\begin{equation}\label{eq:curve_equation}
    \frac{\tan x}{x} = \frac{\tanh y}{y}.
\end{equation}
Together with the real and imaginary axes, this relation characterizes the entire level set $\{\omega\in\C:\Im(\Sinc(\omega))=0\}$. As illustrated in Figure~\ref{fig:Sinc_Roots}(a), the black curves represent the contours where $\Im(\Sinc(\omega))=0$, with the top and front views shown in Figures~\ref{fig:Sinc_Roots}(b) and~\ref{fig:Sinc_Roots}(c), respectively. 
\begin{definition}\label{def:lambda_j}
Let $\{\lambda_j\}_{j=-\infty}^\infty$ denote the strictly increasing sequence of all real roots of the equation $\tan(x)=x$, indexed such that $\lambda_0 = 0$. For $j \neq 0$, these roots satisfy the bounds
$$|j|\pi + \frac{\pi}{4} < |\lambda_j| < |j|\pi + \frac{\pi}{2},$$
with $\mathrm{sgn}(\lambda_j) = \mathrm{sgn}(j)$.
\end{definition}

Using Definition \ref{def:lambda_j}, we can precisely characterize the properties of $\Sinc(\omega)$ along the curve~\eqref{eq:curve_equation}, as summarized in Lemma~\ref{lem3.4}.

\begin{lemma}\label{lem3.4}
For each $j\in\mathbb{Z}$, Eq.~\eqref{eq:curve_equation} defines a curve $\omega_j(y) = x_j(y) + iy$ in the complex plane. Each curve emanates from the point $\lambda_j$ on the real axis, i.e., $x_j(0) = \lambda_j$. Along these curves, the following properties hold:
\begin{enumerate}
    \item If $j$ is even, $\Sinc(\omega_j(y))$ is real and monotonically increasing with $|y|$, satisfying the parity relation $\Sinc(\omega_j(y))=\Sinc(\omega_j(-y))$ for all $y\in\R$.
    \item If $j$ is odd, $\Sinc(\omega_j(y))$ is real and monotonically decreasing with $|y|$, satisfying the same parity relation $\Sinc(\omega_j(y))=\Sinc(\omega_j(-y))$ for all $y\in\R$.
\end{enumerate}
\begin{proof}
    Taking the limit $y \to 0$ in Eq.~\eqref{eq:curve_equation} yields $\tan x = x$, confirming that each curve originates from $x_j(0)=\lambda_j$. To prove monotonicity, we simplify the expression for $\Sinc(\omega)$ along these curves. By substituting the condition $y \sin x \cosh y = x \cos x \sinh y$ into the numerator of $\Re(\Sinc(\omega))$, we obtain:
    $$
    \Sinc(\omega_j(y)) = \frac{\frac{x^2}{y}\cos x \sinh y + y\cos x \sinh y}{x^2+y^2} = \cos(x_j(y)) \frac{\sinh y}{y}.
    $$
    For $y>0$, Eq.~\eqref{eq:curve_equation} can be written as
    \begin{equation}\label{eqn_tanh}
    \frac{\tan x_j(y)}{x_j(y)} = \frac{\tanh y}{y}.
    \end{equation}
    The right-hand side is strictly decreasing from $1$ to $0$ as $y$ increases from $0$ to $+\infty$. On the other hand, on each interval $(|j|\pi,\,|j|\pi+\frac{\pi}{2})$ the function $\tan x/x$ is strictly increasing and positive. Therefore, as $y$ increases, equation \eqref{eqn_tanh} forces $|x_j(y)|$ to decrease monotonically from $|\lambda_j|$ towards $|j|\pi$, and thus $|\cos(x_j(y))|$ strictly increases with
    $\operatorname{sgn}(\cos(x_j(y))) = (-1)^j$.
    Meanwhile, $\sinh y / y$ is strictly increasing and positive for $y>0$. Hence,
    \[
    |\Sinc(\omega_j(y))| = |\cos(x_j(y))|\,\frac{\sinh y}{y}
    \]
    strictly increases. For even $j$, $\cos(\lambda_j)>0$, so $\Sinc(\omega_j(y))$ is positive and strictly increasing; for odd $j$, $\cos(\lambda_j)<0$, so $\Sinc(\omega_j(y))$ is negative and strictly decreasing. Using the evenness in $y$ established above, the same monotonicity holds as a function of $|y|$.
\end{proof}
\end{lemma}

Note that at $\alpha = \pi$ we have $\Sinc(\alpha) = 0$, so Eq.~\eqref{eq:trans2} reduces to $\Sinc(\omega) = 0$, whose real solutions are precisely $\omega = n\pi$, $n\in\mathbb{Z}$. We use these points as anchors when continuing the solution branches in $\alpha$. Now we analyze how a specific root $\omega$ travels along these curves as $\alpha$ varies. Along the curves where $\Im(\Sinc(\omega))=0$, the original transcendental equation reduces to finding the roots of a real-valued equation. To characterize the distribution of these roots and their dependence on the parameter $\alpha$, we explore the properties of the implicit function $\omega(\alpha)$ defined by $G(\omega, \alpha) = \Sinc(\omega) + \mathrm{K}\Sinc(\alpha) = 0$. The implicit function theorem guarantees that for a given point $(\omega_0, \alpha_0)$ satisfying $G(\omega_0, \alpha_0)=0$, we can find a unique analytic function $\omega(\alpha)$ in a neighborhood of $\alpha_0$ as long as $\partial G/\partial \omega \neq 0$. The ``problematic'' points, known as branch points, occur precisely where the following condition holds:
$$ \frac{\partial G}{\partial \omega} = \frac{d}{d\omega}\Sinc(\omega) = \frac{\omega\cos(\omega) - \sin(\omega)}{\omega^2} = 0.$$
This is equivalent to the equation $\tan(\omega)=\omega$, which means branch points in the $\omega$-plane occur exactly at $\omega = \lambda_j$ (as defined in Definition \ref{def:lambda_j}). This naturally leads to the definition of branch points in the physical parameter $\alpha$.

\begin{definition}\label{def:alpha_j}
    A value $\alpha \in (0, 2\pi)$ is called a branch point in the $\alpha$-plane if the solution $\omega(\alpha)$ coincides with a root $\lambda_j$ of $\tan\omega=\omega$. Specifically, for a given branch, the corresponding parameter $\alpha_{\lambda_j}$ is determined by
    \[
    \Sinc(\lambda_j) + \mathrm{K}\Sinc(\alpha_{\lambda_j}) = 0.
    \]
\end{definition}

% \begin{remark}
%     This framework allows us to understand the special values of $\alpha$ mentioned in the following theorems. A special value $\alpha_j$ is simply a value of $\alpha$ that causes the solution $\omega(\alpha)$ to coincide with a branch point, say $\omega_j$. They are determined by the relation $\Sinc(\omega_j) + \mathrm{K}\Sinc(\alpha_j) = 0$.
% \end{remark}

\begin{figure}[!ht]
  \centering
  \begin{minipage}[c]{.45\linewidth}
    \centering
    \includegraphics[width=\linewidth]{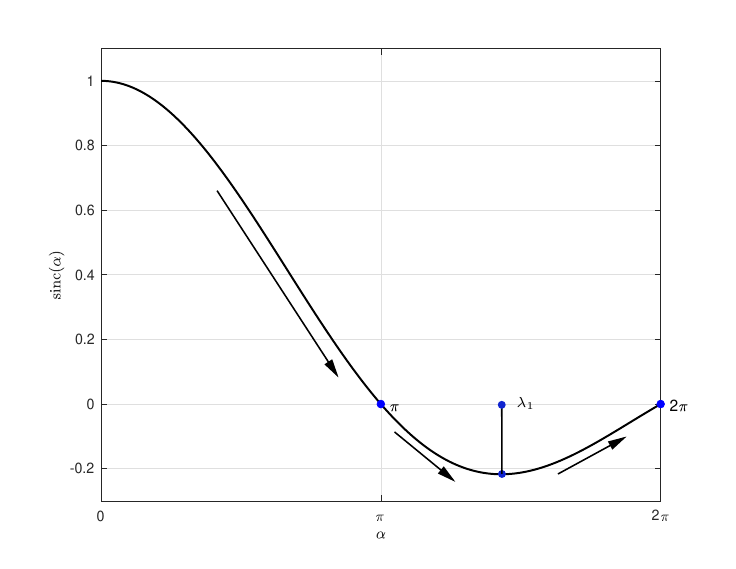}\\
    \footnotesize (a) $\Sinc(\alpha)$ in $(0,2\pi)$
  \end{minipage}\hfill
  \begin{minipage}[c]{.45\linewidth}
    \centering
    \includegraphics[width=\linewidth]{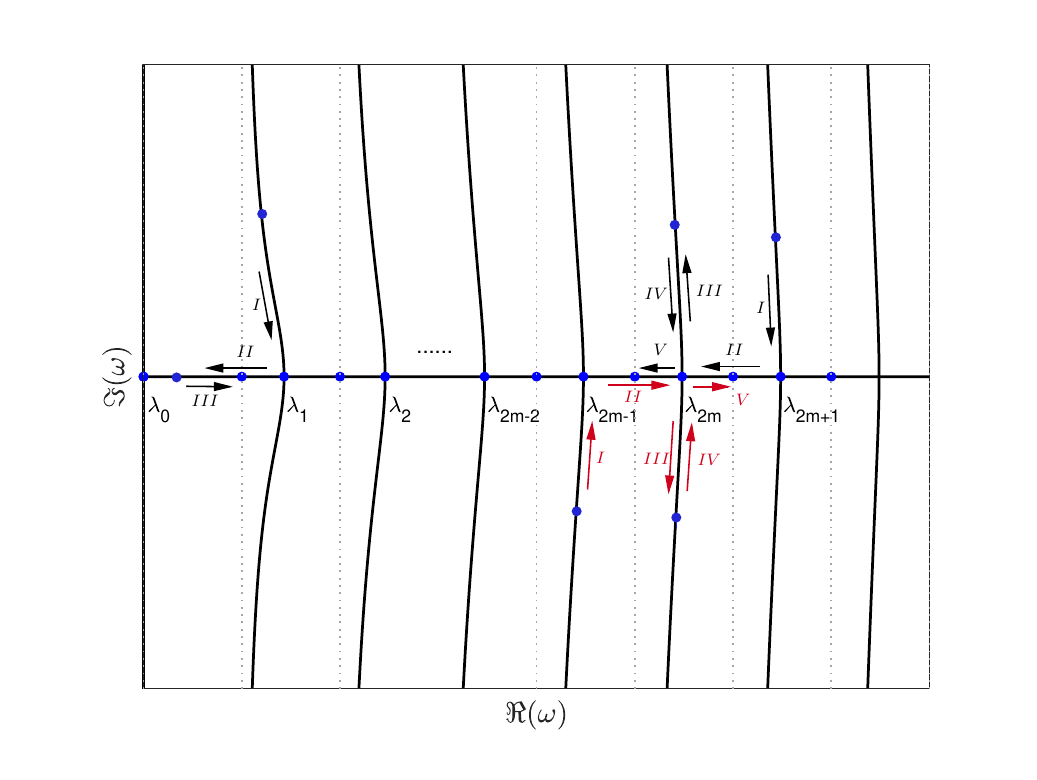}\\
    \footnotesize (b) $\Im(\Sinc(\omega))=0$
  \end{minipage}

  \caption{Implicit branches \texorpdfstring{$\omega(\alpha)$}{omega(alpha)} for
  \texorpdfstring{$\mathrm{K}=1$}{K=1}. (a)  The function \texorpdfstring{$\Sinc(\alpha)$}{sinc(alpha)} in \texorpdfstring{$(0,2\pi)$}{(0,2pi)}. (b) The level set \texorpdfstring{$\Im(\Sinc(\omega))=0$}{Im(sinc(omega))=0}.
  The arrows in (b) show the evolution as \texorpdfstring{$\alpha$}{alpha} increases:
  red for even \texorpdfstring{$n=2m$}{n=2m} (lower half-plane), black for odd
  \texorpdfstring{$n=2m+1$}{n=2m+1} (upper half-plane).}
  \label{fig:roots}
\end{figure}

To gain insight into the global behavior of these solution branches, we visualize the trajectories for the representative case $\mathrm{K}=1$ in Figure~\ref{fig:roots}. The behavior of $\Sinc(\alpha)$ over the interval $(0, 2\pi)$ is illustrated in Figure~\ref{fig:roots}(a). As the parameter $\alpha$ sweeps from $0$ to $2\pi$, the root $\omega(\alpha)$ anchored at $\omega(\pi) = n\pi$ traces a specific path in the complex plane. These trajectories exhibit distinct topological patterns based on the parity of $n$: for even indices $n=2m$ ($m\ge 1$), the roots follow the paths indicated by the red arrows in Figure \ref{fig:roots}(b). Conversely, for odd indices $n=2m+1$ ($m\ge 0$), the roots follow the trajectories marked by black arrows. 
Formalizing the existence and analyticity of these observed branches, we present the following theorem.
\begin{theorem}\label{thm3.5}
For $k\in\mathbb{N}^+$, let $\mathfrak{B}_k\subset(0,2\pi)$ denote the finite set of real branch points for the branch of solutions of $G(\omega,\alpha)=0$ anchored at $\omega(\pi)=k\pi$, i.e.\ those $\alpha$ for which there exists $\omega$ satisfying
\[
G(\omega,\alpha)=0,\qquad \partial_\omega G(\omega,\alpha)=0.
\]
Define
\[
D_k:=
\begin{cases}
\{\alpha\in\C:\ 0<\Re\alpha<2\pi,\ \Im\alpha\le 0\}\setminus \mathfrak{B}_k, & \text{if $k$ is even},\\[0.4em]
\{\alpha\in\C:\ 0<\Re\alpha<2\pi,\ \Im\alpha\ge 0\}\setminus \mathfrak{B}_k, & \text{if $k$ is odd}.
\end{cases}
\]
Then there exists a simply connected open set $V$ with $D_k\subset V\subset\C$ and an analytic function
$\omega_k(\alpha):V\to\C$ such that
\[
G(\omega_k(\alpha),\alpha)=0 \quad\text{for all }\alpha\in V,
\mbox{ and }
\omega_k(\pi)=k\pi.
\]
\end{theorem}

% -----------------------------
% Section 4: Corner analysis 
% -----------------------------

\section{Integral equations near corners}
\label{sec4}

\begin{figure}[!ht]
    \centering
    \includegraphics[width=.4\linewidth]{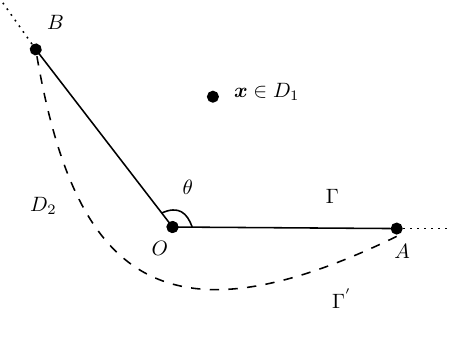} 
    \caption{Wedge domain $\Gamma=OA\cup OB$ with opening angle $\theta$ at $O$, parameterized by \eqref{eq:WedgePara}. Extending the edges to rays partitions $\R^2$ into $D_1$ (interior) and $D_2$ (exterior).}
    \label{fig:wedgecase}
\end{figure}

With all the analytic tools prepared above, we are ready to capture the singularities of Eq.~\eqref{eq:bie_original} on cornered boundaries. By localization and a partition-of-unity technique, the boundary integral operator near each corner reduces to the corresponding operator on a wedge. Suppose the wedge $\Gamma\subset\R^2$ is composed of two edges $OA$ and $OB$, parameterized by 
\begin{equation}\label{eq:WedgePara}
    \gamma(t)=\left\{
    \begin{array}{ll}
    (t, 0), & \ 0 \leq t \leq 1, \\
    (-t\cos\theta, -t\sin\theta), & \ -1 \leq t \leq 0. \\
    \end{array}\right.
\end{equation}
By extending the sides of the wedge to infinity, we divide $\R^2$ into two open sets $D_1$ and $D_2$, as shown in Figure~\ref{fig:wedgecase}.

We will identify a countable set of exponents $\{z_{n,j}\}_{n\ge1,\,j=1,\dots,4}\subset\C$ together with associated coefficient vectors $\mathbf{p}^{n,j}=[p_1^{n,j},p_2^{n,j},p_3^{n,j},p_4^{n,j}]^\top$, depending on $\theta$, such that, for any fixed $(n,j)$, if the density function $\bphi$ is defined by
$$
\bphi(t)=
\begin{cases}
    [p_1^{n, j}, p_2^{n, j}]^\top \cdot t^{z_{n, j}}, &0\leq t\leq 1,\\
    [p_3^{n, j}, p_4^{n, j}]^\top \cdot (-t)^{z_{n, j}}, &-1\leq t\leq 0,
\end{cases}
$$
then the left-hand side of Eq.~\eqref{eq:bie_original} is smooth near $t=0$.

We will also prove the converse. Specifically, let $N$ be a positive integer, and suppose further that $\beta_{n,j} \in \C$, for $j=1,2,3,4$, and $n=0,1,2, \ldots, N-1$, with $\boldsymbol{f}$ satisfying
$$
\boldsymbol{f}(t)=
\begin{cases}
    \sum\limits_{n=0}^{N-1}[\beta_{n, 1} , \beta_{n,2} ]^\top \cdot t^n, &0\leq t\leq 1,\\
    \sum\limits_{n=0}^{N-1}[\beta_{n, 3} , \beta_{n,4} ]^\top \cdot (-t)^n, &-1\leq t\leq 0.
\end{cases}
$$
Then, for all but countably many values of $0 < \theta < 2\pi$, there exist unique coefficients $\alpha_{n,j} \in \C$, for $j=1,2,3,4$ and $n=1, \ldots N$, such that the density function defined by
$$
\bphi(t)=\begin{cases}
    \sum\limits_{n=1}^N\sum\limits_{j=1}^4 \alpha_{n,j}[p_1^{n, j}, p_2^{n, j}]^\top \cdot t^{z_{n, j}}, &0\leq t\leq 1,\\
    \sum\limits_{n=1}^N\sum\limits_{j=1}^4 \alpha_{n,j}[p_3^{n, j}, p_4^{n, j}]^\top \cdot (-t)^{z_{n, j}}, &-1\leq t\leq 0,
\end{cases}
$$
satisfies Eq.~\eqref{eq:bie_original} with error $O\left(|t|^{N}\right)$. 

\subsection{Boundary integral equations on the wedge}
\label{subsec41} 
To derive the boundary integral equations, we place the point $x$ on the edges $OA$ and $OB$ respectively, from which we can derive the corresponding boundary integral operators for $x\in OA$ and $x\in OB$. Assume that $\Re(z)>-1$ and $z\notin\mathbb{Z}$, and set
\begin{equation}\label{eq:density}
    \bphi(t)=
    \begin{cases}
    [p_1,p_2]^\top \cdot t^z, &0\leq t\leq 1,\\
    [p_3,p_4]^\top \cdot (-t)^z, &-1\leq t\leq 0.
    \end{cases}
\end{equation}
To maintain brevity, denote 
$$c_1 \triangleq \frac{\mu}{2(\lambda+2\mu)}, \quad c_2 \triangleq \frac{\lambda+\mu}{\lambda+2\mu},$$ 
$$I_{OA} \triangleq \int_{OA}\mathbf{D}(\bx,\by)\bphi(\by)\mathrm{d}s_{\by}, \quad I_{OB} \triangleq \int_{OB}\mathbf{D}(\bx,\by)\bphi(\by)\mathrm{d}s_{\by}.$$
For $x\in OA$, we have $x=(t,0)$ where $0\leq t\leq 1$, then 
\begin{align*}
    I_{OA} &=\frac{c_1}{\pi}\int_0^1 \frac{s^z}{t-s}\begin{bmatrix}p_2\\-p_1\end{bmatrix}\mathrm{d}s,\\
    I_{OB} &= \frac{c_1}{\pi}\left\{\int_{-1}^0\frac{-t\sin\theta\cdot(-s)^z}{ t^2+s^2+2ts\cos\theta}\begin{bmatrix}p_3\\p_4\end{bmatrix}\mathrm{d}s
    + \int_{-1}^0\frac{(-t\cos\theta-s)\cdot(-s)^z}{t^2+s^2+2ts\cos\theta}\begin{bmatrix}p_4\\-p_3\end{bmatrix}\mathrm{d}s\right\} \\
    & + \frac{c_2}{\pi} \int_{-1}^0\frac{-t\sin\theta\cdot(-s)^z}{(t^2+s^2+2ts\cos\theta)^2} 
    \begin{bmatrix}
    (t+s\cos\theta)^2 & (t+s\cos\theta)s\sin\theta \\
    (t+s\cos\theta)s\sin\theta &s^2\sin^2\theta
    \end{bmatrix}
    \begin{bmatrix}
    p_3\\p_4
    \end{bmatrix}\mathrm{d}s,
\end{align*}
which further gives the boundary double-layer operator on $OA$:
\begin{align*}
    \mathcal{D}_{OA}(t) = c_1 
    \left\{ - I_0(z,t)\cdot \begin{bmatrix}p_2\\-p_1\end{bmatrix} + I_3 \begin{bmatrix}p_3\\p_4\end{bmatrix}
    + I_4 \begin{bmatrix}p_4\\-p_3\end{bmatrix}\right\} 
     + c_2
    \begin{bmatrix}
    I_{11} & I_{12}\\
    I_{21} & I_{22}
    \end{bmatrix}
    \begin{bmatrix}
    p_3\\p_4
    \end{bmatrix},
\end{align*}
where $0\leq t\leq 1$ and
\begin{align*}
    &I_3 \triangleq - t\sin\theta \cdot I_1(z,\theta,t),\\
    &I_4 \triangleq -t\cos\theta\cdot I_1(z,\theta,t)+I_1(z+1,\theta,t),\\
    &I_{11} \triangleq -t^3\sin\theta\cdot I_2(z,\theta,t)+2t^2\sin\theta\cos\theta\cdot I_2(z+1,\theta,t)-t\sin\theta\cos^2\theta\cdot I_2(z+2,\theta,t),\\
    &I_{12} \triangleq I_{21} = t^2\sin^2\theta\cdot I_2(z+1,\theta,t) - t\cos\theta\sin^2\theta\cdot I_2(z+2,\theta,t),\\
    &I_{22} \triangleq -t\sin^3\theta\cdot I_2(z+2,\theta,t).
\end{align*}
For $x\in OB$, we have $x=-t(\cos\theta,\sin\theta)$ where $-1\leq t\leq 0$, then 
\begin{align*}
    I_{OB} &= \frac{c_1}{\pi}\int_{-1}^0 \frac{(-s)^z}{t-s} \begin{bmatrix}p_4\\-p_3\end{bmatrix}\mathrm{d}s, \\
    I_{OA} &= \frac{c_1}{\pi}\left\{\int_0^1\frac{t\sin\theta\cdot s^z}{ t^2+s^2+2ts\cos\theta}\begin{bmatrix}p_1\\p_2\end{bmatrix}\mathrm{d}s
    + \int_0^1\frac{(-t\cos\theta-s)\cdot s^z}{t^2+s^2+2ts\cos\theta}\begin{bmatrix}p_2\\-p_1\end{bmatrix}\mathrm{d}s\right\} \\
    & + \frac{c_2}{\pi} \int_0^1\frac{t\sin\theta\cdot s^z}{(t^2+s^2+2ts\cos\theta)^2} 
    \begin{bmatrix}
    (t\cos\theta+s)^2 & (t\cos\theta+s)t\sin\theta \\
    (t\cos\theta+s)t\sin\theta &t^2\sin^2\theta
    \end{bmatrix}
    \begin{bmatrix}
    p_1\\p_2
    \end{bmatrix}\mathrm{d}s.
\end{align*}
Replace $t$ by $-t$, $0\leq t\leq 1$, then we can obtain the boundary double-layer operator on $OB$:
\begin{align*}
    \mathcal{D}_{OB}(-t) =& c_1\left\{I_0(z,t)\begin{bmatrix}
    p_4\\-p_3
    \end{bmatrix}+ I_3 \begin{bmatrix}p_1\\p_2\end{bmatrix}
    - I_4 \begin{bmatrix}p_2\\-p_1\end{bmatrix}\right\} 
    + c_2
    \begin{bmatrix}
    I_{11}^* & I_{12}^*\\
    I_{21}^* & I_{22}^*
    \end{bmatrix}
    \begin{bmatrix}
    p_1\\p_2
    \end{bmatrix},
\end{align*}
where $0\leq t\leq 1$ and 
\begin{align*}
    I_{11}^* &\triangleq -t^3\sin\theta\cos^2\theta\cdot I_2(z,\theta,t)+2t^2\sin\theta\cos\theta\cdot I_2(z+1,\theta,t)-t\sin\theta\cdot I_2(z+2,\theta,t),\\
    I_{12}^* &\triangleq I_{21}^* = -t^3\cos\theta\sin^2\theta\cdot I_2(z,\theta,t) + t^2\sin^2\theta\cdot I_2(z+1,\theta,t),\\
    I_{22}^* &\triangleq  -t^3\sin^3\theta\cdot I_2(z,\theta,t).
\end{align*}

For the sake of brevity, denote the left-hand side of Eq.~\eqref{eq:bie_original} by $\boldsymbol{h}(\bx)$. Then we derive
\begin{equation}\nonumber
    \boldsymbol{h}(\bx(t)) = -\frac{1}{2}\bphi(\bx(t)) + \mathcal{D}[\bphi](t),\ -1\leq t\leq 1.
\end{equation}
Replacing $t$ by $-t$ on the interval $[-1,0]$, we obtain a linear system of four equations involving the coefficient vector $\mathbf{p}=[p_1, p_2, p_3, p_4]^\top$:
\begin{equation*}
\begin{cases}
    \boldsymbol{h}(\bx(t)) = -\frac{1}{2}[p_1,p_2]^\top\cdot t^z + \mathcal{D}_{OA}(t), \\
    \boldsymbol{h}(\bx(-t)) = -\frac{1}{2}[p_3,p_4]^\top\cdot t^z + \mathcal{D}_{OB}(-t), 
\end{cases} 0\leq t\leq 1.
\end{equation*}
More specifically, for $0\leq t\leq 1$, we have
\begin{equation}\nonumber
    \tilde{\boldsymbol{h}}(t)\triangleq
    \begin{bmatrix}
        h_1(t)\\h_2(t)\\h_1(-t)\\h_2(-t)
    \end{bmatrix}=
    -\frac{1}{2} \begin{bmatrix}
    p_1\\p_2\\p_3\\p_4
    \end{bmatrix} t^z + \mathbf{G}_{\mathrm{DL}}(z,\theta,t)
    \begin{bmatrix}
    p_1\\p_2\\p_3\\p_4
    \end{bmatrix},
\end{equation}
where 
\begin{equation*}
    \mathbf{G}_{\mathrm{DL}}(z,\theta,t)\triangleq
    \begin{bmatrix}
    0&-c_1 I_0&c_1 I_3+c_2 I_{11}&c_1 I_4+c_2 I_{12}\\    
    c_1 I_0&0&-c_1 I_4+c_2 I_{12}&c_1 I_3+c_2 I_{22}\\
    c_1 I_3+c_2 I_{11}^*&-c_1 I_4+c_2 I_{12}^*&0&c_1 I_0\\
    c_1 I_4+c_2 I_{12}^*&c_1 I_3+c_2 I_{22}^*&-c_1 I_0&0
    \end{bmatrix}.
\end{equation*}
Upon inserting the explicit formulas of 
$I_0,I_1,I_2,I_3,I_4,I_{11},I_{12},I_{22},I_{11}^*,I_{12}^*,I_{22}^*$ via Lemmas \ref{lem3.1}, \ref{lem3.2} and \ref{lem3.3}, we obtain
\begin{equation*}
    \tilde{\boldsymbol{h}}(t) = \mathbf{A}(z,\theta)\mathbf{p} \cdot t^z + \sum_{k=0}^{\infty}\mathbf{F}(k,z,\theta) \mathbf{p} \cdot t^k,
\end{equation*}
where the explicit closed forms of the matrix $\mathbf{A}(z,\theta)$ and $\mathbf{F}(k,z,\theta)$ are collected in Appendix~\ref{app:AF_calculation}. The above results are summarized in the following theorem.

\begin{theorem}\label{thm4.1}
    Suppose $z\notin\mathbb{Z}$ and $\Re(z)>-1$. If the density function $\bphi$ is given by the ansatz~\eqref{eq:density}, then the left-hand side of boundary integral equation~\eqref{eq:bie_original} on the wedge~\eqref{eq:WedgePara} can be written as
\begin{equation}\label{eq:MapEquation}
    \tilde{\boldsymbol{h}}(t) = \left[ \mathbf{A}(z,\theta)\cdot
    t^z + \sum_{k=0}^{\infty}\mathbf{F}(k,z,\theta)\cdot t^k \right]\mathbf{p},
\end{equation}
where $0\leq t\leq 1$. Furthermore, if $\boldsymbol{h}$ is smooth, then the singular exponents $z$ must satisfy the condition $\det(\mathbf{A}(z,\theta)) = 0$, which yields the following transcendental equation:
    \begin{equation}\label{eq:trans1}
         \left(z^2\sin^2\theta - \sin^2(\left(2\pi -\theta\right)z)\right)\left(z^2\sin^2\theta-c^2\sin^2(\theta z)\right)=0,
    \end{equation}
where $c=\frac{\lambda+3\mu}{\lambda+\mu}$ and the coefficient vector $\mathbf{p}$ satisfies $\mathbf{p}\in \ker(\mathbf{A}(z,\theta))$.
\end{theorem}
\begin{remark}\label{rmk4.2}
    The factorization in~\eqref{eq:trans1} shows that the characteristic condition for the boundary integral operator splits into two Lam\'e-type transcendental families,
    \[
    H_{2\pi-\theta}(z):\ z^{2}\sin^{2}(2\pi-\theta) = \sin^{2}\!\big((2\pi-\theta)z\big),
    \qquad
    H_{\theta}(z):\ z^{2}\sin^{2}\theta = c^{2}\sin^{2}(\theta z).
    \]
    These two families encode, in a dual manner, the corner spectra of the interior Dirichlet and exterior Neumann problems as seen through the boundary integral representations. Which family actually appears in the boundary equation depends on both the potential representation and the boundary data. Under the artificial point source configuration used in our numerical experiments, the indirect double-layer formulation (IDP) selects precisely the family $H_{2\pi-\theta}(z)$, while by adjoint duality the single-layer formulation for the exterior Neumann problem (ENP) selects the family $H_{\theta}(z)$. However, in the absence of analytic boundary data, both families may contribute. 
\end{remark}

\subsection{Existence of the implicit function \texorpdfstring{$z(\theta)$}{z(theta)}}\label{subsec42}
Now we analyze the implicit functions $z(\theta)$ defined by the transcendental equation \eqref{eq:trans1}, which is closely connected to Eq.~\eqref{eq:trans2}, as shown in the following lemma.
\begin{lemma}\label{lem4.2}
    Through changes of variables, the following relationships can be established:
    \begin{enumerate}
    \item Substituting $\alpha = 2\pi-\theta$ and $\omega = \alpha z$ into the equations $H_1(z,\theta) = z\sin\theta - \sin(\left(2\pi -\theta\right)z) = 0$ and $H_2(z,\theta) = z\sin\theta + \sin(\left(2\pi -\theta\right)z) = 0$, we obtain
    $$G(\omega,\alpha)=\Sinc(\omega)\pm \Sinc(\alpha)=0.$$
    \item Substituting $\alpha = \theta$ and $\omega = \alpha z$ into the equations $H_3(z,\theta) = z\sin\theta + c\sin(\theta z) = 0$ and $H_4(z,\theta) = z\sin\theta - c\sin(\theta z)  = 0$, we obtain
    $$G(\omega,\alpha)=\Sinc(\omega)\pm \frac{\lambda+\mu}{\lambda+3\mu}\Sinc(\alpha)=0.$$
    \end{enumerate}
\end{lemma}

Consider the case $H_1(z, \theta)$. Combining Lemma~\ref{lem4.2} with Theorem~\ref{thm3.5}, we can establish the existence of an implicit function $z_{n,1}(\theta)$ such that $H_1\bigl(z_{n,1}(\theta),\theta\bigr)=0$.
\begin{theorem}\label{thm4.3}
Let $n\in\mathbb{N}^+$. Consider the equation $H_1(z,\theta)=0$. Let $\Theta_n\subset(0,2\pi)$ denote the finite set of real branch-point parameters for the branch anchored at $z(\pi)=n$. Let the domain $D_n$ be:
\[
D_n:=
\begin{cases}
\{\theta\in\C:\ 0<\Re\theta<2\pi,\ \Im\theta\le 0\}\setminus \Theta_n, & \text{if $n$ is even},\\[0.4em]
\{\theta\in\C:\ 0<\Re\theta<2\pi,\ \Im\theta\ge 0\}\setminus \Theta_n, & \text{if $n$ is odd}.
\end{cases}
\]
Then there exists a simply connected open set $V_n$ with $D_n\subset V_n\subset\C$ and an analytic function $z_{n,1}(\theta):V_n\to\C$ such that
\[
H_1\bigl(z_{n,1}(\theta),\theta\bigr)=0 \quad \text{for all } \theta\in V_n, \quad \text{and} \quad z_{n,1}(\pi)=n.
\]
\end{theorem}

\begin{remark}\label{rmk4.5}
    Similar results hold for each $H_j(z,\theta)=0$ with $j=1,2,3,4$. Specifically, for any $n\in\mathbb{N}^+$, there exist analytic functions $z_{n,j}(\theta):V_j \rightarrow \C$ satisfying $H_j(z_{n,j}(\theta),\theta)=0$, for all $\theta \in V_j$ and $z_{n,j}(\pi)=n$, where each $V_j \subset \C$ is a simply connected open set. 
\end{remark}

\begin{example}[Numerical roots]\label{ex:example1}
    To numerically quantify how the opening angle $\theta$ influences the corner singularity of the density function, we compute for each branch $H_1$--$H_4$ the root $z$ with the smallest positive real part. In accordance with the finite-energy requirement for the displacement field, it requires $\Re(z)>0$ for the double-layer density. Using $\lambda=1$, $\mu=2$ and angles $\theta\in\{\pi/4,\ \pi/2,\ 3\pi/4,\ 5\pi/4,\ 3\pi/2,\ 7\pi/4\}$, the resulting dominant exponents are listed in Table~\ref{tab:parametric_roots}.
    \begin{table}[!ht]
      \centering
      \caption{Dominant singular exponents $z_k$ (rounded to 4 decimals) for wedge angle $\theta$. Entries are listed in the order of four characteristic equations $H_1\text{--}H_4$. Bold marks the smallest positive $\Re z$.}
      \label{tab:parametric_roots}
      \setlength{\tabcolsep}{6pt}
      \renewcommand{\arraystretch}{1.15}
      \begin{tabular}{@{} c c c c c @{}}
        \toprule
        $\theta$ & $z_1$ & $z_2$ & $z_3$ & $z_4$ \\
        \midrule
        $\pi/4$   & \(\mathbf{0.5050}\) & \(0.6597\) & \(5.6004 \pm 1.5055i\) & \(2.7474\) \\
        $\pi/2$   & \(\mathbf{0.5445}\) & \(0.9085\) & \(2.8381 \pm 0.4470i\) & \(1.5408\) \\
        $3\pi/4$  & \(\mathbf{0.6736}\) & \(1.0000\) & \(1.5393\)             & \(1.1784\) \\
        $5\pi/4$  & \(1.8854 \pm 0.3607i\) & \(1.0000\) & \(\mathbf{0.7422}\) & \(0.8678\) \\
        $3\pi/2$  & \(2.7396 \pm 1.1190i\) & \(1.0000\) & \(\mathbf{0.6105}\) & \(0.7346\) \\
        $7\pi/4$  & \(5.3905 \pm 2.7204i\) & \(1.0000\) & \(\mathbf{0.5414}\) & \(0.6050\) \\
        \bottomrule
      \end{tabular}
    \end{table}

    Table~\ref{tab:parametric_roots} shows that, for all tested angles, the dominant real parts are positive, so the boundary density exhibits a non-oscillatory algebraic decay $t^{\Re z}$ toward the corner. The dependence on $\theta$ is non-monotone: the dominant values range from $0.5050$ at $\theta=\pi/4$ to $0.6736$ at $\theta=3\pi/4$, and decrease again toward re-entrant configurations. Besides, complex roots appear on some branches (e.g., $\theta=5\pi/4$), but their real parts exceed the dominant one at the same angle and thus remain higher-order with logarithmic oscillations $t^{a}\{\cos(b\ln t),\sin(b\ln t)\}$.

    Especially, the choice of the dominant exponent with $\Re(z)>0$ follows from the link between the PDE field and the boundary density via the boundary integral representations. By Green's representation theorem,
    $$
    \boldsymbol{u}_{in}=D[\gamma \boldsymbol{u}_{in}] - S[\mathcal{T}\boldsymbol{u}_{in}],
    \quad
    \boldsymbol{u}_{ex}=-D[\gamma \boldsymbol{u}_{ex}] + S[\mathcal{T}\boldsymbol{u}_{ex}],
    $$
    the double-layer term preserves the leading exponent (density and field have the same corner power, so $u\sim r^{\,z}$ if $\psi\sim r^{\,z}$), whereas the single-layer term raises it by one (the field gains one power relative to the density, so $u\sim r^{\,z+1}$ if $\varphi\sim r^{\,z}$). Consequently, admissibility of the displacement field in $H^1$ enforces $\Re(z)>0$ for the double-layer density, while for a single-layer representation it translates to $\Re(z)>-1$. The present table adopts the double-layer potential and therefore selects, across the four branches, the smallest positive real part as the dominant contribution for each angle.
\end{example}

\subsection{Invertibility of mapping matrix}\label{subsec43}
Suppose $\mathbf{B}(\theta)$ is the mapping matrix transforming the coefficients of the truncated basis functions
$$(p_1^{n,j} t^{z_{n,j}}, p_2^{n,j} t^{z_{n,j}}, p_3^{n,j} t^{z_{n,j}}, p_4^{n,j} t^{z_{n,j}}),\ n=1,\ldots,N,\ j=1,2,3,4$$
to Taylor expansion coefficients of the vector $\tilde{\boldsymbol{h}}(t)$. The explicit form of $\mathbf{B}(\theta)$ is established in the following theorem.
\begin{theorem}\label{thm4.5}
    Let $N \geq 2$ be an integer. For $n=1,2, \ldots, N$, $j=1,2,3,4$, suppose that $z_{n,j}(\theta)$ are analytic functions satisfying $\det \mathbf{A}\left(z_{n,j}(\theta), \theta\right)=0$ for $\theta \in V_j \subset \C$, with $V_j$ defined in Remark \ref{rmk4.5}, and that $\mathbf{p}^{n,j} = \left(p_1^{n,j}, p_2^{n,j}, p_3^{n,j}, p_4^{n,j}\right)^\top \in \ker(\mathbf{A}\left(z_{n,j}(\theta), \theta\right))$. Moreover, suppose the density $\bphi$ satisfies
    \begin{align*}
        \begin{bmatrix}
            \phi_1(t)\\ \phi_2(t)\\ \phi_1(-t)\\ \phi_2(-t)
        \end{bmatrix}=
        \sum_{n=1}^N\sum_{j=1}^4 \alpha_{n,j}
        \begin{bmatrix}
             p_1^{n,j}\\
             p_2^{n,j}\\
             p_3^{n,j}\\
             p_4^{n,j}
        \end{bmatrix}t^{z_{n,j}},\ 0<t<1,
    \end{align*}
    where $\boldsymbol{\alpha}_n = \left(\alpha_{n,1}, \alpha_{n,2}, \alpha_{n,3}, \alpha_{n,4}\right)^\top \in \C^4, n = 1, \ldots, N$. Then the left-hand side of Eq.~\eqref{eq:bie_original} admits the asymptotic expansion
    \begin{align*}
    \tilde{\boldsymbol{h}}(t) = \sum_{k=0}^{N-1}\boldsymbol{\beta}_k t^k + O(t^{N}), \ 0<t<1,
    \end{align*}
    where $\boldsymbol{\beta}_k = \left(\beta_{k,1}, \beta_{k,2}, \beta_{k,3}, \beta_{k,4}\right)^\top\in\C^4$, $k = 0, \ldots, N-1$ satisfy
    \begin{align*}
        \begin{bmatrix}
            \boldsymbol\beta_0\\ \vdots\\ \boldsymbol\beta_{N-1}
        \end{bmatrix} =
        \begin{bmatrix}
            \mathbf B_{0,1}(\theta)&\cdots&\mathbf B_{0,N}(\theta)\\
            \vdots&\ddots&\vdots\\
            \mathbf B_{N-1,1}(\theta)&\cdots&\mathbf B_{N-1,N}(\theta)
            \end{bmatrix}
            \begin{bmatrix}
            \boldsymbol\alpha_1\\ \vdots\\ \boldsymbol\alpha_N
        \end{bmatrix},
    \end{align*}
    More explicitly, the $4N\times 4N$ mapping matrix $\mathbf{B}(\theta)$ consists of $4 \times 4$ block matrices $\mathbf{B}_{k,n}(\theta)$:
    $$\mathbf{B}_{k,n}(\theta)=
    \left[
    \begin{array}{c|c|c}
     \mathbf{F}\!\left(k,z_{n,1}(\theta),\theta\right)\mathbf{p}^{n,1}
     & \cdots &
     \mathbf{F}\!\left(k,z_{n,4}(\theta),\theta\right)\mathbf{p}^{n,4}
    \end{array}
    \right].$$
\end{theorem}

% 2. 清楚映射的具体形式之后，为了构建定理，我们还需要证明矩阵B的可逆性
The following two lemmas establish the invertibility of $\mathbf{B}(\theta)$ for almost every angle.
\begin{lemma}\label{lem4.6}
    Let $V = \cap_{j=1}^4 V_j$ be the domain where the root functions $z_{n,j}(\theta)$ are analytic. The function $\det(\mathbf{B}(\theta))$ is analytic for $\theta \in V$. For the special case $\theta=\pi$, the matrix is nonsingular, i.e., $\det(\mathbf{B}(\pi)) \neq 0$.
\begin{proof}
    The determinant of a matrix is a polynomial of its entries, so its analyticity depends on the analyticity of the matrix elements themselves. The entries of $\mathbf{B}(\theta)$ are constructed from the matrix-valued function $\mathbf{F}$ and the null vectors $\mathbf{p}^{n,j}$, evaluated at the root functions $z_{n,j}(\theta)$. The analyticity of the root functions $z_{n,j}(\theta)$ has already been established in Theorem~\ref{thm4.3}. Since the entries of $\mathbf{F}(k, z, \theta)$ are analytic functions of $z$ and $\theta$, and the null vector $\mathbf{p}^{n,j}(\theta)$ of the analytic matrix $\mathbf{A}(z_{n,j}(\theta), \theta)$ can also be chosen to be analytic, the composition of these functions yields entries for $\mathbf{B}(\theta)$ that are analytic in $\theta$. Consequently, $\det \mathbf{B}(\theta)$ is an analytic function for $\theta \in V$.
    
    For the special case $\theta=\pi$, the wedge degenerates to a straight line. This geometric simplification induces a block-diagonal structure in the matrix $\mathbf{A}(z, \pi)$, which allows for a direct calculation of the null vectors and the matrix $\mathbf{B}(\pi)$. As detailed in Appendix~\ref{app:invertibility}, this calculation confirms that $\det(\mathbf{B}(\pi))$ is non-zero.
\end{proof}
\end{lemma}

\begin{lemma}\label{lem4.7}
    The matrix $\mathbf{B}(\theta)$ is invertible for all $\theta \in (0, 2\pi)$ except for at most a countable set of values. Moreover, this set of exceptional angles has no limit points within $(0, 2\pi)$.
\begin{proof}
    By the identity theorem, this result is a direct consequence of the properties of analytic functions established in Lemma \ref{lem4.6}.
    %We have shown that the function $\det \mathbf{B}(\theta)$ is analytic on a domain containing the real interval $(0, 2\pi)$ and $\det \mathbf{B}(\theta) \neq 0$. By the identity theorem for analytic functions, the zeros of a non-trivial analytic function must be isolated. This implies that the set of angles $\theta \in (0, 2\pi)$ for which $\det \mathbf{B}(\theta)=0$ is at most countable and have no limit points within the open interval $(0, 2\pi)$. Therefore, the matrix $\mathbf{B}(\theta)$ is invertible for all but a countable set of exceptional values.
\end{proof}
\end{lemma}

Now we obtain the final result of this section, which shows that for most opening angles $\theta$, the density function $\bphi$ can be constructed to match arbitrary polynomial boundary data up to order $N-1$.
\begin{theorem}\label{thm4.8}
    Suppose that $N \geq 2$ is an integer. Then for each $\theta \in(0,2\pi)$ except for countably many values, there exist $z_{n,j}\in\C$, $\mathbf{p}^{n,j}\in\C^4$, $n=1,2,\ldots,N,\ j=1,2,3,4$, such that the following holds: Suppose $\boldsymbol{\beta}_k = \left(\beta_{k,1}, \beta_{k,2}, \beta_{k,3}, \beta_{k,4}\right)^\top \in \C^4, k=0,1, \ldots, N-1$ and the boundary data $\boldsymbol{f}(t)$ is given by
    \begin{align*}
    \tilde{\boldsymbol{f}}(t) = 
        \begin{bmatrix}
            f_1(t) & f_2(t) & f_1(-t) & f_2(-t)
        \end{bmatrix}^\top
    = \sum_{k=0}^{N-1} \boldsymbol{\beta}_k t^k + O(t^{N}), \ 0<t<1.
    \end{align*}
    Then there exist unique vectors $\boldsymbol{\alpha}_n = \left(\alpha_{n,1}, \alpha_{n,2}, \alpha_{n,3}, \alpha_{n,4}\right)^\top \in \C^4, n = 1, \ldots, N$, such that if the density $\bphi$ is defined by
    \begin{align*}
    \tilde{\bphi}(t) =
        \begin{bmatrix}
            \phi_1(t) & \phi_2(t) & \phi_1(-t) & \phi_2(-t)
        \end{bmatrix}^\top=
        \sum_{n=1}^N\sum_{j=1}^4 \alpha_{n,j}
             \mathbf{p}^{n,j}\cdot t^{z_{n,j}},\ 0<t<1,
    \end{align*}
    then $\bphi$ satisfies the following boundary integral equation on the wedge:
    $$\boldsymbol{f}(\bx(t)) = -\frac{1}{2}\bphi(\bx(t)) + \mathcal{D}[\bphi](t),\ -1\leq t\leq 1$$
    with an error of order $O\left(|t|^{N}\right)$.
\end{theorem}

\subsection{Hybrid basis: limitations on sharp corners}
It seems natural to use the functions obtained from Theorem \ref{thm4.5} as a basis for the numerical solutions of BIEs on cornered domains. However, we will show the analytic basis functions are unsuitable for direct use in numerical algorithms. For that purpose, we design a hybrid semi-analytic boundary element method (called HybridBEM) on an equilateral triangle to explore the possibilities of using corner singularities as basis functions. Specifically, on each corner panel $\Gamma_c$ the density is expanded in the first $M$ singular powers $\{t^{z_k}\}_{k=1}^M$ with graded collocation, while the remaining boundary $\Gamma_s$ uses $p$-th order Legendre expansion. Panel interactions are assembled by factoring the kernel-density product into a ``singular weight $\times$ smooth factor'' and evaluating the resulting panel moments semi-analytically. Particularly, cross-corner terms follow the explicit mapping series in subsection~\ref{subsec43}. 

\begin{figure}[!htbp]
  \centering
  \begin{minipage}[b]{.46\linewidth}
    \centering
    \includegraphics[width=\linewidth]{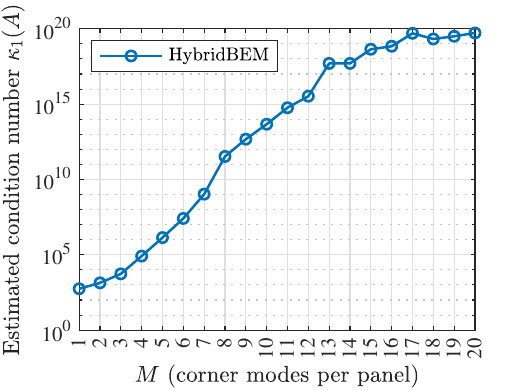}\\
    \footnotesize(a) Condition number vs.\ corner modes
  \end{minipage}\hfill
  \begin{minipage}[b]{.46\linewidth}
    \centering
    \includegraphics[width=\linewidth]{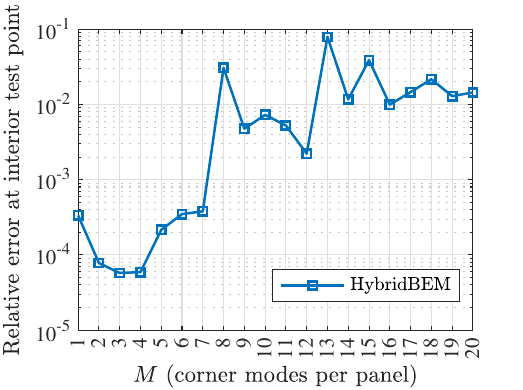}\\
    \footnotesize(b) Relative error vs.\ corner modes
  \end{minipage}
  \caption{HybridBEM on an equilateral triangle. (a) The condition number of the discretized system versus the number of corner modes $M$. (b) The relative error at the test point versus $M$.}
  \label{fig:hybrid-tri}
\end{figure}

Figure~\ref{fig:hybrid-tri} tests the performance as we vary the number of corner modes $M$. As $M$ increases, the number of unknowns grows linearly, yet the conditioning deteriorates: $\kappa_1(A)$ rises from $\mathcal{O}(10^2)$ at $M=1$ to $\mathcal{O}(10^{19})$ at $M=20$. The accuracy does not benefit from the extra corner modes, and even in the best case, the relative error levels off near $O(10^{-5})$, which is already worse than the uniform-panel baseline $O(10^{-6})$. This breakdown is algebraic rather than quadrature-driven. In particular, sampling the functions $\{t^{z_k}\}$ on a finite corner panel produces columns with a Vandermonde-like structure. When several real parts $\Re z_k$ are close, or when a conjugate pair has a very small imaginary part $\Im z$, these columns become nearly linearly dependent, and the local Gram matrices are severely ill-conditioned, which explains the poor stability and lack of accuracy improvement. To overcome these issues, we introduce a singularity guided Nystr\"om method in the next section.

% -------------------------------
% Section 5: Numerical algorithm
% -------------------------------
\section{Singularity guided Nystr\"om method}\label{sec5}
%% Section 5 给出了SGN算法及其误差分析，最后给出了3个数值实验
% 1. SGN 原理：使用奇异性指标指导生成相应的网格；
% 2. 误差分析：根据wellposedness做相应的分析；
% 3. 数值实验：基准实验+凹角实验+奇异幂函数实验。后续在基准实验中考虑加入RCIP做对比！！！
Solving the elastic integral equation~\eqref{eq:bie_original} in domains with corners requires specialized techniques to handle both kernel and geometric singularities effectively. Based on the analytical results above, we propose a singularity guided high-order Nystr\"om (SGN) scheme. Specifically, corner singular powers are extracted from the transcendental equation $H(z,\theta) = 0$ via an analytic root-finding algorithm in the complex plane~\cite{kravanja2000computing} and then used to guide panelization and accuracy control, while the unknown density is still discretized at Gauss-Legendre (GL) nodes. Near or self interactions are evaluated by kernel splitting and a forward-stable Legendre recurrence~\cite{ma1996generalized,yao2024robust, martinsson2007interpolation}, effectively addressing Cauchy-type singularities. This yields a square and well-conditioned linear system, so a dense linear solver suffices to obtain stable results. Details of the SGN scheme are given as follows.

\subsection{Algorithm design}
We parameterize each boundary arc by arclength $t\in[0, L]$ originating from the corner at $t=0$. For the opening angle $\theta$, let $z_1, \ldots, z_n$ be the leading singularity exponents with $\Re z > 0$, computed via the contour integration method~\cite{kravanja2000computing}. To resolve the corner singularities, we construct an adaptive panel mesh by recursive bisection. On each panel $I$, we define the multi-exponent Legendre-tail indicator at degree $p$ by
\begin{equation}\label{eq:legendre_indicator}
    \eta(I;p):=\max_{1\leq k\leq n}\Bigg(\sum_{j=p}^{2p-1}\big|\alpha_j(z_k;I)\big|^2\Bigg)^{1/2},
\end{equation}
where $\alpha_j(z_k;I)$ is the $j$-th Legendre coefficient of the singular function $t^{z_k}$ restricted to $I$. A panel is accepted if \begin{equation}\label{acceptrule}
    \eta(I;p)\le\varepsilon_{\mathrm{pan}},
\end{equation} where $\varepsilon_{\mathrm{pan}}$ is a given error tolerance. Otherwise, it is bisected. Panels away from corners are kept uniform with the same degree $p$. In this way, the mesh clusters automatically wherever the computed singular family $\{t^{z_k}\}$ demands resolution. 

Now we consider the quadrature rules for the boundary integral:
\begin{equation}\label{eq:double_layer}
    \calK[\bphi](\bx) = \int \mathbf{D}\left(\bx, \by\right)
    \bphi(\by) \mathrm{d}s_{\by},
\end{equation}
where $\mathbf{D}\left(\bx, \by\right)$ is given in~\eqref{eq:fundamental_solution_traction}. Far-field interactions are computed via the standard $p$-point Gauss-Legendre rule. For near-field and self interactions, we utilize the kernel splitting technique:
$$\mathbf D(t,s) = \frac{1}{t-s} \mathbf D_1(t,s) + \mathbf D_2(t,s),$$
where $\mathbf D_1, \mathbf D_2$ are analytic on $[-1,1]^2$:
$$\begin{aligned}
    \mathbf{D}_1(t, s) &= 
    \frac{\mu}{\lambda+2\mu}\frac{(t-s)\,(\btau(\bx)\cdot\br)}{2\pi r^2}\mathbf{L},\\
    \mathbf{D}_2(t, s) &= \frac{\mu}{\lambda+2\mu}\frac{(\btau(\by)-\btau(\bx))\cdot\br}{2\pi r^2}\mathbf{L}+\frac{\mu}{\lambda+2\mu}\frac{\bn(\by)\cdot\br}{2\pi r^2}\mathbf{I}
    +\frac{\lambda+\mu}{\lambda+2\mu}\frac{(\bn(\by)\cdot\br)\br\otimes\br}{\pi r^4}.
\end{aligned}$$
Moreover, let $J(s) = |\gamma'(s)|$ be the Jacobian of the panel parameterization at the variable $s$, i.e., $d s_{\by}=J(s)\,\mathrm{d}s$. Interpolating $\mathbf D_j\,\bphi\,J$, $j=1,2$ by Legendre polynomials
$$(\mathbf D_1\,\bphi\,J)(t,s)\approx\sum_{j=0}^{p-1}\boldsymbol d_j(t)\,P_j(s),\qquad
(\mathbf D_2\,\bphi\,J)(t,s)\approx\sum_{j=0}^{p-1}\boldsymbol f_j(t)\,P_j(s),$$
leads to panel moments
$$C_j(t):=\text{p.v.}\int_{-1}^{1}\frac{P_j(s)}{t-s}\,\mathrm{d}s,\qquad    I_j:=\int_{-1}^{1}P_j(s)\,\mathrm{d}s .$$
Since the Gauss-Legendre nodes never include the endpoints, the Cauchy moments $C_j$ can be evaluated via the three-term recurrence
$$
\bigg\{
\begin{aligned}
&C_0(t)=-\ln\Big|\frac{t-1}{t+1}\Big|,\quad
C_1(t)=-2+t\,C_0(t), \\
&j\,C_j=(2j-1)\,t\,C_{j-1}-(j-1)\,C_{j-2}, \qquad j\ge2,
\end{aligned}
\bigg.
$$
for any real $t\neq\pm1$. Thus, near-field or self interactions can be approximated by $\sum_{j=0}^{p-1}\big[C_j(t)$ $\boldsymbol d_j(t)+I_j\,\boldsymbol f_j(t)\big]$. Enforcing the equation at all Gauss-Legendre nodes produces a square Nystr\"om system, which is then solved by a direct solver.

\subsection{Error analysis}
% 接下来对提出的SGN算法做误差分析～
Having detailed the algorithmic framework, we now analyze the error in the weighted Sobolev space $X := H_\nu^s(\Gamma)^2$ with $s\in(1/2,3/2)$ and $s<\nu<s+1/2$. Assuming $s-\nu\notin\{\Re z_{n,j}\} \cup \mathbb{Z}$,
Theorem~\ref{thm:wellposedness} guarantees that the interior double-layer operator $$\mathcal A:=\Big(-\frac{\theta(\bx)}{2\pi}\mathcal I+\calK\Big):X\to X$$
is an isomorphism with $C_{\mathcal A}:=\|\mathcal A^{-1}\|_{\mathcal L(X)}<\infty$, where $\|\cdot\|_{\mathcal L(X)}$ denotes the operator norm from $X$ to $X$. Throughout this subsection, $C>0$ denotes a generic constant independent of $h$ and $\varepsilon_{\rm pan}$. 
\begin{lemma}[Decay of Legendre coefficients]\label{lem:legendre_tail}
Let $I=[0,h]$ and $f(t)=t^z$ with $\Re z>0$. The Legendre coefficients $\alpha_j(z;I)$ of $f$ on $I$ satisfy
\begin{equation}\label{eq:legendre_decay}
|\alpha_j(z;I)| \le C(z)\,h^{\Re z}\,j^{-(2\Re z+1)}, \quad j\ge 1.
\end{equation}
\begin{proof}
Under the change of variables $t=\tfrac{h}{2}(1+x)$, we obtain $g(x)=(\frac{h}{2})^z(1+x)^z, x\in[-1,1]$. By the asymptotic theory for Legendre polynomials~\cite[Thm.~3.1]{Sidi2011}, the $j$-th coefficient of $(1+x)^z$ decays as $O(j^{-(2\Re z+1)})$. Combining the scaling factor $h^z$, we obtain the result.
\end{proof}
\end{lemma}

\begin{lemma}[Weighted best approximation]\label{lem:best}
Let $\Pi_h$ be the piecewise polynomial projection on the SGN mesh $\mathcal T_h$, and $V_h$ the space of piecewise polynomials of degree $p$ on $\mathcal T_h$. For $\phi\in X$ with singular expansion $\phi(t)=\sum_k c_k\,t^{z_k}+\phi_{\rm reg}(t)$ with $\Re z_k>0$ and $\phi_{\rm reg}$ analytic near corners, there exist constants $C,\kappa>0$ such that
$$\inf_{v_h\in V_h}\|\phi-v_h\|_{X}
\ \le\ C\Big(\varepsilon_{\rm pan}+e^{-\kappa p}\Big),$$
where $\varepsilon_{\rm pan}$ is a given tolerance that satisfies inequality \eqref{acceptrule}.
\begin{proof}
    For the analytic part $\phi_{\rm reg}$, since the weight exponent satisfies $\nu-s > 0$, the weighted norm is bounded by the standard Sobolev norm on each panel. Thus, standard approximation theory~\cite{CanutoSpectral} implies exponential convergence:
    $$
    \|\phi_{\rm reg}-\Pi_h\phi_{\rm reg}\|_{X} \le C\,e^{-\kappa p}.
    $$
    For the singular part $t^{z_k}$, consider a corner panel $I=[0,h]$. Mapping $I$ to $[-1,1]$, the projection error is given by $e_p := t^{z_k} - \Pi_h t^{z_k} = \sum_{j \geq p} \alpha_j P_j$.
    The mesh construction ensures $\eta(I;p) \le \varepsilon_{\mathrm{pan}}$. On the interval $[-1,1]$, the $H^1$-norm of the error is bounded by the weighted $\ell^2$-norm of the coefficients~\cite{CanutoSpectral}:
    $$
    \|e_p\|_{H^1(-1,1)}^2
    \lesssim \sum_{j\ge p} j^2 |\alpha_j|^2.
    $$
    Thus, by Lemma~\ref{lem:legendre_tail} together with the detailed asymptotics in~\cite{Sidi2011}, the infinite tail satisfies 
    $$
    \sum_{j\ge p} j^2 |\alpha_j|^2
    \le C\,p^2 \sum_{j=p}^{2p-1} |\alpha_j|^2
    \le C\,p^2 \varepsilon_{\rm pan}^2.
    $$
    This yields $\|e_p\|_{H^1(-1,1)} \le C\,p\,\varepsilon_{\rm pan}$.
    Applying a standard scaling argument for the weighted norm on $I$~\cite{Babuska1994pFEM}, we obtain
    $$
    \|t^{z_k}-\Pi_h t^{z_k}\|_{X(I)} \le C\,h^{\nu-s+\frac12}\,\|e_p\|_{H^1(-1,1)}
    \le C\,h^{\nu-s+\frac12}\,p\,\varepsilon_{\rm pan}.
    $$
    Summing over all corner panels and over the finitely many exponents $z_k$, and using $\nu-s > 0$ so that \(\sum_{I\in\mathcal T_h^{\rm corner}} h_I^{2(\nu-s)+1}\le C\), yields
    \[
    \|\phi_{\rm sing}-\Pi_h\phi_{\rm sing}\|_X \le C\,p\,\varepsilon_{\rm pan},
    \]
    where $p$ can be absorbed into the constant $C$. Combining the estimates for $\phi_{\rm reg}$ and $\phi_{\rm sing}$ and minimizing over $v_h\in V_h$ proves the claim.
\end{proof}
\end{lemma}

\begin{lemma}[Operator consistency]\label{lem:consistency} 
Let $\mathcal A_h:X\to X$ be the SGN operator and define $\Delta_h:=\mathcal A-\mathcal A_h$ with $\delta_h:=\|\Delta_h\|_{\mathcal L(X)}$. Then there exist $C_0,\kappa>0$ and a function $c_{\rm far}(p)$ with at most polynomial growth on $p$ such that
$$\|\mathcal A-\mathcal A_h\|_{\mathcal L(X)}
\ \le\ C\Big(\varepsilon_{\rm pan}+c_{\rm far}e^{-\kappa p}\Big).$$
\begin{proof}
    For near and self interactions, the local error arises solely from replacing the density in $(\mathbf D_k \bphi J)$ with its $p$-th order Legendre interpolant. The analyticity of $\mathbf D_k$ and the tail condition $\eta(I;p)\le\varepsilon_{\rm pan}$ ensure an interpolation error of $O(\varepsilon_{\rm pan})$. By the boundedness of singular integral operators on $X$, this translates to an $O(\varepsilon_{\rm pan})$ operator error. For far-field interactions, the kernel is analytic. The $p$-point Gauss--Legendre quadrature introduces an error of order $O(e^{-\kappa p})$. Summing over interactions yields the claimed bound.
\end{proof} 
\end{lemma}

Combining the stability of invertible operators, we derive the main error bound.
\begin{theorem}[Nystr\"om error bound]\label{thm:nystrom_error}
Assume $\delta_h C_{\mathcal A}<1$.
Then $\mathcal A_h$ is invertible on $X$ and its inverse satisfies
$\|\mathcal A_h^{-1}\|_{\mathcal L(X)}\le C_{\mathcal A}/(1-\delta_h C_{\mathcal A})$.
Let $\phi_h\in V_h$ be the SGN solution. It holds
\begin{equation}\label{errorbound1}
 \|\phi-\phi_h\|_{X}
\ \le\ \frac{C\,C_{\mathcal A}}{1-\delta_h C_{\mathcal A}}\,
\Big(\inf_{v_h\in V_h}\|\phi-v_h\|_{X}+\|\Delta_h\phi\|_{X}\Big)   
%\end{equation}
%Consequently, it holds
%\begin{equation}\label{errorbound2}
%    \|\phi-\phi_h\|_{X}
\ \le\ \frac{C\,C_{\mathcal A}}{1-\delta_h C_{\mathcal A}}\,
\Big(\varepsilon_{\rm pan}+c_{\rm far}e^{-\kappa p}\Big).
\end{equation}
\begin{proof}
    First, $\mathcal A_h=\mathcal A-\Delta_h=\mathcal A\big(I-\mathcal A^{-1}\Delta_h\big)$ and
    $\|\mathcal A^{-1}\Delta_h\|\le C_{\mathcal A}\delta_h<1$, hence
    $\mathcal A_h^{-1}=\big(I-\mathcal A^{-1}\Delta_h\big)^{-1}\mathcal A^{-1}$ with $\|\mathcal A_h^{-1}\|\le C_{\mathcal A}/(1-\delta_h C_{\mathcal A})$. For any $v_h\in V_h$, we have
    $$\mathcal A_h(\phi_h-v_h)
    =\mathcal A\phi-\mathcal A_h v_h
    =\mathcal A(\phi-v_h)+(\mathcal A-\mathcal A_h)v_h
    =\mathcal A(\phi-v_h)+\Delta_h v_h.$$
    Applying $\mathcal A_h^{-1}$ yields
    $$\|\phi_h-v_h\|
    \le \|\mathcal A_h^{-1}\|\big(\|\mathcal A\|\,\|\phi-v_h\|+\|\Delta_h v_h\|\big).$$
    Moreover, it holds that
    $$\|\Delta_h v_h\|\le \|\Delta_h\phi\|+\|\Delta_h\|\,\|\phi-v_h\|
    =\|\Delta_h\phi\|+\delta_h\,\|\phi-v_h\|.$$
    Combining the above two inequalities yields
    $$\|\phi_h - v_h\| \leq \|\mathcal A_h^{-1}\| \Big( (\|\mathcal A\| + \delta_h)\,\|\phi-v_h\| + \|\Delta_h\phi\| \Big).$$
    Using $\|\mathcal A_h^{-1}\|\le C_{\mathcal A}/(1-\delta_h C_{\mathcal A})$, we can obtain
    $$\|\phi-\phi_h\| \le \|\phi-v_h\| + \|\phi_h-v_h\|
    \le \frac{C\,C_{\mathcal A}}{1-\delta_h C_{\mathcal A}}
    \big(\|\phi-v_h\|+\|\Delta_h\phi\|\big).$$
    Minimizing over $v_h\in V_h$  and using  Lemmas~\ref{lem:best} and \ref{lem:consistency} complete the proof of \eqref{errorbound1}. 
\end{proof}
\end{theorem}

\begin{corollary}\label{cor:u_field}
For any compact set $\Omega_0\Subset\Omega$, the interior double-layer
potential operator~\eqref{eq:DLP} maps $X$ boundedly into $L^2(\Omega_0)^2$. Therefore,
\[
\|\boldsymbol{u}-\boldsymbol{u}_h\|_{L^2(\Omega_0)^2}
\ \le\ \frac{C\,C_{\mathcal A}}{1-\delta_h C_{\mathcal A}}\,
\Big(\varepsilon_{\rm pan}+c_{\rm far}e^{-\kappa p}\Big).
\]
\begin{proof}
Let $D$ denote the DLP with kernel
$K(\bx,\by):=(\mathcal T_{\by}\mathbf G(\bx,\by))^{\top}$, which is smooth and bounded on $\Omega_0\times\Gamma$ since $\Omega_0\Subset\Omega$ has a positive distance from $\Gamma$. Choose neighborhoods $U_k$ of the corner vertices and set $\Gamma_{\rm reg}:=\Gamma\setminus\bigcup_{k=1}^M U_k$.
Write
\[
D\boldsymbol\phi = \int_{\Gamma_{\rm reg}} K(\bx,\by)\boldsymbol\phi(\by)\,ds_{\by}
\;+\;\sum_{k=1}^M \int_{\Gamma\cap U_k} K(\bx,\by)\boldsymbol\phi(\by)\,ds_{\by}
=: D_{\rm reg}\boldsymbol\phi+\sum_{k=1}^M D_k\boldsymbol\phi .
\]

Applying the Cauchy--Schwarz inequality to $D_k$ yields
\[
|(D_k\boldsymbol\phi)(\bx)|
\le
\Big(\int_{\Gamma\cap U_k}|K(\bx,\by)|^2\,r_k(\by)^{-2(\nu-s)}\,ds_{\by}\Big)^{1/2}
\Big(\int_{\Gamma\cap U_k}r_k(\by)^{2(\nu-s)}|\boldsymbol\phi(\by)|^2\,ds_{\by}\Big)^{1/2}.
\]
Integrating over $\bx\in\Omega_0$ and using Fubini's theorem gives
\[
\|D_k\boldsymbol\phi\|_{L^2(\Omega_0)^2}^2
\le
\Big(\int_{\Gamma\cap U_k}r_k(\by)^{2(\nu-s)}|\boldsymbol\phi|^2\,ds_{\by}\Big)
\int_{\Gamma\cap U_k} r_k(\by)^{-2(\nu-s)}
\Big(\int_{\Omega_0}|K(\bx,\by)|^2\,d\bx\Big)\,ds_{\by}.
\]
Since $\sup\limits_{\by\in\Gamma}\int_{\Omega_0}|K(\bx,\by)|^2\,d\bx < \infty$ and $\nu-s<1/2$ implies $\int_{\Gamma\cap U_k} r_k(\by)^{-2(\nu-s)}ds_{\by}<\infty$, we obtain
\[
\|D_k\boldsymbol\phi\|_{L^2(\Omega_0)^2}
\le C\Big(\int_{\Gamma\cap U_k}r_k(\by)^{2(\nu-s)}|\boldsymbol\phi|^2\,ds_{\by}\Big)^{1/2}
\le C\|\boldsymbol\phi\|_{X}.
\]

For the remaining part, Cauchy--Schwarz inequality and boundedness of $K$ also yield
\[
|(D_{\rm reg}\boldsymbol\phi)(\bx)|
\le \Big(\int_{\Gamma_{\rm reg}}|K(\bx,\by)|^2\,ds_{\by}\Big)^{1/2}\,
\|\boldsymbol\phi\|_{L^2(\Gamma_{\rm reg})^2}
\le C\,\|\boldsymbol\phi\|_{L^2(\Gamma_{\rm reg})^2}.
\]
Hence, integrating over $\Omega_0$ gives 
\[
\|D_{\rm reg}\boldsymbol\phi\|_{L^2(\Omega_0)^2}\le C\|\boldsymbol\phi\|_{L^2(\Gamma_{\rm reg})^2}\le C\|\boldsymbol\phi\|_{X}.
\] 

Combining the above two estimates yields $\|D\boldsymbol\phi\|_{L^2(\Omega_0)^2}\le C\|\boldsymbol\phi\|_X$.
Applying this bound to $\boldsymbol\phi-\boldsymbol\phi_h$ and using
Theorem~\ref{thm:nystrom_error} complete the proof.
\end{proof}
\end{corollary}

\subsection{Numerical experiments}
% 这里所谓的accuracy plateau指的是：relative error无法再进一步降低时所对应的精度阈值。
We now test the performance of the SGN scheme through several numerical experiments. Unless stated otherwise, we set the Lam\'e parameters $(\lambda,\mu)=(1,2)$, polynomial degree $p=16$, and length factor $\sigma=0.1$. Boundary data are synthesized using a few exterior sources placed away from the domain. The relative error is measured at interior target points $\boldsymbol x_m$ by
\begin{equation}
    E=\Bigg(\frac{\sum_m\|\boldsymbol u_{\text{num}}(\boldsymbol x_m)-\boldsymbol u_{\text{ex}}(\boldsymbol x_m)\|^2}{\sum_m\|\boldsymbol u_{\text{ex}}(\boldsymbol x_m)\|^2}\Bigg)^{1/2},
\end{equation}
where $\boldsymbol u_{\text{num}}$ and $\boldsymbol u_{\text{ex}}$ denote the numerical and exact displacement fields, and $\|\cdot\|$ is the Euclidean norm in $\mathbb R^2$.  All the experiments are implemented in \textsc{Matlab} on a personal computer equipped with an ARM processor. 

\subsubsection{Singularity guided panelization}
\begin{figure}[!ht]
  \centering
  \begin{minipage}[b]{.33\linewidth}
    \centering
    \includegraphics[width=\linewidth]{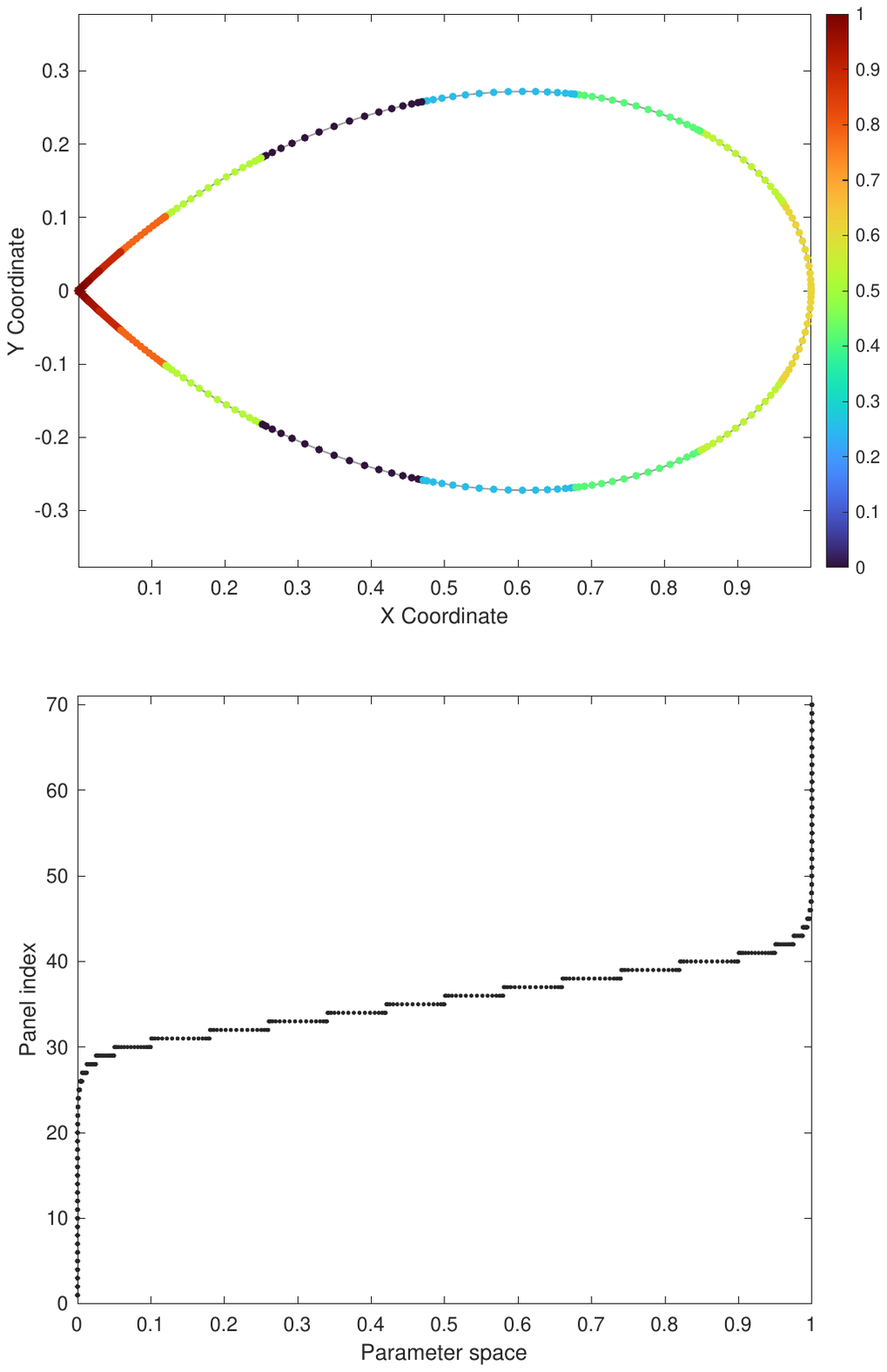}\\
    \footnotesize (a) Droplet
  \end{minipage}\hfill
  \begin{minipage}[b]{.33\linewidth}
    \centering
    \includegraphics[width=\linewidth]{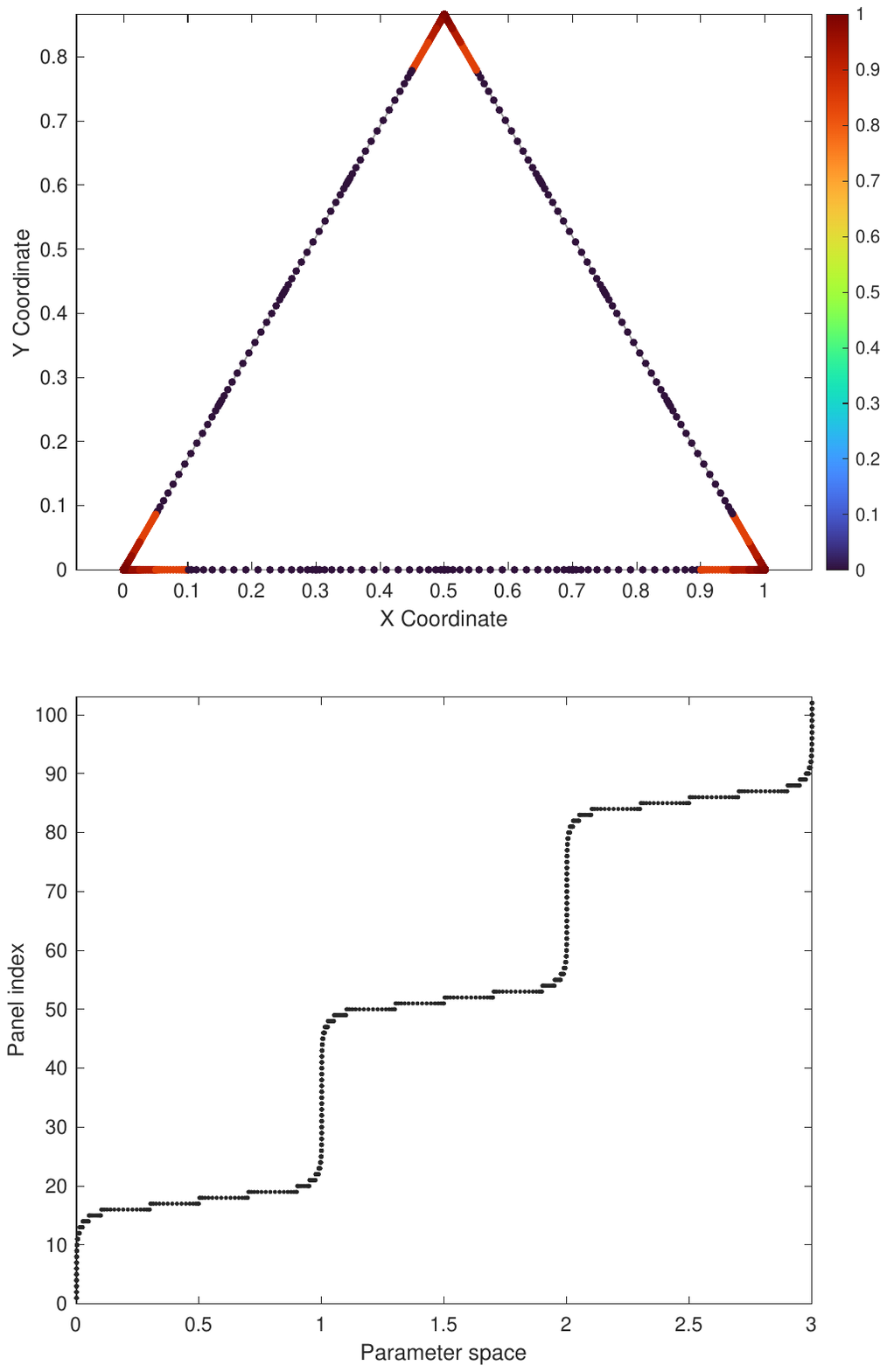}\\
    \footnotesize (b) Triangle
  \end{minipage}\hfill
  \begin{minipage}[b]{.33\linewidth}
    \centering
    \includegraphics[width=\linewidth]{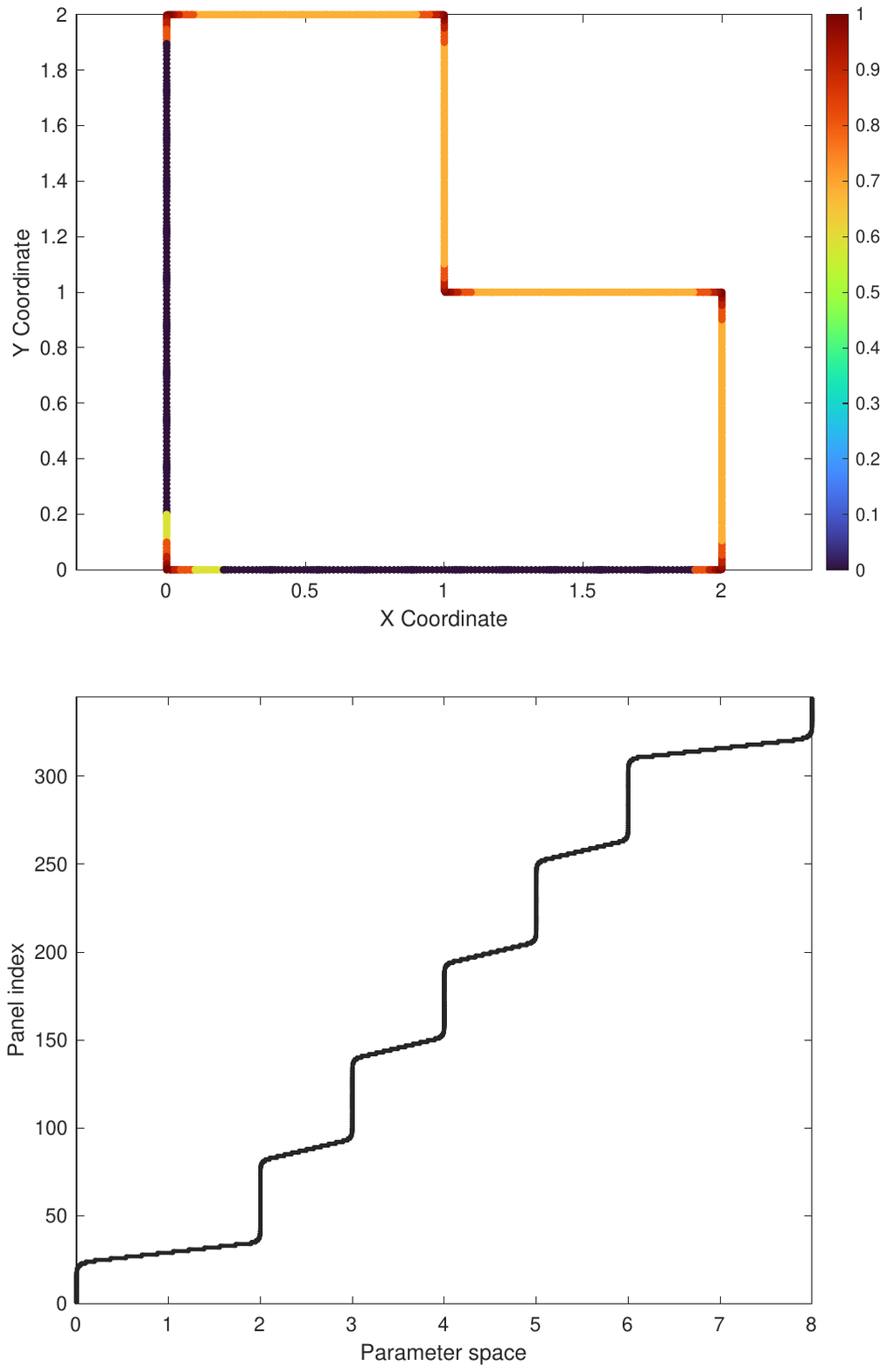}\\
    \footnotesize (c) L-shaped domain
  \end{minipage}
  \caption{Graded boundary discretizations and parameter-space panel distributions for the droplet, triangle, and L-shaped domains. The top row shows boundaries colored by relative panel length, and the bottom row shows the corresponding distributions of panels in the parameter variable $t$.}
  \label{fig:example1}
\end{figure}

In this example, we examine how SGN distributes panels near corners and how the error depends on the panel tolerance. We consider three representative geometries: a droplet, a triangle, and an L-shaped domain. For each geometry, we apply SGN varying the tolerance $\varepsilon_{\rm pan}$. Figure~\ref{fig:example1} illustrates the resulting graded boundary discretizations and their parameter-space distributions. 

\begin{figure}[!ht]
  \centering
   \begin{minipage}[b]{.45\linewidth}
    \centering
    \includegraphics[width=\linewidth]{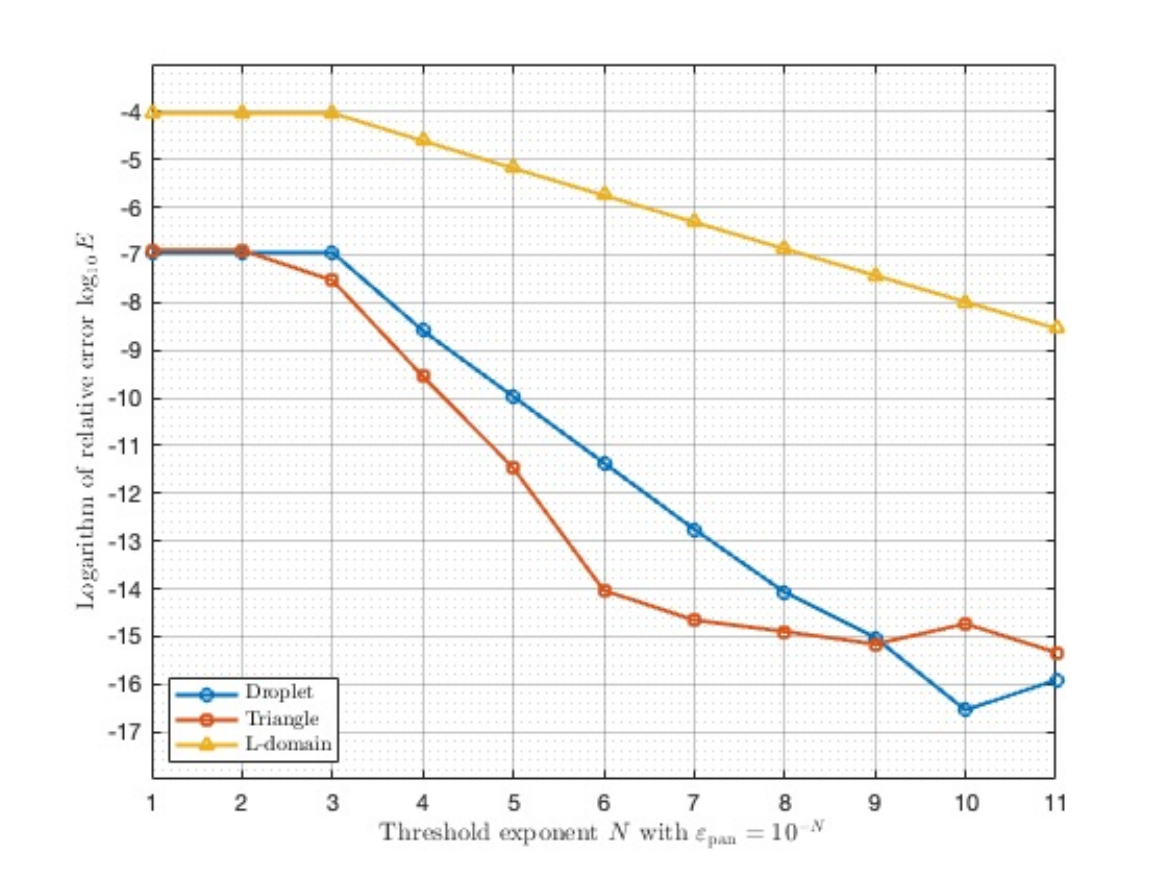}\\
    \footnotesize (a) Different geometries
  \end{minipage}
  \begin{minipage}[b]{.45\linewidth}
    \centering
    \includegraphics[width=\linewidth]{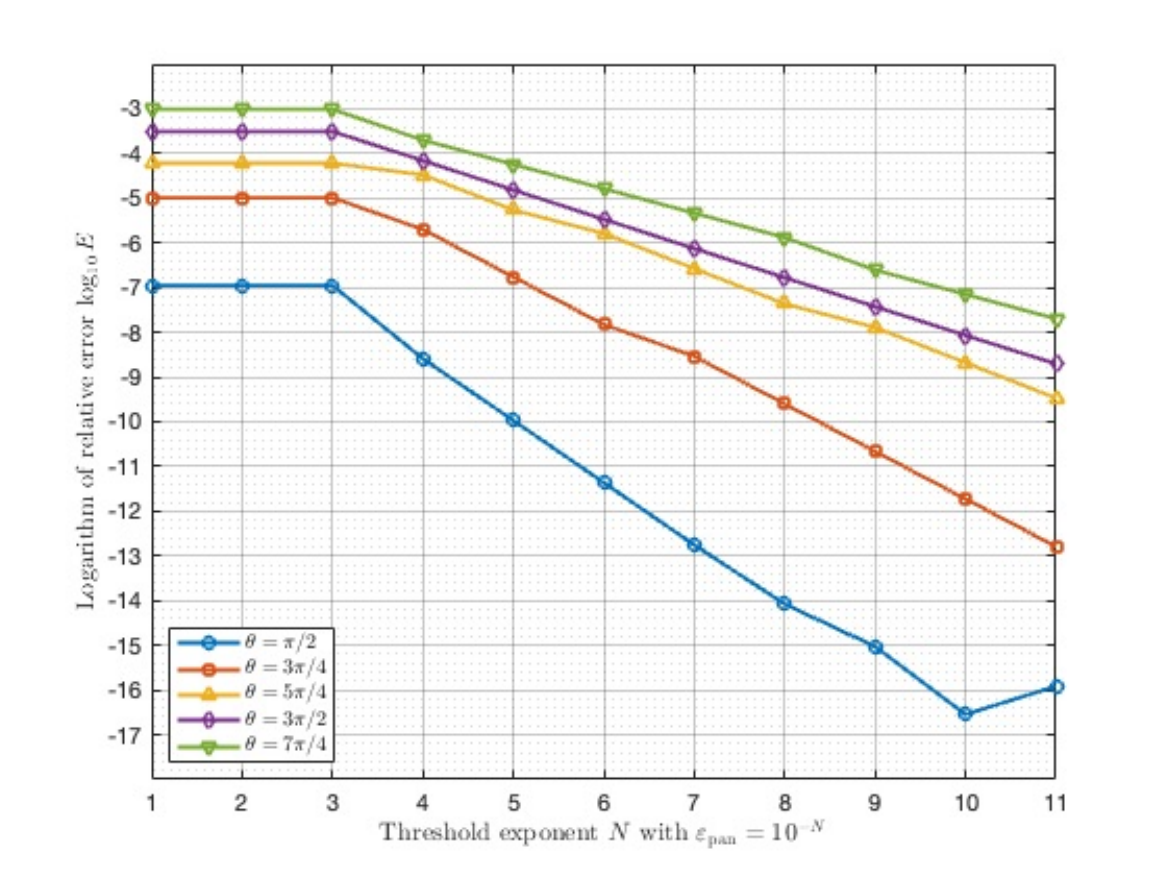}\\
    \footnotesize (b) Different corner angles
  \end{minipage}
  \caption{Relative error $\log_{10}\big(\|\boldsymbol u_{\text{num}}-\boldsymbol u_{\text{ex}}\|/\|\boldsymbol u_{\text{ex}}\|\big)$ versus the panel threshold $\varepsilon_{\text{pan}}$ in the SGN scheme. (a) Different geometries, including droplet, triangle, and L-shaped domains. (b) Single-corner droplets with varying opening angle $\theta$.}
  \label{fig:exp1-error-threshold}
\end{figure}

For the same three geometries, we vary the threshold exponent $N$ so that $\varepsilon_{\rm pan} = 10^{-N}$, and solve the scattering problem by SGN. The results are plotted in Figure~\ref{fig:exp1-error-threshold}(a), which shows the relative error versus $N$ for each shape in logarithmic scale. As $\varepsilon_{\text{pan}}$ decreases, the relative error decays almost linearly with the threshold in the regime where panel refinement dominates, which is in agreement with the error estimate given by Corollary~\ref{cor:u_field}.

\subsubsection{Comparison with a uniform mesh}
\begin{figure}[!ht]
  \centering
  %------------------- 第一行：UM -------------------%
  \begin{minipage}[b]{.34\linewidth}
    \centering
    \includegraphics[width=\linewidth]{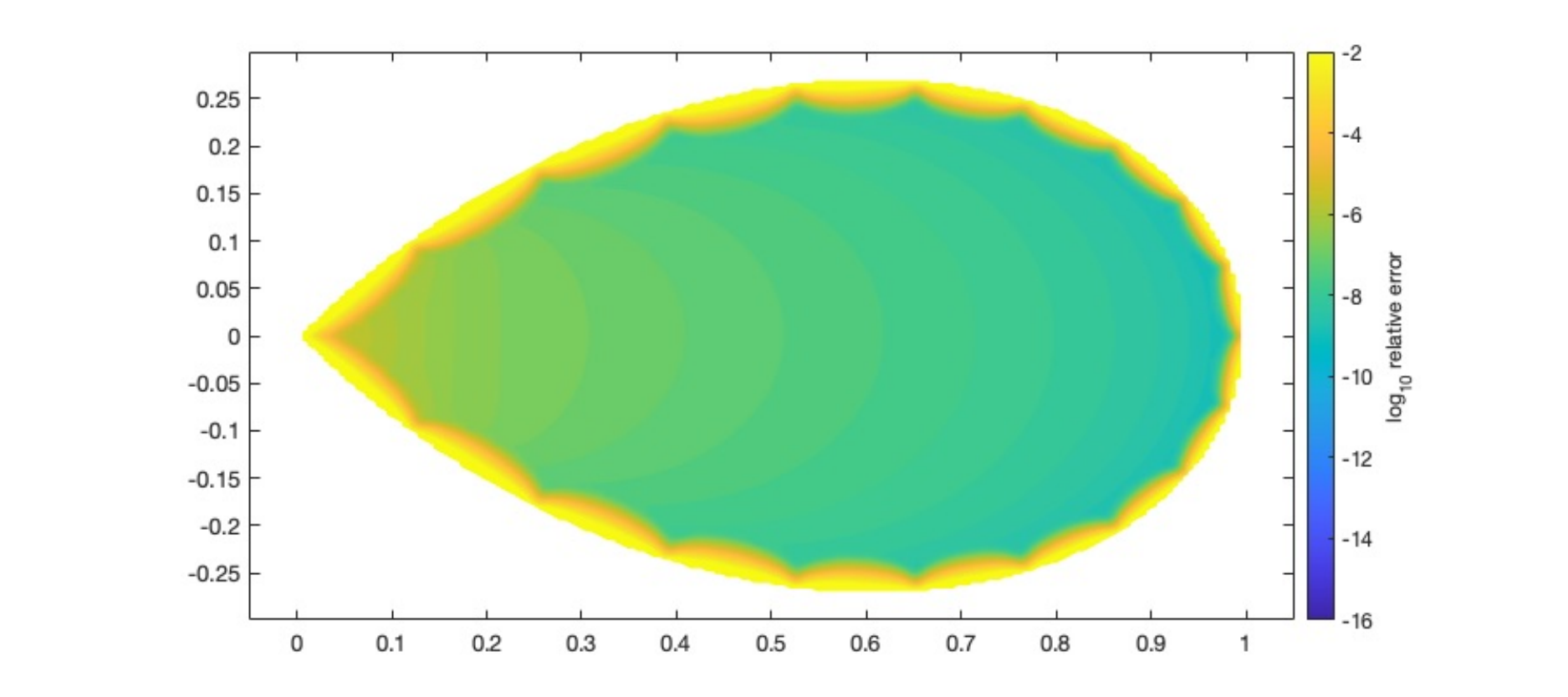}
  \end{minipage}\hspace{-0.03\textwidth}
  \begin{minipage}[b]{.34\linewidth}
    \centering
    \includegraphics[width=\linewidth]{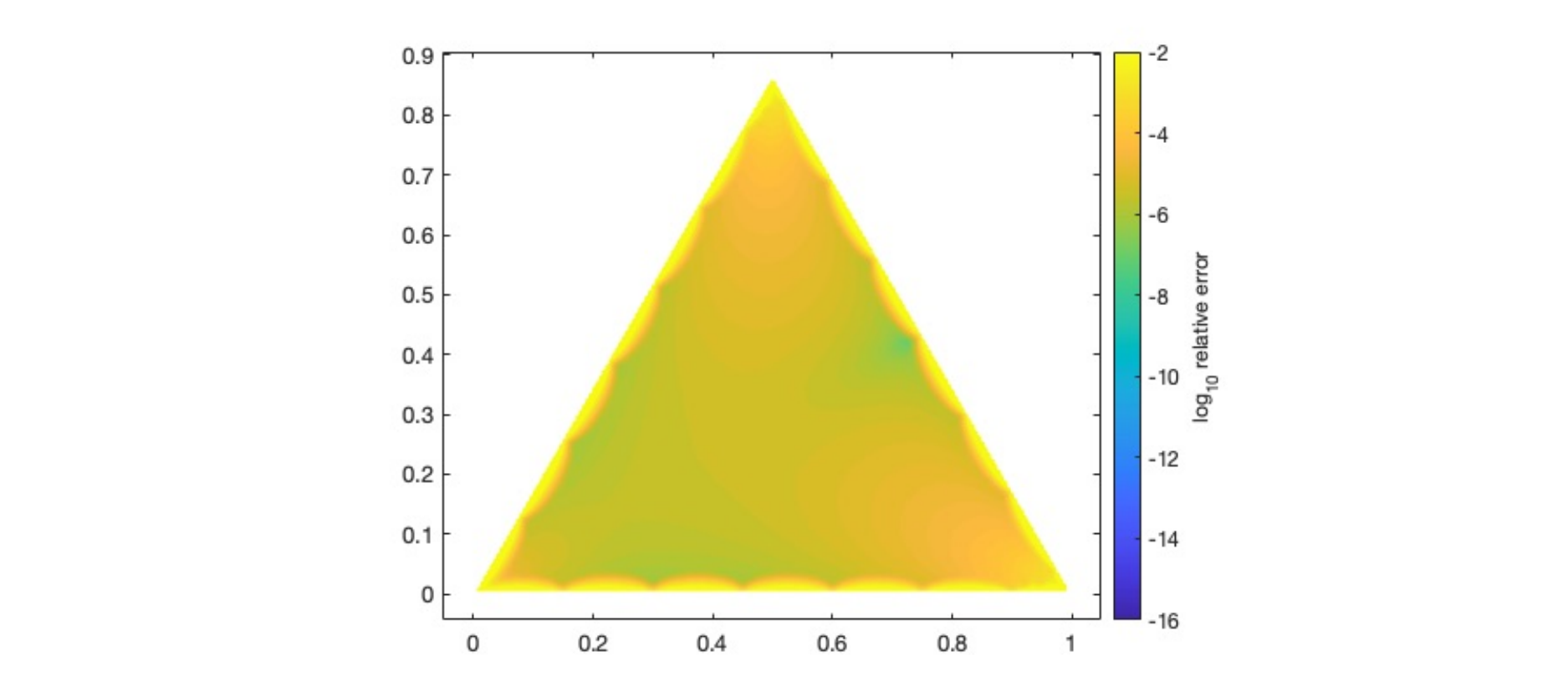}
  \end{minipage}\hspace{-0.08\textwidth}
  \begin{minipage}[b]{.34\linewidth}
    \centering
    \includegraphics[width=\linewidth]{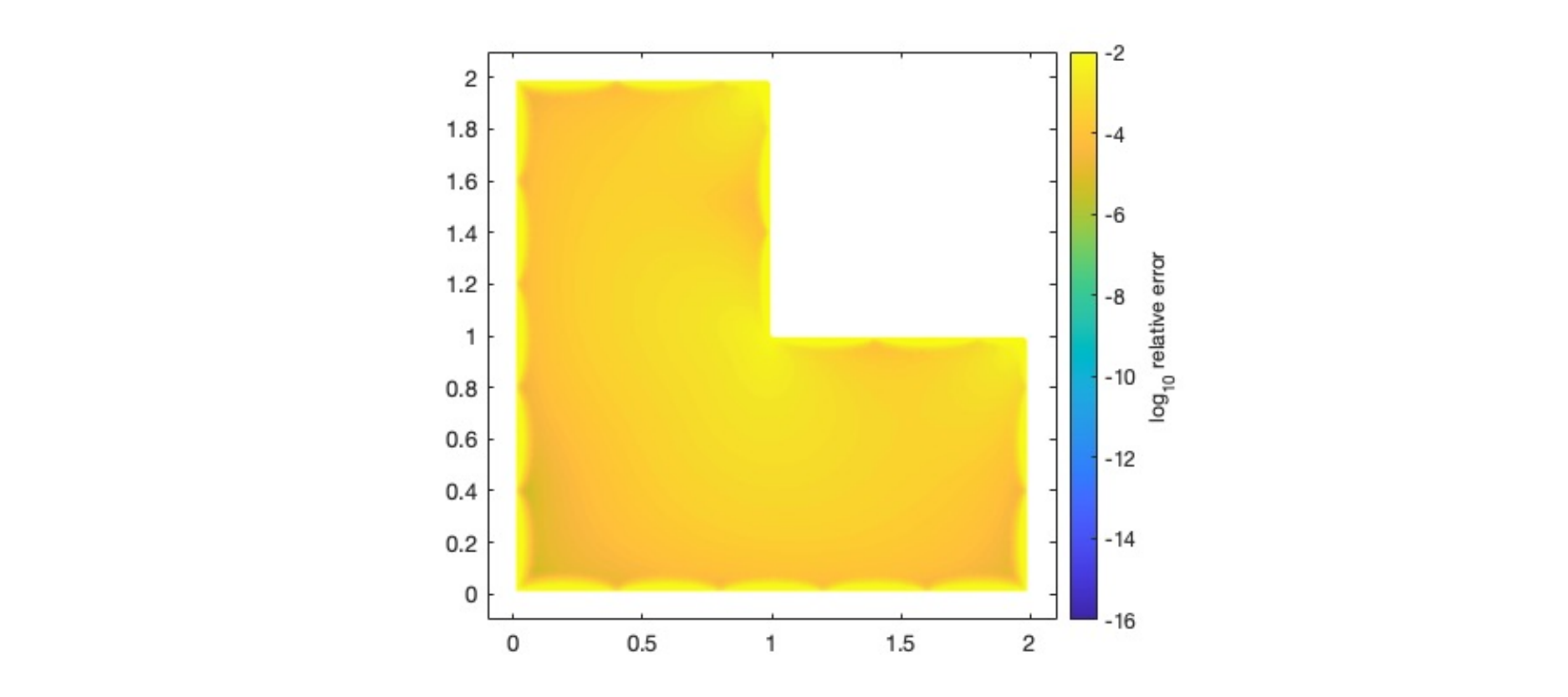}
  \end{minipage}
  %------------------- 第二行：SGN -------------------%
  \begin{minipage}[b]{.34\linewidth}
    \centering
    \includegraphics[width=\linewidth]{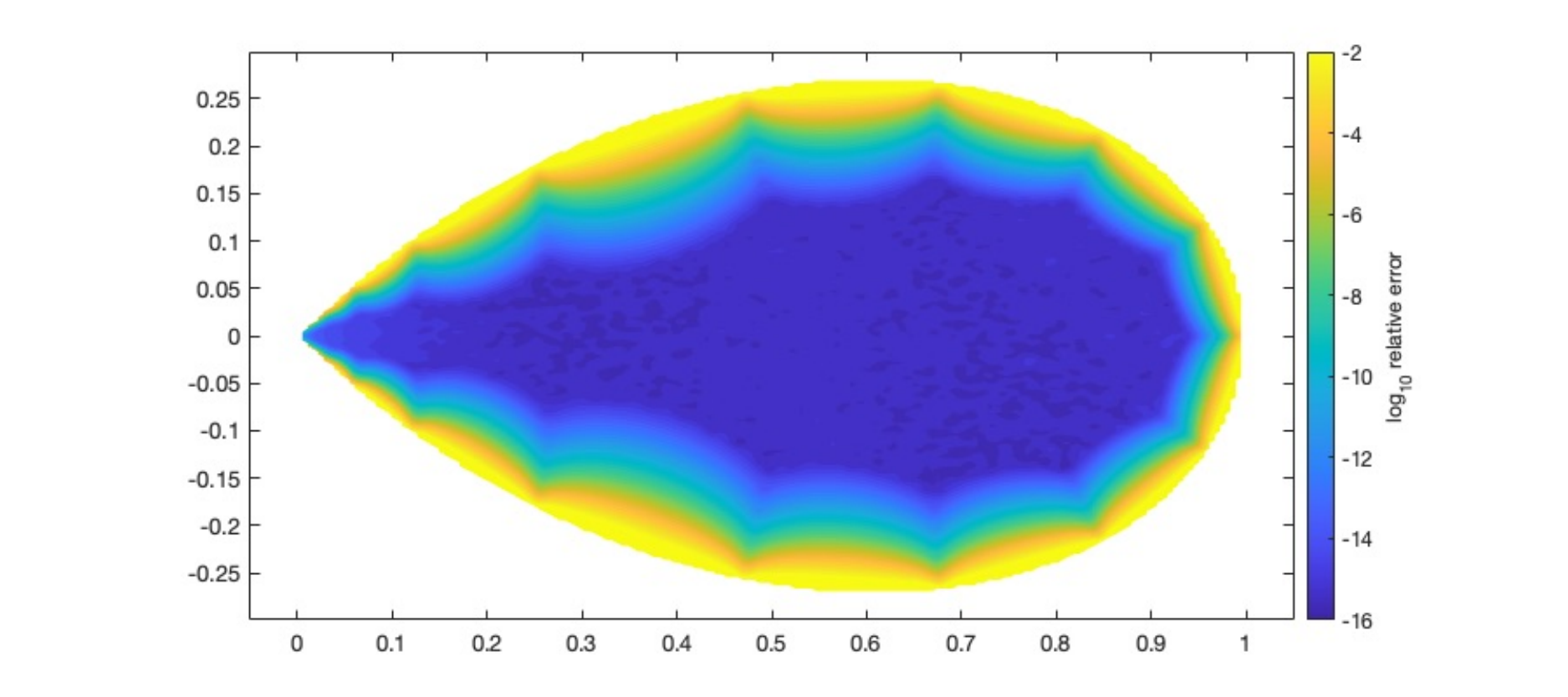}
  \end{minipage}\hspace{-0.03\textwidth}
  \begin{minipage}[b]{.34\linewidth}
    \centering
    \includegraphics[width=\linewidth]{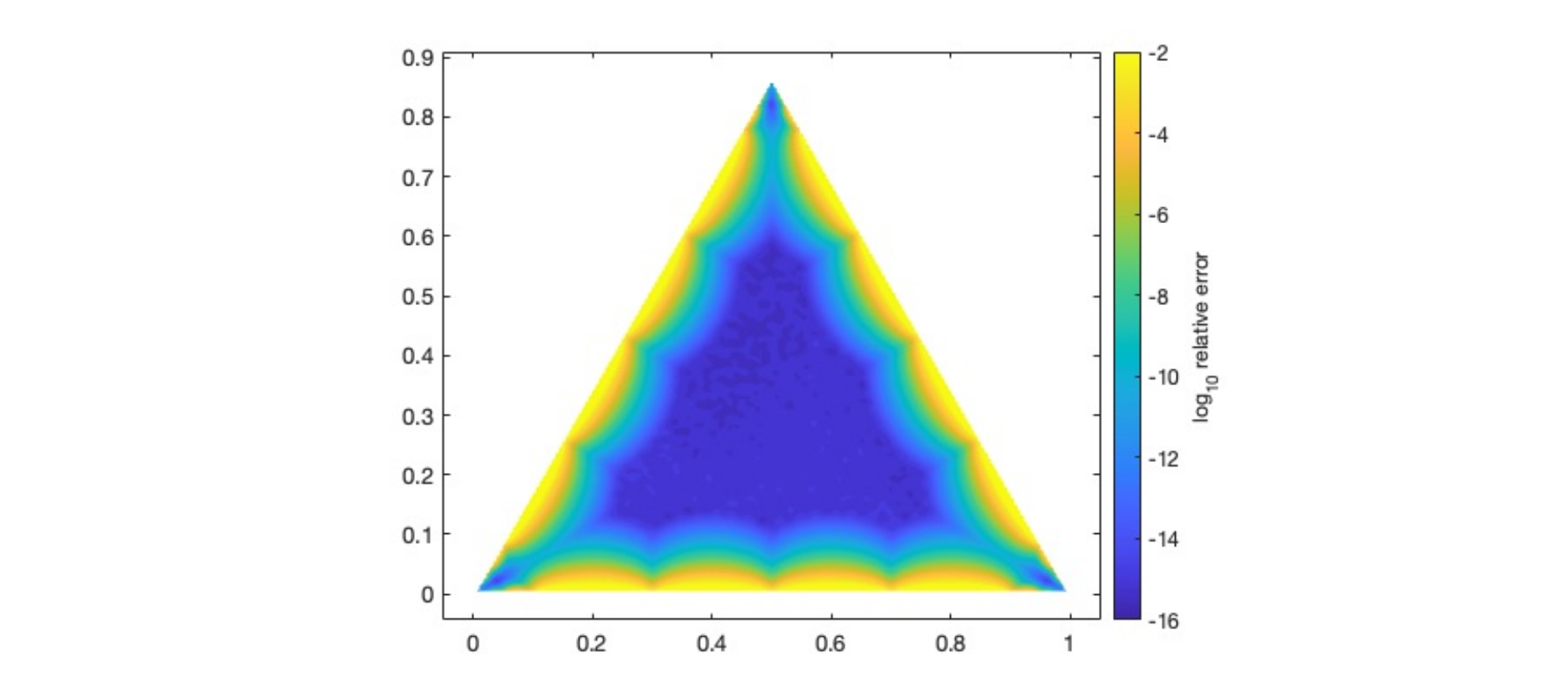}
  \end{minipage}\hspace{-0.08\textwidth}
  \begin{minipage}[b]{.34\linewidth}
    \centering
    \includegraphics[width=\linewidth]{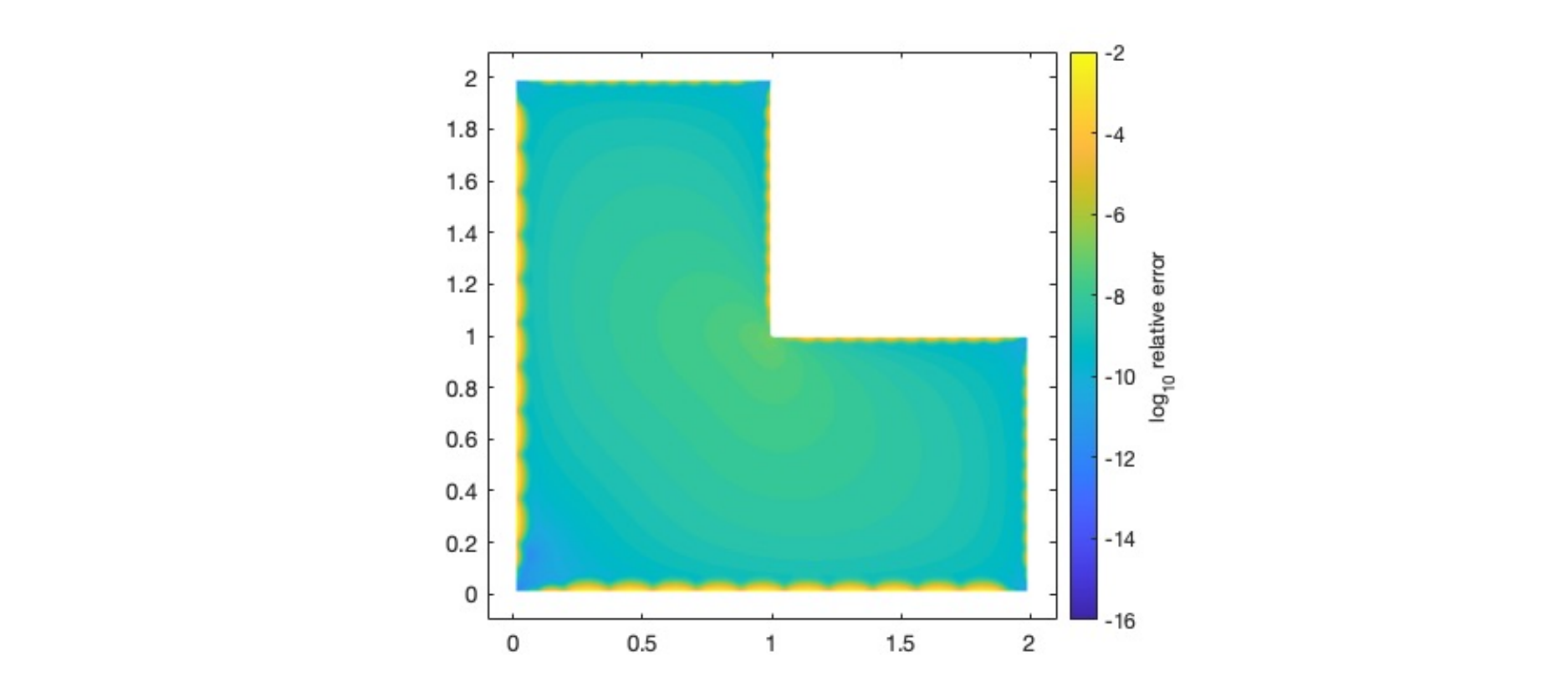}
  \end{minipage}
  \caption{Scattering fields for three geometries. The top and bottom rows represent UM and SGN respectively.}
  \label{fig:experiment2_wholefield}
\end{figure}

Guided by the first example, we choose $\varepsilon_{\rm pan}$ such that the accuracy stops increasing for each geometry and compare SGN with a uniform-panel Nystr\"om method (UM). In UM, each edge is discretized by uniform panels with the same Gauss--Legendre degree $p$. Figure~\ref{fig:experiment2_wholefield} shows the corresponding scattering fields for the droplet, triangle, and L-shaped domains, which provides a qualitative figure of the accuracy gains obtained by singularity guided panelization. Table~\ref{tab:experiment2} summarizes the number of degrees of freedom (DoFs) and the relative error $E$ for a fixed target point. Across all shapes, SGN achieves substantially higher accuracy than UM. In particular, for the droplet and triangle, SGN reaches near machine precision, while UM stagnates at $10^{-8}$--$10^{-6}$. For the re-entrant L-shaped domain, SGN still improves the accuracy by several orders, but the error is noticeably higher. We call this the re-entrant angle effect. The detailed reason for this accuracy loss will be studied in the next example.

\begin{table}[!ht]
    \centering
    \caption{Baseline comparison between UM and SGN.}
    \label{tab:experiment2}
    \begin{tabular}{@{}l l c c c@{}}
    \toprule
    \textbf{Domain} & \textbf{Method} & \textbf{$\varepsilon_{\rm pan}$} & \textbf{DoFs} & \textbf{Relative error} \\ \midrule

    \multirow{2}{*}[-0.5ex]{Droplet}       
        & UM  & -- & 320 &  5.2073e-08 \\ 
        & SGN  & $10^{-9}$ & 832 & 9.3983e-16 \\ \midrule \addlinespace 
    
    \multirow{2}{*}[-0.5ex]{Triangle}      
        & UM  & -- & 320 & 4.0626e-06 \\
        & SGN  & $10^{-7}$ & 2208 & 8.1924e-16 \\ \midrule \addlinespace 
    
    \multirow{2}{*}[-0.5ex]{L-shaped domain} 
        & UM  & -- & 320 & 2.1345e-04 \\
        & SGN  & $10^{-12}$ & 6656 &  8.0433e-10  \\ 
    \bottomrule
    \end{tabular}
\end{table}

\subsubsection{Re-entrant angle effect}
To rigorously isolate the role of the corner geometry, we consider a family of single-corner droplet domains. Figure~\ref{fig:exp1-error-threshold}(b) shows $\log_{10}E$ versus the threshold exponent $N$ for several opening angles $\theta$. For any fixed $N$ in the refinement-dominated regime, the error increases monotonically as $\theta$ moves from convex to re-entrant. Specifically, while the error remains near machine precision for acute and right angles, it rises sharply as $\theta$ approaches $2\pi$. This degradation occurs despite the SGN algorithm automatically clustering panels near the singularity, suggesting an intrinsic difficulty in the re-entrant regime that simple mesh refinement cannot easily overcome.

To investigate the underlying cause, we select two representative cases for a detailed comparison: a convex corner with $\theta=\pi/2$ and a re-entrant corner with $\theta=3\pi/2$. Both runs use SGN with the four-branch singular sets and an adaptive $\varepsilon_{\rm pan}$. We also list the dominant corner exponents to explore the differences, as shown in Table~\ref{tab:experiment3}. The results show that the re-entrant case requires markedly more panels near the corner and still exhibits a larger error. A possible explanation is that the performance gap is driven by complex roots: in the re-entrant case the characteristic set contains $z=a\pm ib$, which may induce logarithmic oscillations and make geometric clustering less effective. However, we observe that the same pattern also appears for Laplace problems, where all singular exponents are real, suggesting that complex exponents are not the primary cause.

\begin{figure}[!ht]
  \centering
  \begin{minipage}[b]{.46\linewidth}
    \centering
    \includegraphics[width=\linewidth]{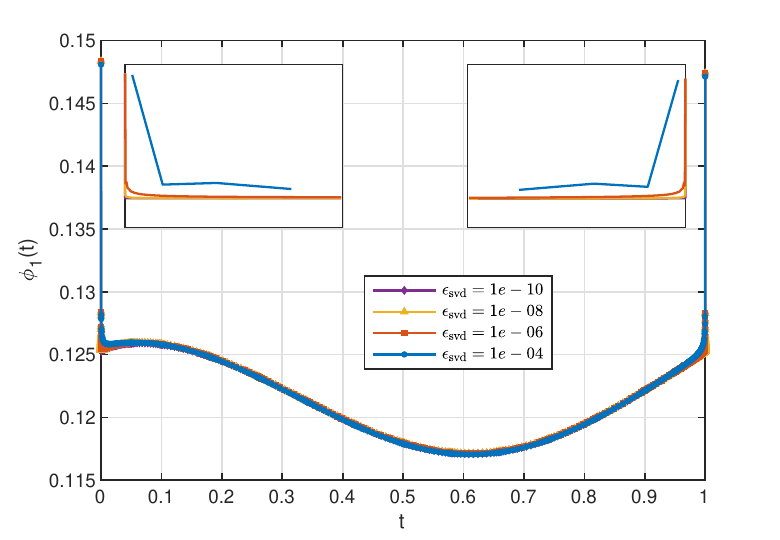}\\
    \footnotesize (a) Density function $\phi_1$ at angle $\theta=\pi/2$
  \end{minipage}\hfill
  \begin{minipage}[b]{.46\linewidth}
    \centering
    \includegraphics[width=\linewidth]{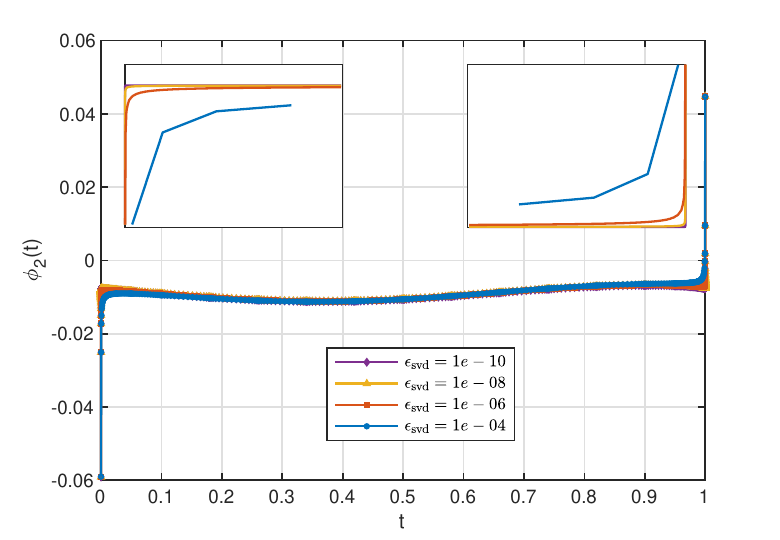}\\
    \footnotesize (b) Density function $\phi_2$ at angle $\theta=\pi/2$
  \end{minipage}
  \caption{Near-corner behavior at opening angle $\theta=\pi/2$. (a) and (b) show the behavior of computed $\phi_1$ and $\phi_2$, respectively. In each panel, the small insets are zoomed views of the innermost boundary panels adjacent to the corner (near $t=0$ and $t=1$).}
  \label{fig:exp2-phi12}
\end{figure}

\begin{table}[!ht]
\centering
\caption{Angle isolation on single-corner droplets.}
\label{tab:experiment3}
\begin{tabular}{@{}l c c c c@{}}
\toprule
\textbf{Domain} & \textbf{Representative exponents} & \textbf{$\varepsilon_{\rm pan}$} & \textbf{DoFs} & \textbf{$E$} \\
\midrule
$\theta=\pi/2$ & $ z_1 \approx 0.5445$, $ z_2 \approx 0.9085$ & $10^{-9}$
& 832 & 9.3983e-16 \\
$\theta=3\pi/2$ 
& $ z_1 \approx 2.7396 \pm 1.1190i,  z_2 \approx 1$ & $10^{-15}$ & 1536 & 4.9599e-12 \\
\bottomrule
\end{tabular}
\end{table}
To visualize what happens near the corner, we plot for $\theta=\pi/2$ the two density components in Figure~\ref{fig:exp2-phi12}. Even in the convex case, the innermost panels show mild numerical artifacts—slight flattening and tiny ripples—which are consistent with parameter-space crowding. Locally, a wedge of opening $\theta$ is conformally flattened by the power map $w=z^{\pi/\theta}$, whose Jacobian scales like $|dw/dz|\sim r^{\,\pi/\theta-1}$. When $\theta>\pi$, the Jacobian blows up as $r\to0$, forcing many parameter nodes to squeeze into a very short arc in physical space. This amplifies near singular quadrature error and the condition number of the discretized system, forcing the error to plateau earlier than expected given the refined clustering.

% --------------------------
% Section 6: Conclusions
% --------------------------
\section{Conclusions and future work}\label{sec6}
This paper investigates the two-dimensional Lam\'e system on domains with corners through boundary integral equations. On the analytical side, a local Mellin analysis on a wedge yields a factorizable characteristic equation for the boundary density's singular exponents, cleanly separating physical branches and explaining which exponents are compatible with the corresponding problems. In weighted Sobolev spaces $H_\nu^s(\Gamma)^2$, we derive a computable Fredholm criterion for the interior double-layer formulation by excluding the real parts of these exponents and the integer poles. Moreover, we prove the invertibility of an explicit density-to-Taylor mapping that links singular density to smooth boundary data. On the numerical side, we design a singularity guided high-order Nystr\"om (SGN) scheme: corner exponents drive panel adaptivity through a multi-exponent Legendre-tail indicator, and near/self interactions are evaluated with a Cauchy-Legendre three-term recurrence. We prove an a priori error bound that couples the panel indicator with exponentially small far-field quadrature error. Numerical experiments on cornered geometries show that SGN achieves higher accuracy than a uniform Nystr\"om baseline, identify re-entrant corners as crowding-limited, and demonstrate that hybrid singular-basis expansions are severely ill-conditioned.

Based on the present results, we propose several directions for future research. We intend to apply the SGN framework to time-harmonic and quasi-static elasticity to quantify corner behaviors, as well as to transmission problems involving interfacial cracks. Geometrically, the method shows promise for generalization to three-dimensional axisymmetric models and polyhedral structures, where edge and vertex singularities interact. Finally, we will investigate shape derivatives in cornered domains for optimal design, and incorporate the close evaluation~\cite{baohualai25} to ensure accuracy near boundaries.

\backmatter
% -------------------------------
% Statements and Declarations
% -------------------------------

% If you prefer the CRediT taxonomy wording, you may alternatively use:
% Conceptualization — First Author; Methodology — First Author (analysis) and Second Author (algorithmic design); 
% Software — First Author; Validation — First Author; Formal analysis — First Author; Investigation — First Author; 
% Resources — Second Author; Data curation — First Author; Writing—original draft — First Author; 
% Writing—review \& editing — Both authors; Visualization — First Author; Supervision — Second Author; 
% Project administration — Second Author; Funding acquisition — None.

% --------------------------
% Appendices
% --------------------------
\begin{appendices}
\section{Explicit formulas for the matrices \texorpdfstring{$\mathbf{A}(z,\theta)$}{A(z,theta)} and \texorpdfstring{$\mathbf{F}(k,z,\theta)$}{F(k,z,theta)}}
\label{app:AF_calculation}
    This appendix lists the explicit formulas for the matrix $\mathbf{A}(z,\theta)$ and $\mathbf{F}(k,z,\theta)$ that appear in Eq.~\eqref{eq:MapEquation}.
    
\subsection*{A.1. Columns of $\mathbf{A}(z,\theta)$}
    The column vectors of $\mathbf{A}(z,\theta)$ are
    {\footnotesize
    \begin{align*}
    \mathbf{A}_1(z,\theta) =
    \begin{bmatrix}
    -\frac{1}{2} \\[3pt]
    -c_1 \cot(\pi z)  \\[3pt]
    -c_1 \csc(\pi z) \sin(z (\pi - \theta)) + \frac{1}{4} c_2 \csc(\pi z) \left[-(z+2) \sin(z (\pi - \theta)) + z \sin(\pi z-(z-2)\theta)\right] \\[3pt]
    -c_1 \cos(z (\pi - \theta)) \csc(\pi z) + \frac{1}{2} c_2 z \csc(\pi z) \sin(\theta) \sin(\pi z - (z-1) \theta)
    \end{bmatrix},
    \end{align*}
    \begin{align*}
    \mathbf{A}_2(z,\theta) =
    \begin{bmatrix}
    c_1 \cot(\pi z)\\[3pt]
    -\frac{1}{2}\\[3pt]
    c_1 \cos(z (\pi - \theta)) \csc(\pi z) + \frac{1}{2} c_2 z \csc(\pi z) \sin(\theta) \sin(\pi z - (z-1) \theta)\\[3pt]
    -c_1 \csc(\pi z) \sin(z (\pi - \theta)) + \frac{1}{4} c_2 \csc(\pi z) \left[(z-2) \sin(z (\pi - \theta)) - z \sin(\pi z - (z-2) \theta)\right]
    \end{bmatrix},
    \end{align*}
    \begin{align*}
    \mathbf{A}_3(z,\theta) =
    \begin{bmatrix}
    -c_1 \csc(\pi z) \sin(z (\pi - \theta)) + \frac{1}{4} c_2 \csc(\pi z) \left[(z-2) \sin(z (\pi - \theta)) - z \sin(\pi z - (z+2) \theta)\right]\\[3pt]
    c_1 \cos(z (\pi - \theta)) \csc(\pi z) - \frac{1}{2} c_2 z \csc(\pi z) \sin(\theta) \sin(\pi z - (z+1) \theta)\\[3pt]
    -\frac{1}{2}\\[3pt]
    c_1 \cot(\pi z)
    \end{bmatrix},
    \end{align*}
    \begin{align*}
    \mathbf{A}_4(z,\theta) =
    \begin{bmatrix}
    -c_1 \cos(z (\pi - \theta)) \csc(\pi z) - \frac{1}{2} c_2 z \csc(\pi z) \sin(\theta) \sin(\pi z - (z+1) \theta) \\[3pt]
    -c_1 \csc(\pi z) \sin(z (\pi - \theta)) + \frac{1}{4} c_2 \csc(\pi z) \left[-(z+2) \sin(z (\pi - \theta)) + z \sin(\pi z - (z+2)\theta)\right] \\[3pt]
    -c_1 \cot(\pi z) \\[3pt]
    -\frac{1}{2}
    \end{bmatrix}.
    \end{align*}
    }
        
\subsection*{A.2. Coefficients of $\mathbf{F}(k,z,\theta)$}
    For $k\ge4$, the four column vectors are
    {\footnotesize
    \begin{align*}
    \mathbf{F}_1(k,z,\theta) &=
    \begin{bmatrix}
    0\\[3pt]
    -\dfrac{c_1}{\pi(k-z)}\\[3pt]
    \dfrac{c_2 k\sin((-2+k)\theta) + \big[4c_1 - c_2(2+k)\big]\sin(k\theta)}{4\pi(k - z)}\\[3pt]
    -\dfrac{2c_1\cos(k\theta) + c_2 k\sin(\theta)\sin((1-k)\theta)}{2\pi(k - z)}
    \end{bmatrix},\\
    \mathbf{F}_2(k,z,\theta) &=
    \begin{bmatrix}
    \dfrac{c_1}{\pi(k-z)}\\[3pt]
    0\\[3pt]
    \dfrac{2c_1\cos(k\theta) - c_2 k\sin(\theta)\sin((1-k)\theta)}{2\pi(k- z)}\\[3pt]
    \dfrac{-c_2 k\sin((-2 + k)\theta) + \big[4c_1 + c_2(-2 + k)\big]\sin(k\theta)}{4\pi(k - z)}
    \end{bmatrix},\\
    \mathbf{F}_3(k,z,\theta) &=
    \begin{bmatrix}
    \dfrac{\big[4c_1 + c_2(-2 + k)\big]\sin(k\theta) - c_2 k\sin((2 + k)\theta)}{4\pi(k - z)}\\[3pt]
    \dfrac{2c_1\cos(k\theta) - c_2 k\sin(\theta)\sin((1 + k)\theta)}{2\pi(k-z)}\\[3pt]
    0\\[3pt]
    \dfrac{c_1}{\pi(k-z)}
    \end{bmatrix},\\
    \mathbf{F}_4(k,z,\theta) &=
    \begin{bmatrix}
    -\dfrac{2c_1\cos(k\theta) + c_2 k\sin(\theta)\sin((1 + k)\theta)}{2\pi(k-z)}\\[3pt]
    \dfrac{\big[4c_1 - c_2(2 + k)\big]\sin(k\theta) + c_2 k\sin((2 + k)\theta)}{4\pi(k - z)}\\[3pt]
    -\dfrac{c_1}{\pi(k - z)}\\[3pt]
    0
    \end{bmatrix}.
    \end{align*}
    }

    \noindent For the low-order cases, we list the explicit block for $k=0,1,2$ and the four columns for $k=3$:
    {\footnotesize
    \begin{align*}
         \mathbf{F}(0,z,\theta) =
         \begin{bmatrix}
        0 & -\frac{c_1}{\pi z} & 0 & \frac{c_1}{\pi z} \\[3pt]
        \frac{c_1}{\pi z} & 0 & -\frac{c_1}{\pi z} & 0 \\[3pt]
        0 & -\frac{c_1}{\pi z} & 0 & \frac{c_1}{\pi z} \\[3pt]
        \frac{c_1}{\pi z} & 0 & -\frac{c_1}{\pi z} & 0
        \end{bmatrix}.
    \end{align*}
    \begin{align*}
        \mathbf{F}(1,z,\theta) =\begin{bmatrix}
        0 & \frac{c_1}{\pi(1-z)} & \frac{c_1 \sin(\theta)}{\pi(1-z)} & \frac{c_1 \cos(\theta)}{\pi (-1 + z)} \\[3pt]
        \frac{c_1}{\pi (-1 + z)} & 0 & \frac{c_1 \cos(\theta)}{\pi(1-z)} & \frac{c_1 \sin(\theta)}{\pi (1 - z)} \\[3pt]
        \frac{c_1 \sin(\theta)}{\pi (1-z)} & \frac{c_1 \cos(\theta)}{\pi (1-z)} & 0 & \frac{c_1}{\pi (-1 + z)} \\[3pt]
        \frac{c_1 \cos(\theta)}{\pi (-1 + z)} & \frac{c_1 \sin(\theta)}{\pi(1-z)} & \frac{c_1}{\pi (1-z)} & 0
        \end{bmatrix}.
    \end{align*}
    \begin{align*}
        \mathbf{F}(2,z,\theta) =\begin{bmatrix}
        0& -\frac{c_1}{\pi(-2+z)} & \frac{(-c_1 + c_2 + c_2 \cos(2\theta)) \sin(2\theta)}{\pi (-2 + z)} & \frac{c_1 \cos(2\theta) + c_2 \sin^2(2\theta)}{\pi (-2 + z)} \\[3pt]
        \frac{c_1}{\pi (-2 + z)} & 0 & \frac{-c_1 \cos(2\theta)+c_2\sin^2(2\theta)}{\pi (-2+z)} &
        \frac{(-c_1 + c_2 - c_2 \cos(2\theta)) \sin(2\theta)}{\pi (-2 + z)} \\[3pt]
        \frac{(-c_1 + 2c_2) \sin(2\theta)}{\pi (-2 + z)} & -\frac{c_1 \cos(2\theta)}{\pi (-2+z)} & 0 & \frac{c_1}{\pi (-2+z)} \\[3pt]
        \frac{c_1 \cos(2\theta)}{\pi (-2 + z)} & -\frac{c_1 \sin(2\theta)}{\pi (-2 + z)} & -\frac{c_1}{\pi (-2 + z)} & 0
        \end{bmatrix}.
    \end{align*}
    \begin{align*}
        \mathbf{F}_1(3,z,\theta) &= 
        \begin{bmatrix}
        0 \\[3pt]
        \frac{c_1}{\pi (-3 + z)} \\[3pt]
        \frac{-c_2 \sin(\theta) + (-c_1 + c_2) \sin(3 \theta)}{\pi (-3 + z)} \\[3pt]
        \frac{-c_2 \cos(\theta) + (c_1 + c_2) \cos(3 \theta)}{\pi (-3 + z)}
        \end{bmatrix},
        \mathbf{F}_2(3,z,\theta) =
        \begin{bmatrix}
        -\frac{c_1}{\pi(-3+z)} \\[3pt]
        0 \\[3pt]
        \frac{-c_2 \cos(\theta) + (-c_1 + c_2)\cos(3\theta)}{\pi (-3+ z)} \\[3pt]
        -\frac{c_1 \sin(3 \theta)}{\pi (-3 + z)}
        \end{bmatrix},\\
        \mathbf{F}_3(3,z,\theta)&=
        \begin{bmatrix}
        \frac{-4 c_2\sin(\theta)-(4 c_1 + c_2) \sin(3 \theta) + 3 c_2 \sin(5 \theta)}{4 \pi (-3 + z)} \\[3pt]
        \frac{-2 c_1 \cos(3 \theta) + 3 c_2 \sin(\theta) \sin(4 \theta)}{2 \pi(-3+z)} \\[3pt]
        0\\[3pt]
        -\frac{c_1}{\pi (-3+z)}
        \end{bmatrix},
        \mathbf{F}_4(3,z,\theta) =
        \begin{bmatrix}
        \frac{2 c_1 \cos(3 \theta) + 3 c_2 \sin(\theta) \sin(4 \theta)}{2 \pi (-3 + z)} \\[3pt]
        \frac{(-4 c_1 + 5 c_2) \sin(3 \theta) - 3 c_2 \sin(5 \theta)}{4 \pi (-3 + z)} \\[3pt]
        \frac{c_1}{\pi (-3 + z)} \\[3pt]
        0
        \end{bmatrix}.
    \end{align*}
    }
    
\section{Invertibility of \texorpdfstring{$\mathbf{B}(\theta)$}{B(theta)} at the \texorpdfstring{$\theta=\pi$}{theta=pi} limit}
\label{app:invertibility}
    This appendix proves that the $4N\times4N$ mapping matrix $\mathbf{B}(\theta)$ is nonsingular at $\theta=\pi$ and gives an explicit blockwise formula for $\mathbf{B}(\pi)$. Throughout, we fix an ordered list of integer exponents $m_n:=n-1$ for $n=1,\ldots,N$, so that Taylor rows are $k=0,1,\ldots,N-1$ and the $n$-th column block resonates exactly at the row $k=m_n$.
    
\subsection*{B.1 Closed forms of $\mathbf{A}(z,\pi)$ and $\mathbf{F}(k,z,\pi)$}
    % Remark 1：这里 $\theta=\pi$ 的讨论思路跟一般情形不太一样，需要从极限的角度来理解：A 在 $z=m$ 处是一个亚纯矩阵函数，有洛朗展开（奇异主部矩阵+有限部分+极限部分）。奇异主部矩阵完全决定了允许的极限核子空间（否则等式无法成立），有限项并不会改变这个核子空间，只会影响后续我们在"共振行"上得到的有限主部系数～
    The singular exponents $z$ satisfy
    $$
    \Big(z^2\sin^2\theta-\sin^2\big((2\pi-\theta)z\big)\Big)
    \Big(z^2\sin^2\theta-c^2\sin^2(\theta z)\Big)=0.
    $$
    At $\theta=\pi$ this reduces to $\sin(\pi z)=0$, hence $z=m\in\mathbb{Z}$. Substituting $\theta=\pi$ into the explicit columns $\mathbf{A}_j(z,\theta)$ gives
    $$
    \mathbf{A}(z,\pi)=
    \begin{bmatrix}
    -\tfrac12 & c_1\cot(\pi z) & 0 & -c_1\csc(\pi z)\\
    -c_1\cot(\pi z) & -\tfrac12 & c_1\csc(\pi z) & 0\\
    0 & c_1\csc(\pi z) & -\tfrac12 & -c_1\cot(\pi z)\\
    -c_1\csc(\pi z) & 0 & c_1\cot(\pi z) & -\tfrac12
    \end{bmatrix}.
    $$
    Let $s_m:=(-1)^m$. Near $z=m$, it holds
    $$
    \cot(\pi z)=\frac{1}{\pi(z-m)}+O(z-m),\qquad
    \csc(\pi z)=\frac{s_m}{\pi(z-m)}+O(z-m),
    $$
    so $\mathbf{A}(z,\pi)$ is meromorphic with a simple pole $z=m$. Its Laurent expansion is
    \begin{equation}\label{eq:Laurent-A}
    \mathbf{A}(z,\pi)\;=\;\frac{1}{\pi(z-m)}\,\mathbf{S}_m\;+\;\mathbf{A}^{(0)}(m)\;+\;O(z-m),
    \end{equation}
    with the singular principal-part matrix
    $$
    \mathbf{S}_m=
    \begin{bmatrix}
    0 & c_1 & 0 & -c_1 s_m\\
    -c_1 & 0 & c_1 s_m & 0\\
    0 & c_1 s_m & 0 & -c_1\\
    -c_1 s_m & 0 & c_1 & 0
    \end{bmatrix}.
    $$
    If $\mathbf{p}(z)$ is an analytic branch with $\mathbf{A}(z,\pi)\mathbf{p}(z)=\mathbf{0}$, inserting \eqref{eq:Laurent-A} and taking $z\to m$ forces the leading equation $\mathbf{S}_m\,\mathbf{p}(m)=\mathbf{0}$. Thus the limiting kernel at $z=m$ equals $\ker \mathbf{S}_m$, which is two-dimensional with basis
    \begin{equation}\label{eq:kernel-basis}
    \mathbf{p}^{(+)}(m)=\begin{bmatrix}1\\ 0\\ s_m\\ 0\end{bmatrix},\qquad
    \mathbf{p}^{(-)}(m)=\begin{bmatrix}0\\ 1\\ 0\\ s_m\end{bmatrix}.
    \end{equation}
    The finite part $\mathbf{A}^{(0)}(m)=-\tfrac12\mathbf{I}_4$ does not change the kernel directions selected by $\mathbf{S}_m$, but it contributes the jump term $-\tfrac12\,\mathbf{p}(m)$ to the resonant Taylor coefficient in the mapping. 
    
    From the closed formulas of $\mathbf{F}(k,z,\theta)$, one obtains that for all $k\ge0$,
    \begin{equation}\label{eq:F-pi}
    \mathbf{F}(k,z,\pi)=\frac{c_1}{\pi(k-z)}\,\mathbf{M}(k),\quad
    \mathbf{M}(k)=
    \begin{bmatrix}
    0 & 1 & 0 & -(-1)^k\\
    -1 & 0 & (-1)^k & 0\\
    0 & (-1)^k & 0 & -1\\
    -(-1)^k & 0 & 1 & 0
    \end{bmatrix}.
    \end{equation} 
    Acting on \eqref{eq:kernel-basis} at $m_n$ and writing $s_n:=(-1)^{m_n}$, it yields
    \begin{align*}
    \mathbf{M}(k)\,\mathbf{p}^{(+)}(m_n)=
    \begin{bmatrix}0\\ -1+s_n(-1)^k\\ 0\\ s_n-(-1)^k\end{bmatrix},\quad
    \mathbf{M}(k)\,\mathbf{p}^{(-)}(m_n)=
    \begin{bmatrix}1-s_n(-1)^k\\ 0\\ (-1)^k-s_n\\ 0\end{bmatrix}.
    \end{align*}
    Hence, these two columns vanish iff $k$ and $m_n$ have the same parity.
    
\subsection*{B.2 Explicit blocks of $\mathbf{B}_{k,n}(\pi)$ via pole analysis}
    For a single branch $(z,\mathbf{p})$, the mapping in Sec.~\ref{sec4} gives
    $$
    \tilde{\boldsymbol{h}}(t)=\mathbf{A}(z,\pi)\,\mathbf{p}\,t^z
    \;+\;
    \sum_{k=0}^{\infty}\mathbf{F}(k,z,\pi)\,\mathbf{p}\,t^{k}.
    $$
    \emph{(i) Non-resonant rows $k\neq z$.} By Eq.~\eqref{eq:F-pi}, the coefficient of $t^k$ is
    $$
    \mathbf{F}(k,z,\pi)\,\mathbf{p}
    =\frac{c_1}{\pi(k-z)}\,\mathbf{M}(k)\,\mathbf{p}.
    $$
    \emph{(ii) Resonant row $k=z$.} Write $z = m + \varepsilon$, $\varepsilon \to 0$, and choose an analytic branch
    $$\mathbf{p}(z)=\mathbf{p}(m)+\varepsilon\,\mathbf{s}+O(\varepsilon^2)$$ with
    $\mathbf{p}(m)\in\ker\mathbf{A}(m,\pi)$ and
    $\mathbf{s}=\lim_{z\to m}(\mathbf{p}(z)-\mathbf{p}(m))/(z-m)$. Inserting \eqref{eq:Laurent-A} into $\mathbf{A}(z,\pi)\mathbf{p}(z)\equiv 0$ and matching orders in $(z-m)$ yields the differentiated solvability condition
    $$
    \frac{1}{\pi}\,\mathbf{S}_m\,\mathbf{s}\;-\;\frac12\,\mathbf{p}(m)\;=\;0 \quad\Longleftrightarrow\quad
    \mathbf{S}_m\,\mathbf{s}\;=\;\frac{\pi}{2}\,\mathbf{p}(m).
    $$
    Using this and Eq.\eqref{eq:F-pi}, the resonant $t^m$-coefficient becomes
    \begin{equation}\label{eq:res-row}
    -\frac12\,\mathbf{p}(m)\;+\;\frac{1}{\pi}\mathbf{S}_m\mathbf{s}\;-\;\frac{c_1}{\pi}\,\mathbf{M}(m)\,\mathbf{s}
    \;=\;-\frac{c_1}{\pi}\,\mathbf{M}(m)\,\mathbf{s}.
    \end{equation}
    
    Now fix the ordered integers $\{m_n\}_{n=1}^{N}$, and for each $n$ pick four analytic coefficient vectors $\mathbf{p}^{n,j}(z,\theta)$, $j=1,\ldots,4$, with $z_{n,j}(\theta)\to m_n$ as $\theta\to\pi$. At $\theta=\pi$, define
    $$
    \mathbf{p}^{n,j}:=\mathbf{p}^{n,j}(m_n,\pi),\qquad
    \mathbf{s}^{\,n,j}:=\lim_{\theta\to\pi}\frac{\mathbf{p}^{n,j}(z_{n,j}(\theta),\theta)-\mathbf{p}^{n,j}}{\,z_{n,j}(\theta)-m_n\,}.
    $$
    Collect them in the $4\times4$ matrices
    $\mathcal{C}_n=[\mathbf{p}^{n,1}\ \mathbf{p}^{n,2}\ \mathbf{p}^{n,3}\ \mathbf{p}^{n,4}]$,
    $\mathcal{S}_n=[\mathbf{s}^{\,n,1}\ \mathbf{s}^{\,n,2}\ \mathbf{s}^{\,n,3}\ \mathbf{s}^{\,n,4}]$.
    Summing branchwise contributions, the $4\times4$ blocks of $\mathbf{B}(\pi)$ are:
    \begin{equation}\label{eq:Bknpi}
    \mathbf{B}_{k,n}(\pi)=
    \begin{cases}
    \dfrac{c_1}{\pi(k-m_n)}\,\mathbf{M}(k)\,\mathcal{C}_n, & k\neq m_n,\\[8pt]
    -\dfrac{c_1}{\pi}\,\mathbf{M}(m_n)\,\mathcal{S}_n, & k=m_n,
    \end{cases}
    \end{equation}
    with $\mathbf{M}(k)$ given in \eqref{eq:F-pi}. In particular, we may take
    $$
    \mathbf{p}^{n,1}=\mathbf{p}^{(+)}(m_n),\qquad
    \mathbf{p}^{n,2}=\mathbf{p}^{(-)}(m_n),
    $$
    and choose $\{\mathbf{p}^{n,3},\mathbf{p}^{n,4}\}$ so that $\mathcal{C}_n$ is invertible.
    
\subsection*{B.3 Parity reduction and invertibility}
    %% 在这个章节中，考虑把正交变换P显式作用到每个4*4小块上面，变为偶/奇组合（正如kirill serkh中的文章中那样，对称/反对称组合）。直观理解的话，就是把"两条边"的变量换成"偶/奇"两个独立通道；在 $\theta=\pi$ 的对称几何里，这是最自然的分解。
    % 1. 对称化：使用 P 将系统拆成两个互不耦合的扇区；
    % 2. 下三角：由核-核外的幂级数结构可知，当行按 k=0,1,... 递增，列按 m_1,...递增排序后，每个扇区都是块下三角。其对角块正是各自的共振 2*2 区域；
    % 3. 非零对角：把每个对角写成显式的矩阵；
    Introduce the orthogonal parity map
    $$
    P_4:=\frac{1}{\sqrt2}\begin{bmatrix}
    1&0&1&0\\ 0&1&0&1\\ 1&0&-1&0\\ 0&1&0&-1
    \end{bmatrix},\qquad \mathcal P:=\mathbf I_N\otimes P_4.
    $$
    Let $\mathcal Q$ be the block-diagonal permutation which, for each odd row-index $k$ (respectively, odd column resonance index $m_n$), swaps the two $2$-dimensional parity subspaces so that the active $2\times2$ subblock of $P_4^{\top}\mathbf M(k)P_4$ always occupies the same position. Using
    \begin{equation}\label{eq:Mhat}
    P_4^{\top}\mathbf M(k)P_4=
    \begin{bmatrix}
    0&0&0&1+(-1)^k\\
    0&0&-(1+(-1)^k)&0\\
    0&1-(-1)^k&0&0\\
    (-1)^k-1&0&0&0
    \end{bmatrix},
    \end{equation}
    we obtain the orthogonal similarity
    $$
    \widehat{\mathbf B}(\pi):=(\mathcal P\mathcal Q)^{\top}\,\mathbf B(\pi)\,(\mathcal P\mathcal Q)
    =\begin{bmatrix}\mathbf B_{+}&\mathbf 0\\ \mathbf 0&\mathbf B_{-}\end{bmatrix},
    $$
    i.e., the matrix decouples into two independent parity sectors. After the parity reordering $\mathcal Q$, the active $2\times2$ block of $P_4^{\top}\mathbf M(k)P_4$ equals $\sigma_k\,\mathbf L$ with $\sigma_k\in\{\pm2\}$ depending only on the parity of $k$. Similarly, the resonant selector equals $\tau_{m_n}\,\mathbf L$ with $\tau_{m_n}\in\{\pm2\}$ depending only on the parity of $m_n$.

    Fix one sector and move the resonant rows $k=m_n$ to the bottom by a permutation. For every non-resonant row $k\neq m_n$, by \eqref{eq:Bknpi} and \eqref{eq:Mhat} the corresponding $2\times2$ block factors as
    $$
    \widehat{\mathbf B}_{k,n}(\pi)
    =\frac{c_1}{\pi(k-m_n)}\,\sigma_k\,\mathbf L\;\widehat{\mathcal C}_n^{(\mathrm{act})},
    $$
    where $\widehat{\mathcal C}_n^{(\mathrm{act})}$ is the $2\times2$ restriction of $P_4^{\top}\mathcal C_n P_4$ to the active parity subspace and is invertible. For each $n$, the resonant block at $k=m_n$ reduces to
    $$
    \mathbf R_n \;=\; -\frac{c_1}{\pi}\,\tau_{m_n}\,\mathbf L\;\widehat{\mathcal S}_n^{(\mathrm{act})},
    $$
    where $\widehat{\mathcal S}_n^{(\mathrm{act})}$ is the $2\times2$ restriction of $P_4^{\top}\mathcal S_n P_4$ to the active subspace. By choosing analytic kernel branches, we may assume $\widehat{\mathcal S}_n^{(\mathrm{act})}$ invertible. Under the canonical normalization $\widehat{\mathcal S}_n^{(\mathrm{act})}=\mathbf I_2$, one has $\mathbf R_n=\pm(c_1/\pi)\mathbf L$.

    Because in this sector the set of row indices coincides with the set of resonance indices of the same parity, there is exactly one invertible pivot  $\mathbf R_n$ in each column block. Perform block Gaussian elimination with these $2\times2$ pivots:
    $$(\text{row }k)\leftarrow(\text{row }k)-\widehat{\mathbf B}_{k,n}(\pi)\,\mathbf R_n^{-1}\cdot(\text{row }m_n),\qquad k\neq m_n.$$
    Each step is a unit lower-triangular left multiplication and preserves invertibility. Equivalently, this eliminates non-resonant couplings via the Schur complement of the resonant diagonal, yielding a block lower-triangular matrix whose bottom-right diagonal is $\mathrm{diag}(\mathbf R_n)$ and whose top-left block is, by \eqref{eq:Bknpi}–\eqref{eq:Mhat}, a Cauchy-type matrix with kernel entries $1/(k-m_n)$ up to invertible block-diagonal factors. Hence it has full rank. As each $\mathbf R_n$ is invertible, the whole sector is invertible. Orthogonality of $(\mathcal P\mathcal Q)$ gives $\det\mathbf B(\pi)=\det(\mathbf B^+)\cdot\det(\mathbf B^-)\neq0$.
\end{appendices}

%%===========================================================================================%%
%% If you are submitting to one of the Nature Portfolio journals, using the eJP submission   %%
%% system, please include the references within the manuscript file itself. You may do this  %%
%% by copying the reference list from your .bbl file, paste it into the main manuscript .tex %%
%% file, and delete the associated \verb+\bibliography+ commands.                            %%
%%===========================================================================================%%

\bibliography{sn-bibliography}% common bib file

%% if required, the content of .bbl file can be included here once bbl is generated
%%\input sn-article.bbl

\end{document}